%% file: FFDIC_v3_submission.tex
\documentclass{amsart}
\usepackage[utf8]{inputenc}
\usepackage{amsmath}
\usepackage{palatino}
\usepackage{mathpazo}
\usepackage{float}
\usepackage[shortlabels]{enumitem}
\usepackage[makeroom]{cancel}
\usepackage{tikz-cd} 
\usepackage{quiver}
\usepackage[nopatch=footnote]{microtype}
\usepackage{amsfonts, amsthm}
\usepackage{subfig}
\usepackage{mathtools}
\usepackage{amssymb, hyperref, color, graphicx, caption}
\usepackage[top=2cm,bottom=3cm,textwidth=460pt]{geometry}
\usepackage{fullpage}
\usepackage{setspace}
\usepackage{tikzsymbols}
\usepackage{braket}
\usepackage{tcolorbox}
\usepackage{xcolor}   
\hypersetup{
    colorlinks=true, 
    linktoc=all,     
    linkcolor=blue,  
    citecolor=blue
}
\usepackage{calligra,mathrsfs}

\usepackage{subfiles}
\usepackage{verbatim}

\makeatletter
\newcommand{\appendixref}[2]{%
  \expandafter\gdef\csname appref@#1\endcsname{#2}%
}

\newcommand{\smartthmref}[1]{%
  Theorem~\ref{#1}%
  \@ifundefined{appref@#1}{}{/\csname appref@#1\endcsname}%
}
\makeatother
\usepackage{booktabs}
\usepackage{longtable}
\usepackage{array}
\numberwithin{equation}{section}
\newtheorem{theorem}{Theorem}[section]
\newtheorem{remark}{Remark}[section]
\newtheorem{definition}[theorem]{Definition}
\newtheorem{lemma}[theorem]{Lemma}

\newtheorem{corollary}[theorem]{Corollary}
\newtheorem{notation}[theorem]{Notation}
\newtheorem{prop}[theorem]{Proposition}

\newtheorem{example}[theorem]{Example}

\newcommand{\Hom}{\mathrm{Hom}}
\newcommand{\End}{\mathrm{End}}
\newcommand{\C}{\mathbb{C}}
\newcommand{\cC}{\mathcal{C}}
\newcommand{\cD}{\mathcal{D}}
\newcommand{\cV}{\mathcal{V}}

\newcommand{\cW}{\mathcal{W}}
\newcommand{\Z}{\mathbb{Z}}

\newcommand{\F}{\mathbb{F}}
\newcommand{\IC}{\mathrm{Ind}(\mathcal{C})}

\newcommand{\Id}{\mathrm{Id}}

\newcommand\ov[1]{\overline{#1}}

\newcommand{\fr}[1]{\mathfrak{#1}}

\newcommand{\cA}{\mathcal{A}
}
\newcommand{\ox}{\otimes}
\newcommand{\pr}{\partial}

\newcommand{\cH}{\mathcal{H}}

\newcommand{\cv}{\mathrm{coev}}
\newcommand{\irep}{\underline{\mathrm{Re}}\mathrm{p}}
\newcommand{\GL}{\mathrm{GL}}
\newcommand{\Or}{\mathrm{O}}
\newcommand{\Sp}{\mathrm{Sp}}
\newcommand{\gul}[1]{\underline{\mathfrak{g}}_{ {\raisebox{1.5pt}{$\scriptstyle  #1$}}}}
\newcommand{\Lul}[2]{\underline{\mathfrak{#1}}_{{\raisebox{1.5pt}{$\scriptstyle  #2$}}}}
\newcommand{\Lulhat}[2]{\widehat{\underline{\mathfrak{#1}}}_{{\raisebox{1.5pt}{$\scriptstyle  #2$}}}}

\newcommand{\gulhat}[1]{\widehat{\underline{\mathfrak{g}}}_{{\raisebox{1.5pt}{$\scriptstyle  #1$}}}}

\newtheorem*{inner}{\innerheader}
\newcommand{\innerheader}{}
\newenvironment{defi}[1]
 {\renewcommand\innerheader{#1}\begin{inner}}
 {\end{inner}}

\newcommand{\ssF}[1]{\mathrm{ss}_{#1}}
 
\newcommand{\ssW}[1]{\mathrm{ss}^{\mathcal{W}}_{#1}}

\newcommand{\ssIW}[1]{\mathrm{I}_{#1}^{\mathcal{W}}}

\newcommand{\ssV}[1]{\mathrm{ss}^{\mathfrak{z}}_{#1}}

\newcommand{\ssIV}[1]{\mathrm{I}_{#1}^{\mathfrak{z}}}

\newcommand{\genheading}[1]{%
  \medskip
  \noindent\textbf{#1.}\par
  \nobreak
  \smallskip
}

\title{Interpolating Feigin--Frenkel Duality at the Critical Level to Matrices of Complex Size}
\author{Andrew Riesen}
\date{\today}

\begin{document}
\begin{abstract}

In this paper, we extend Feigin--Frenkel duality at the critical level to complex rank by identifying two seemingly unrelated constructions in complex rank. On the affine side, we interpolate Molev's construction of higher Segal--Sugawara vectors and thereby describe the centers of universal affine vertex algebras at the critical level in Deligne's interpolating categories. On the $\cW$-side, we construct the classical $\cW$-algebras associated with Feigin's Lie algebras of complex rank $\fr{gl}_{\lambda}$ and $\fr{po}_{\lambda}$ as Poisson vertex algebras, realizing their Drinfeld--Sokolov reduction via an interpolated Adler--Gelfand--Dickey bracket. Upon specialization to positive integer rank in types A, B, and C, this recovers the usual Feigin--Frenkel duality at the critical level. As applications, we obtain a uniform construction of several families of higher Segal--Sugawara vectors for Lie superalgebras and recover a complex-rank analogue of the universal Bethe algebra.
\end{abstract}
\maketitle

\tableofcontents
\section{Introduction}
Feigin--Frenkel duality at the critical level is one of the key structural results underlying modern representation theory. It identifies the center of the universal affine vertex algebra at the critical level with a classical $\cW$-algebra and thereby supplies a key algebraic input to the local geometric Langlands correspondence for loop groups $G((t))$, first developed by Edward Frenkel and Dennis Gaitsgory in \cite{Frenkel_Gatisgory_local_geometric_langlands_loop} (see also \cite{Frenkel_Loop_Group}). The same center, called the Feigin--Frenkel center, also appears naturally in quantum integrable systems, where its elements give rise to the higher Hamiltonians of the Gaudin model and the associated Bethe algebras \cite{Gaudin_FF_center_first_ref}.

More precisely, for a reductive Lie algebra $\fr{g}$ with dual Coxeter number $h^{\vee}$, lacing number $r$, and Langlands dual ${}^L\fr{g}$ with dual Coxeter number ${}^Lh^{\vee}$, Feigin and Frenkel proved in \cite{Feigin_Frenkel_Duality} that the principal affine $W$-algebras satisfy an isomorphism of vertex algebras:
\begin{equation}
    W^k(\fr{g})\cong W^{\ell}({}^L\fr{g})
\end{equation}
where:
\begin{equation}
    r(k+h^{\vee})(\ell+{}^Lh^{\vee})=1
\end{equation}
Passing to the critical limit $k+h^{\vee}\to 0$ yields an isomorphism of Poisson vertex algebras:
\begin{equation}
    \fr{f}_{\fr{g}}:\fr{z}(\widehat{\fr{g}})\xrightarrow{\sim} \cW({}^L\fr{g})
\end{equation}
where $\fr{z}(\widehat{\fr{g}})$ denotes the center of $V^{-h^{\vee}}(\fr{g})$ and the right-hand side is the classical $\cW$-algebra of ${}^L\fr{g}$. In Type $A$, this classical $\cW$-algebra structure arises from the second Adler--Gelfand--Dickey bracket.

The goal of this paper is to transport this picture to complex rank. This is not merely a matter of analytically continuing formulas in the rank parameter. The two sides of Feigin--Frenkel duality naturally live in different complex-rank formalisms: the center at the critical level, or the affine side, is most naturally constructed in Deligne's interpolating categories, while the $\cW$-side is governed by Feigin's Lie algebras of complex rank and their Drinfeld--Sokolov reductions. One of the main points of this paper is that these two a priori different complex-rank worlds fit together via an interpolated Feigin--Frenkel isomorphism at the critical level.

In \cite{Pavel_Complex_Rank_2}, Pavel Etingof proposed studying algebraic structures in interpolating categories under the slogan ``representation theory in complex rank.'' Deligne's interpolating categories have since become a robust source of new phenomena in representation theory and mathematical physics \cite{Interpolated_Yangian, Interpolated_Harish_Chandra,wan2025fockspacetensorproduct,Binder_Rychkov_Deligne_in_interpolating,zeng2025largenvertexalgebras}. The present paper places Feigin--Frenkel duality at the critical level naturally in this framework. The main contributions of this paper are the following.
\begin{enumerate}
    \item We construct the affine side in complex rank: the centers of universal affine vertex algebras at the critical level in Deligne's interpolating categories, together with explicit generating families obtained by interpolating higher Segal--Sugawara vectors.

    \item We construct the corresponding $\cW$-side in complex rank: the classical $\cW$-algebras attached to Feigin's Lie algebras $\fr{gl}_{\lambda}$ and $\fr{po}_{\lambda}$, equipped with Poisson vertex algebra structures defined by an interpolated Adler--Gelfand--Dickey bracket.

    \item We identify the affine and $\cW$ constructions, obtaining complex-rank Feigin--Frenkel isomorphisms at the critical level in Types $A$, $B$, and $C$.
\end{enumerate}

As applications, we construct uniform families of higher Segal--Sugawara vectors for the Lie superalgebras $\widehat{\fr{gl}}_{m|n}$ and $\widehat{\fr{osp}}_{M|2n}$, recovering the constructions of \cite{Molev_Segal_Sugawara_Vectors_super_gl, Molev_Segal_Sugawara_Vectors_super_osp}. We also recover the interpolated universal Bethe algebra of \cite{Rybnikov} directly from the interpolated Feigin--Frenkel center.

\subsection*{The affine side: the center at the critical level in interpolating categories}

Deligne's interpolating categories $\irep(\GL_{\alpha},\F)$, $\irep(\Or_{\alpha},\F)$, and $\irep(\Sp_{\alpha},\F)$ generalize the tensor categories of finite-dimensional representations of their respective classical groups. Concretely, they are rigid symmetric monoidal categories, or more precisely symmetric pseudo-tensor categories, generated by a distinguished object whose diagrammatic endomorphism algebras depend polynomially on the rank $\alpha$. Working internally in these categories produces canonical Lie algebras $\Lul{gl}{\alpha}$, $\Lul{so}{\alpha}$, and $\Lul{sp}{\alpha}$. At positive integer rank, the corresponding specialization functors produce the classical Lie algebras and their representation categories. We refer to these interpolating settings as Types $A$, $B$, and $C$, respectively. Throughout the paper, we use finite-rank to refer to the ordinary classical Lie algebras $\fr{gl}_n$, $\fr{so}_{2n+1}$, and $\fr{sp}_{2n}$, together with their associated affine vertex algebras, centers, and classical $\cW$-algebras, as opposed to
their interpolated counterparts in Deligne's categories or Feigin's complex-rank Lie algebras.

Using the notion of a vertex algebra in a symmetric pseudo-tensor category introduced in \cite{zeng2025largenvertexalgebras}, we define the universal affine vertex algebra $V^k(\gul{\alpha})$ associated with the canonical Lie algebra in the corresponding interpolating category (Theorem \ref{theorem:vertex_algebra_structure_universal_affine}). At the critical level, we define the center as the $\gul{\alpha}[t]$-invariants:
\begin{equation}
    \fr{z}(\gulhat{\alpha}):=\Hom_{\gul{\alpha}[t]}\left(\mathbf{1},V^{-h^{\vee}}\left(\gul{\alpha}\right)\right)
\end{equation}
We prove that this is an ordinary Poisson vertex algebra and that, at positive integer rank, it specializes naturally to the classical Feigin--Frenkel center $\fr{z}(\widehat{\fr{g}}_n)$ (Theorem \ref{theorem:invariant_ind_PVA_structure}). We then interpolate Molev's higher Segal--Sugawara constructions to obtain complete sets of Segal--Sugawara vectors for these centers in Types $A$, $B$, and $C$ (Theorems \ref{theorem:interpolate_segal_sugawara_type_A} and \ref{thm:interpolated_segal_sugawara_Type_B_C}).

A particularly attractive feature of representation theory in complex rank is that it connects naturally with super representation theory. In particular, when $\alpha=m-n$ in Type $A$, and when $\alpha=M-2n$ in Type $B$, the corresponding interpolating categories admit full functors to the representation categories of the supergroups $\GL_{m|n}$ and $\mathrm{OSp}_{M|2n}$, respectively. Consequently, our description of the interpolated Feigin--Frenkel center allows us to recover large families of higher Segal--Sugawara vectors for $\widehat{\fr{gl}}_{m|n}$ and $\widehat{\fr{osp}}_{M|2n}$, and to prove an asymptotic version of Molev and Nurcombe's differential generator conjecture for the even component of the $\widehat{\fr{osp}}$ center (Theorems \ref{theorem:segal_sugawara_super_case}, \ref{theorem:Segal_Sugawara_Super_osp}, \ref{theorem:asymptotic_description}).
This super viewpoint also explains the appearance of interpolated parity functors:
\begin{equation}
    \Pi_{\alpha}:\irep(\GL_{\alpha},\F)\rightarrow \irep(\GL_{-\alpha},\F)^{\mathrm{Tw}}
\end{equation}
\begin{equation}
    \Pi_{\alpha}:\irep(\Sp_{\alpha},\F)\rightarrow \irep(\Or_{-\alpha},\F)^{\mathrm{Tw}}
\end{equation}
These parity dualities organize Molev's higher Segal--Sugawara vectors into a single uniform construction and induce anti-isomorphisms on the interpolated centers (Proposition \ref{prop:interpolated_parity_functor_anti}).

\subsection*{The $\cW$-side: Drinfeld--Sokolov reduction for matrices of complex size}

The $\cW$-side is built from Feigin's Lie algebras of matrices of complex size. In \cite{Feigin_gl_lambda}, Feigin defined $\fr{gl}_{\lambda}$ as a quotient of $U(\fr{sl}_2)$ by imposing the relation that the Casimir equals $\frac{\lambda^2-1}{2}$ and defined $\fr{po}_{\lambda}$ as the fixed-point Lie algebra of a natural involution. Here $\fr{gl}_{\lambda}$ interpolates the Lie algebras $\fr{gl}_n$, while $\fr{po}_{\lambda}$ interpolates the symplectic and odd orthogonal series $\fr{sp}_{2n}$ and $\fr{so}_{2n+1}$ at the expected specializations. Khesin and Malikov showed in \cite{Khesin_Malikov_DS_reduction} that the Poisson geometry of pseudodifferential symbols provides the complex-rank analogue of Drinfeld--Sokolov reduction for these Lie algebras, yielding $\cW(\fr{gl}_{\lambda}), \cW(\fr{po}_{\lambda})$.

We recast the constructions of $\cW(\fr{gl}_{\lambda}), \cW(\fr{po}_{\lambda})$ in the language of Poisson vertex algebras using the formalism of Adler--Gelfand--Dickey brackets, introduced by De Sole, Kac, and Valeri \cite{ADG_SKV}, to define the relevant $\lambda$-brackets. This reformulation appears to be new, and it requires extending the Adler-map formalism to suitable localizations of the polynomial ring $\C[T]$. The resulting Poisson vertex algebras $\cW(\fr{gl}_{\lambda})$ and $\cW(\fr{po}_{\lambda})$ interpolate the classical principal $\cW$-algebras in Types $A$, $B$, and $C$ (Propositions \ref{prop:interpolate_gelfand_dickey_projection} and \ref{prop:interpolate_gelfand_dickey_projection_symplectic_orthogonal}). Moreover, the involution $\lambda\leftrightarrow -\lambda$ induces parity anti-isomorphisms:

\begin{equation}
    \Pi_{\lambda}:\cW(\fr{gl}_{\lambda})\rightarrow \cW(\fr{gl}_{-\lambda})
\end{equation}
\begin{equation}
    \Pi_{\lambda}:\cW(\fr{po}_{\lambda})\rightarrow \cW(\fr{po}_{-\lambda})
\end{equation}
(Propositions \ref{prop:anti_iso_complex_rank_Type_A} and \ref{proposition:anti_iso_complex_PVA_side}).

The proof of the main theorem proceeds by comparing explicit generators produced by the two constructions. On the affine side, these are the interpolated higher Segal--Sugawara elements. On the $\cW$-side, they are the coefficients of the pseudodifferential operators appearing in the Adler--Gelfand--Dickey presentation. In the classical case, the corresponding generators match under the Feigin--Frenkel isomorphism at the critical level. We prove that the same identities hold uniformly in the rank parameter, giving the complex-rank isomorphisms.

\subsection*{Main theorem}

The precise complex-rank Feigin--Frenkel isomorphisms are as follows.

\begin{defi}{Theorem \ref{theorem:FF_inter_type_A}}[Interpolated Feigin--Frenkel Isomorphism at the Critical Level: Type $A$]

For every $\alpha\in \C^{\times}$, there is an isomorphism of Poisson vertex algebras:
\begin{equation}
    \fr{f}_{\alpha}:\fr{z}\left(\Lulhat{gl}{\alpha}\right)\rightarrow \cW(\fr{gl}_{\lambda=\alpha})
    \tag{Type A}
\end{equation}
This isomorphism specializes at positive integer rank to the classical Feigin--Frenkel isomorphism at the critical level for $\fr{gl}_n$. Moreover, $\fr{f}_{\alpha}$ is compatible with the interpolated parity functor up to an interpolated Cartan automorphism $\nu_{\alpha}:\fr{z}\left(\Lulhat{gl}{\alpha}\right)\rightarrow \fr{z}\left(\Lulhat{gl}{\alpha}\right)$, in the sense that:
\begin{equation}
    \Pi_{\alpha}\circ \fr{f}_{\alpha}=\fr{f}_{-\alpha}\circ (\Pi_{\alpha}\circ \nu_{\alpha}).
\end{equation}
\end{defi}

\begin{defi}{Theorem \ref{theorem:FF_inter_type_BC}}[Interpolated Feigin--Frenkel Isomorphism at the Critical Level: Types $B$ and $C$]

For every $\alpha\in \C$, there are isomorphisms of Poisson vertex algebras:
\begin{equation}
    \fr{f}_{\alpha}^{B}:\fr{z}\left(\Lulhat{so}{\alpha}\right)\rightarrow \cW(\fr{po}_{\lambda=\alpha-1})
    \tag{Type B}
\end{equation}
\begin{equation}
    \fr{f}_{\alpha}^{C}:\fr{z}\left(\Lulhat{sp}{\alpha}\right)\rightarrow \cW(\fr{po}_{\lambda=\alpha+1})
    \tag{Type C}
\end{equation}
These isomorphisms specialize at positive integer rank to the classical Feigin--Frenkel isomorphisms at the critical level in Types $B$ and $C$. Moreover, these isomorphisms are compatible with the interpolated parity functor in the sense that:
\begin{equation}
    \Pi_{\alpha+1}\circ \fr{f}_{\alpha}^C= \fr{f}^B_{-\alpha}\circ \Pi_{\alpha}.
\end{equation}
\end{defi}

\subsection*{Organization of the paper}

Sections \ref{section:categorical_background}-\ref{section:interpolating_center_critical_level} develop the categorical side of the story: we review interpolating categories, construct categorical vertex algebras and universal affine vertex algebras, and define the center at the critical level in complex rank together with its interpolated higher Segal--Sugawara generators. Section \ref{section:applications} records the applications to Lie superalgebras and interpolated Bethe algebras. Sections \ref{section:DS_complex_motivation}--\ref{section:universal_adler_maps} develop the pseudodifferential side. We review Feigin's Lie algebras of complex size and establish the formalism needed to define $\cW(\fr{gl}_{\lambda})$ and $\cW(\fr{po}_{\lambda})$ as Poisson vertex algebras arising from pseudodifferential symbols. Section \ref{section:finite_rank_generators_and_isomorphism_W_algebras} recalls the classical presentations needed for comparison, and the final section proves the interpolated Feigin--Frenkel isomorphisms at the critical level in Types $A$, $B$, and $C$.

\section*{Acknowledgments}
The author is grateful to Pavel Etingof for suggesting the project, and to both him and Victor Kac for their detailed and valuable feedback.

\part{The Center of \texorpdfstring{$V^k(\fr{g})$}{V-k(g)} at the Critical Level}
\section{Categorical Background}
\label{section:categorical_background}

This section recalls the categorical background used throughout the paper and fixes conventions. For background on monoidal categories and tensor categories, we refer to \cite{MacLane, Tensor_Cat}.

The section is organized as follows. Subsection~\ref{subsection:tensor-category-basics} collects the categorical constructions needed in the paper. This is the most technically dense part of the section, and readers familiar with this formalism may skim it and return to it as needed. Subsection~\ref{subsection:Del_interpolating_category} is the core of the section: it recalls the interpolating categories in types $A$, $B$, and $C$ used throughout the paper, together with their finite-rank specialization functors, their canonical Lie algebras, and the parity dualities relating them. Subsection~\ref{subsection:Ultraproduct_interpolating_category} recalls the ultraproduct realization of interpolating categories, which provides a technical tool used later to compare interpolating constructions with their finite-rank counterparts.
\subsection{Tensor Category Basics}
\label{subsection:tensor-category-basics}
Throughout this subsection, let $\F$ be a field and let $R$ be a commutative $\F$-algebra. We use the following terminology.

An $R$-linear category $\cC$ is called \emph{pseudo-Abelian} if it is additive and every idempotent splits. 
\begin{definition}
A \emph{pseudo-tensor category} over $R$ is an $R$-linear monoidal category $(\cC,\ox,\alpha,l,r,\mathbf{1}_{\cC})$
such that:
\begin{enumerate}
    \item $\cC$ is rigid.
    \item $\cC$ is pseudo-Abelian.
    \item $\End_{\cC}(\mathbf{1}_{\cC})\cong R$.
\end{enumerate}
A \emph{symmetric pseudo-tensor category} over $R$ is a pseudo-tensor category equipped with a symmetric braiding:
\[
c_{X,Y}:X\ox Y\rightarrow Y\ox X \qquad c_{Y,X}\circ c_{X,Y}=\mathrm{Id}_{X\otimes Y}
\]
\end{definition}
By the coherence theorem for symmetric monoidal categories, we suppress associators and unitors from the notation. Thus we write a symmetric pseudo-tensor category by $(\cC,\ox,\mathbf{1}_{\cC},c)$. When it is clear from context, we denote a symmetric pseudo-tensor category simply by $\cC$.

\begin{definition}
  A \emph{pseudo-tensor functor} $F:\cC\rightarrow \cD$ between pseudo-tensor categories over $R$ is a strong monoidal, additive, $R$-linear functor. We denote by $\mathrm{PsTen}_R$ the resulting category of pseudo-tensor categories over $R$.
\end{definition}

The categories considered in this paper are obtained from pre-additive categories by two standard completions: additive completion and Karoubian completion. The additive completion of a pre-additive category $\cC$, denoted by $\cC^{\mathrm{Add}}$, has objects given by finite tuples of objects of $\cC$, with morphisms given by matrices of morphisms in $\cC$ and composition given by matrix multiplication. See \cite[Section 2.5]{Comes_Wilson_GL_T_definition} for more details. Recall that an additive category is called idempotent complete if every idempotent splits. Given an additive category, its \emph{Karoubian completion}, also known as its \emph{idempotent completion}, is defined as follows.
\begin{definition}
The Karoubian completion of an additive category $\cC$ is the category $\cC^{\mathrm{Kar}}$ defined as follows. Objects of $\cC^{\mathrm{Kar}}$ are pairs $(X,e)$ where $X\in \mathrm{Obj}(\cC)$ and $e\in \End_{\cC}(X)$ is idempotent. The morphisms of $\cC^{\mathrm{Kar}}$ are defined as: 
\begin{equation}\Hom_{\cC^{\mathrm{Kar}}}((X,e),(Y,r)):=\{f\in \Hom_{\cC}(X,Y)\mid r\circ f=f=f\circ e\}\label{eqn:morphism_karoubi_completion}
\end{equation} 
\end{definition}
Following \cite{pseudo_tensor_category}, we call $(\cC^{\mathrm{Add}})^{\mathrm{Kar}}$ \emph{the pseudo-Abelian envelope} of $\cC$. 

Throughout this paper, we frequently change the base ring.
\begin{definition}

    If $f:R\rightarrow S$ is a morphism of commutative $\mathbb{F}$-algebras, then this induces a functor:
    \begin{equation}
        f^*:\mathrm{PsTen}_R\rightarrow \mathrm{PsTen}_S
    \end{equation}
    defined on objects $\cC\in \mathrm{PsTen}_R$  as:
   \begin{equation}
       f^*(\cC):= (\cC\otimes_RS )^{\mathrm{Kar}}
   \end{equation}
     where $\cC\ox_R S$ is the following category. The objects of $\cC\ox_R S$ are the same as $\cC$. The morphisms of $\cC\ox_R S$ are given by base change to $S$, i.e., for $X,Y\in \cC$:
    \begin{equation}
        \mathrm{Hom}_{\cC\ox_R S}(X,Y):= \Hom_{\cC}(X,Y)\ox_R S
    \end{equation} 
    The functor $f^*$ acts on pseudo-tensor functors $F:\cC\rightarrow \cD$ as follows.  On morphisms, it is given by:
    \begin{equation}
        (f^*(F))(h):= (F\otimes_R\Id_S)(h) \quad \text{for} \quad  h\in \mathrm{Hom}_{\cC}(X,Y)\otimes_R S
    \end{equation}
    For $X\in \cC$ and an idempotent $e\in \mathrm{End}_{\cC}(X)\otimes_R S$, we have:
    \begin{equation}
        (f^*(F))(X,e):= (F(X), (F\otimes_R\Id_S)(e))
    \end{equation}
\label{definition:base_change_categories}
\end{definition}
It is important to note that taking the Karoubian envelope of $\cC\ox_R S$ is necessary. Changing the base ring can introduce new idempotents which may not split.  The following is routine to check.
\begin{prop}

    Let $f$, $R$, $S$ be as in Definition \ref{definition:base_change_categories}. The functor $f^*:\mathrm{PsTen}_R\rightarrow \mathrm{PsTen}_S$ restricts to a functor between the categories of symmetric pseudo-tensor categories.
\end{prop}

The final technical construction needed in this subsection is the \emph{$\mathrm{Ind}$-completion} of a symmetric pseudo-tensor category $\cC$. This construction allows us to pass from $\cC$ to a category containing the filtered colimits needed to define categorical vertex algebras. As a basic example, the $\mathrm{Ind}$-completion of the category of finite-dimensional vector spaces over $\F$ is equivalent to the category of all vector spaces over $\F$. For a detailed discussion, see \cite[Section 6]{Categories_Sheaves}.

For the reader's convenience, we recall the definition of a filtered category.

\begin{definition}
A category $\mathcal I$ is called \emph{filtered} if:
\begin{enumerate}
    \item $\mathcal I$ is nonempty.
    \item For any two objects $i,j\in \mathcal I$, there exists an object $k\in \mathcal I$ and morphisms $i\to k$ and $j\to k$.
    \item For any two parallel morphisms $f,g:i\to j$, there exists a morphism $h:j\to k$ such that $h\circ f=h\circ g$.
\end{enumerate}
\end{definition}
\begin{definition}
The \emph{Ind-completion} of a category $\cC$, denoted by $\mathrm{Ind}(\cC)$, has as objects filtered diagrams in $\cC$. More precisely, the objects of $\IC$ are functors $F:I\rightarrow \cC$, where $I$ is a  small filtered category. Morphisms between two objects $F:I\rightarrow \cC,G:J\rightarrow \cC$ are given by:
\begin{equation}
    \mathrm{Hom}_{\IC}(F,G):=\varprojlim_{i\in I}\ \varinjlim_{j\in J}\ \mathrm{Hom}(F(i),G(j))
\end{equation}
\end{definition}
One can show that if $\cC$ is a symmetric pseudo-tensor category over $R$, then $\IC$ satisfies all the conditions of a symmetric pseudo-tensor category over $R$ except rigidity. Moreover, there is a natural embedding of $R$-linear additive symmetric monoidal categories $\iota:\cC\rightarrow \IC$ induced by the Yoneda embedding. 
\begin{definition}
An object of $\IC$ isomorphic to some $\iota(X)$ for $X\in \cC$ is called a \emph{finite object of $\IC$}.
\end{definition}
 We use this embedding $\iota:\cC\rightarrow \IC$ without further mention for the rest of the paper. 

Particularly relevant $\mathrm{Ind}$-objects come from considering filtered categories built out of directed sets $(D,\leq)$. 

\begin{example}
\label{example:finite_subset}
Let $S$ be a set. Define $I_S:=\mathrm{Fin}(S)$ to be the category whose objects are finite subsets of $S$ and there is a unique morphism $T\rightarrow T'$ if and only if $T\subset T'$.
In this case $I_S$ is a filtered category. Notice that the map $s\mapsto \{s\}$ identifies the discrete category of $S$, denoted $S_{\mathrm{disc}}$, as a full subcategory of $I_S$. We use this identification without mention.
\end{example}

\begin{definition}[$S$-Graded Object]

\label{defn:S_graded_object}
  Let $S$ be a set and $\cC$ a category. An $S$-graded object in $\IC$ is an object $X$ isomorphic to the colimit of some filtered diagram $F:I_S\rightarrow \cC$ such that if $T\subset S$ is a finite subset, then the maps $\{F(\{t\}\rightarrow T):t\in T\}$ exhibit $F(T)$ as the coproduct of the objects $F(\{t\})$, for $t\in T$.  
\end{definition}
If $X\in \IC$ is an $S$-graded object represented by the filtered diagram $F:I_S\rightarrow \cC$, then one can show that $X$ is the coproduct of the objects $F(\{s\})$, for $s\in S$, in $\IC$. For the sake of brevity, we suppress our choice of filtered diagram and simply write: 
\begin{equation}
   X_s:=F(\{s\}) \quad s\in S, \qquad X=\bigoplus_{s\in S}X_s
\end{equation}

Morphisms of $S$-graded $\IC$ objects can be easily described:
\begin{equation}
    \Hom_{\IC}\left(\bigoplus_{i\in S}X_i,\bigoplus_{j\in S} Y_j\right)=\prod_{i\in S}\bigoplus_{j\in S}\Hom(X_i,Y_j)
\end{equation}
If $\cC$ has a monoidal structure, then tensor product of two $\Z_{\geq 0}$-graded objects $\bigoplus_{i=0}^{\infty} X_i$ and $\bigoplus_{j=0}^{\infty}Y_j$ is another $\Z_{\geq 0}$-graded object, with $\Z_{\geq0}$-grading given by:
\begin{equation}
    \left(\bigoplus_{i=0}^{\infty} X_i\right)\otimes \left(\bigoplus_{j=0}^{\infty}Y_j\right):=\bigoplus_{t=0}^{\infty}\left(\bigoplus_{\substack{(i,j)\in \Z_{\geq 0}^2\\ i+j=t}}X_i\ox Y_j \right)
\end{equation}

The $\mathrm{Ind}$-completion allows us to form tensor algebras.
\begin{definition}

    Let $X\in \IC$. The tensor algebra $(\mathcal{T}(X),m,i)$ is the unital associative $\IC$-algebra given by $\mathcal{T}(X):=\bigoplus_{i=0}^{\infty}X^{\ox i}$ with multiplication induced by the tensor product and the obvious unit $i:\mathbf{1}_{\cC}\rightarrow \mathcal{T}(X)$.
\end{definition}
The algebra structure can be understood via the $\mathrm{Hom}$-set:
\begin{equation}
    \Hom(\mathcal{T}(X)\ox \mathcal{T}(X),\mathcal{T}(X))=\prod_{m,n}\bigoplus_{s} \Hom(X^{\ox m}\ox X^{\ox n}, X^{\ox s})
\end{equation}
Multiplication is defined by specifying that for each $m,n$, the corresponding map is zero everywhere except for $s=m+n$, in which case it is the identity. 

 Symmetric pseudo-tensor categories allow for more interesting internal algebraic constructions. Namely, a symmetric pseudo-tensor category $(\cC,\ox,\mathbf{1}_{\cC},c)$ has a natural action of the symmetric group $S_n$ on every object $X^{\ox n}$ through the braiding. By abuse of notation, we suppress the braiding when talking about this symmetric group action and denote it simply by the corresponding element of $S_n$. This allows us to define algebraic objects in a symmetric pseudo-tensor category, such as Lie algebras and commutative algebras. See \cite[Section 9.9]{Tensor_Cat} for more details.

If $\cC$ is a symmetric pseudo-tensor category over a $\mathbb{F}$-algebra of characteristic $0$, then one can define symmetric and exterior algebras. More precisely, for $X\in \IC$ and $n\in \Z_{+}$, define idempotents on $X^{\ox n}$ by:
\begin{equation}
    e_{X,n}:=\frac{1}{n!}\sum_{\sigma\in S_n}\sigma
\qquad    h_{X,n}:=\frac{1}{n!}\sum_{\sigma\in S_n}\mathrm{sgn}(\sigma)\cdot \sigma
\label{eqn:symmetrizer_anti_symmetrizer}
\end{equation}
\begin{definition}
    
The \emph{$n$th symmetric power of $X$} and \emph{$n$th exterior power of $X$} are defined respectively as:
\begin{equation}
    S^nX:=(X^{\ox n},e_{X,n})\qquad
\bigwedge^nX:= (X^{\ox n}, h_{X,n})
\end{equation}
The \emph{symmetric algebra of $X$} and \emph{exterior algebra of $X$} are defined respectively as:
\begin{equation}
    S(X):=\bigoplus_{n=0}^{\infty}S^nX \ \ \ \  \ \ \bigwedge(X):= \bigoplus_{n=0}^{\infty}\bigwedge^nX 
\end{equation}
\end{definition}
We denote the idempotent associated to $S(X)$ by $\mathrm{sym}_X$ and call it the symmetrization map. Using the $\IC$-algebra structure of $\mathcal{T}(X)$, one can use the idempotents to define an $\IC$-algebra structure on $S(X)$. A similar construction holds for the exterior algebra. 

\subsection{Deligne's Interpolating Categories $\irep(\GL_{\alpha},\F), \irep(\Or_{\alpha},\F), \irep(\Sp_{\alpha},\F)$}
\label{subsection:Del_interpolating_category}

We now recall the interpolating categories that provide the categorical setting for the rest of the paper.

\subsubsection{Interpolating Type $\GL_n$}

As mentioned in the introduction, the category $\irep(\GL_{\alpha},\F)$ should be viewed as an interpolation of the categories $\mathrm{Rep}(\GL_n,\F)$ of finite-dimensional representations of $\GL_n$.  

For the sake of motivation, we restrict to the case $\F=R=\C$. In this setting, let $V:=\C^{n}$, $W:=V^{*}$, and $[r,s]:=V^{\ox r}\ox W^{\ox s}$ for $r,s\in \Z_{\geq0}$. From classical representation theory, we know that every object of $\mathrm{Rep}(\GL_n,\C)$ is a direct summand of some $[r,s]$. In other words, if we take the full subcategory $\mathrm{Rep}_0(\GL_n,\C)$ consisting of objects $[r,s]$, then the pseudo-Abelian envelope will be $\mathrm{Rep}(\GL_n,\C)$. Note that the parametrization of objects in $\mathrm{Rep}_0(\GL_n,\C)$ is independent of $n$. This suggests that if we can describe the morphism spaces in a way that depends polynomially on $n$, then we can interpolate. By mixed Schur--Weyl
duality, there is a surjective $\C$-algebra homomorphism:
\begin{equation}
    B_{r,s}(n)\rightarrow \mathrm{End}_{\GL_n}([r,s])
\end{equation}
such that it is an isomorphism when $r+s\leq n$. Here $B_{r,s}(n)$ denotes the walled Brauer algebra, a diagrammatic algebra whose basis consists of walled $(r,s)$-tangles, with the wall placed between $r$ and $r+1$. The algebra structure of $B_{r,s}(\alpha)$ is given by stacking diagrams and multiplying by $\alpha^{\ell}$ where  $\ell$ is the number of connected components in the middle of the two stacked diagrams. In other words, the algebra structure of $B_{r,s}(\alpha)$ depends polynomially on $\alpha$. To interpolate, we should use a structure similar to the walled Brauer algebra. To make the category sufficiently nice, we then take the pseudo-Abelian envelope.

Interpolating categories were first defined by Deligne in \cite{St_Interpolating_Deligne}. For our purposes, we recall the form given in \cite[Section 3.1]{Comes_Wilson_GL_T_definition}.

\begin{definition}
\label{definition:word_walled_Brauer_diagram}
    Let $w,w'$ be finite words in the alphabet $\{\bullet,\circ\}$. A $(w,w')$-diagram is a graph satisfying the following conditions:
    \begin{enumerate}
        \item The vertices are positioned in two horizontal rows.
        \item The vertices in the bottom row are labeled by the letters of $w$, from left to right.
        \item The vertices in the top row are labeled by the letters of $w'$, from left to right. 
        \item If an edge connects two vertices in the same row, then those vertices have opposite labels.
    \end{enumerate}
\end{definition}

If $X$ is a $(w,w')$-diagram and $Y$ is a $(w',w'')$-diagram, then one can stack the graphs $Y\star X$ by putting the graph $Y$ on top of $X$, and then restricting the graph to the top and bottom words to obtain a $(w,w'')$-graph denoted by $ Y\cdot X$. Let $\ell(X,Y)$ denote the number of connected components of $Y\star X$ minus the number of connected components of $Y\cdot X$. See \cite[Section 3.1]{Comes_Wilson_GL_T_definition} for examples of such calculations. 

The following category is the one whose pseudo-Abelian envelope we consider. Here, $\F$ is any field, $R$ is a commutative $\F$-algebra, and $\alpha\in R$.

\begin{definition}[Interpolating Skeleton in Type $A$]
\label{definition:interpolating_skeleton_type_A}

    The category $\irep_0(\GL_{\alpha},R)$ consists of the following:
    \begin{enumerate}
        \item Objects: Finite words $w$ in $\bullet$ and $\circ$.
        \item Morphisms: For objects $w,w'$, the $\mathrm{Hom}$-set $\Hom(w,w')$ is the free $R$-module with basis labeled by $(w,w')$-diagrams. 
        \item Composition: If $X\in \Hom(w,w'), Y\in \Hom(w',w'')$ are $(w,w')$-diagrams and $(w',w'')$-diagrams respectively, then $Y\circ X:=\alpha^{\ell(X,Y)}Y\cdot X$. Composition is defined by extending through linearity.
        \item Tensor Product: The tensor product of objects is the concatenation of words. The tensor product of morphisms is the concatenation of diagrams. 
        \item Braiding: For words $w,w'$, the braiding $c_{w,w'}$ is the $(ww',w'w)$-diagram that interchanges the two blocks while preserving the order within each block.
    \end{enumerate}
\end{definition}
This will be an $R$-linear symmetric monoidal category that is pre-additive and rigid. One should think of $\bullet$ as the standard representation $V$ and $\circ$ as the dual representation $V^*$. The interpolating category is given by taking the pseudo-Abelian envelope. 
\begin{definition}[Deligne's Interpolating Category $\irep(\GL_{\alpha},R)$]
\label{definition:interpolating_type_A}
\begin{equation}
    \irep(\GL_{\alpha},R):= \left(\irep_0(\GL_{\alpha},R)^{\mathrm{Add}}\right)^{\mathrm{Kar}}
\end{equation}
\end{definition}
Here are some relevant properties we will use:
\begin{prop}[\cite{St_Interpolating_Deligne}]
\label{prop:properties_type_A_interpolating}
\begin{enumerate}
    \item $\irep(\GL_{\alpha},R)$ is a symmetric pseudo-tensor category over $R$.
    \item If $R=\mathbb{F}=\C$ and $\alpha=n\in \Z_+$, then there is a full braided pseudo-tensor functor: $\ssF{n} :\irep(\GL_{n},\mathbb{C})\rightarrow \mathrm{Rep}(\GL_n,\C)$ such that $\ssF{n}(\bullet)=\C^n$.
    \item If $\alpha\in \F$ and the walled Brauer algebra $B_{r,s}(\alpha)$ is semisimple for all $r,s\in \Z_{\geq 0}$, then $\irep(\GL_{\alpha},\F)$ is a semisimple symmetric tensor category over $\mathbb{F}$. In particular, if $\F=\C$ and $\alpha\notin \Z$ then $\irep(\GL_{\alpha},\C)$ is semisimple. 
\end{enumerate}
\end{prop}
\begin{remark}
In the case $R=\mathbb{F}=\C$, one can show that $\mathrm{Rep}(\GL_n,\C)$ may be identified with the semisimplification of $\irep(\GL_n,\C)$. For this reason, we denote the functor in Proposition~\ref{prop:properties_type_A_interpolating} by $\mathrm{ss}_n$.
\end{remark}
Notice that there will be a canonical Lie algebra $\Lul{gl}{\alpha}:=(\bullet\circ, [-,-])$ where $[-,-]: (\bullet\circ)\ox (\bullet \circ)\rightarrow \bullet \circ $ is given by:
\begin{equation}
    [-,-]_{\Lul{gl}{\alpha}}:= (\Id\ox \mathrm{ev}_{\bullet}\ox \Id)\circ (\Id_{\Lul{gl}{\alpha}^{\ox 2}}-c_{\Lul{gl}{\alpha},\Lul{gl}{\alpha}})
\end{equation}
 If $\ssF{n}$ is the functor from Proposition \ref{prop:properties_type_A_interpolating}, then $\ssF{n}(\Lul{gl}{n})=\fr{gl}_n$. If $\alpha\in R$ is invertible, then there exists a Lie algebra $\Lul{sl}{\alpha}:=(\Lul{gl}{\alpha},\Id_{\bullet\circ}-\frac{1}{\alpha}(\mathrm{coev}_{\bullet}\circ  \mathrm{ev}_{\circ} ))$ in $\irep(\GL_{\alpha},R)$, and this decomposes $\Lul{gl}{\alpha}$ as:
\begin{equation}
    \Lul{sl}{\alpha}\oplus \mathbf{1}_{\irep(\GL_{\alpha},R)}=\Lul{gl}{\alpha}
\end{equation}

As previously mentioned, one should think of $\bullet$ as the standard representation of $\Lul{gl}{\alpha}$ and $\circ$ the dual. This can be made precise by defining the following actions of $\Lul{gl}{\alpha}$:
\begin{equation}
    \rho_{\bullet}:\Lul{gl}{\alpha}\otimes \bullet\rightarrow \bullet \qquad \rho_{\bullet}:= \Id_{\bullet}\otimes \mathrm{ev}_{\bullet}
\label{equation:standard_rep_action_type_A}
\end{equation}
\begin{equation}
    \rho_{\circ}:\Lul{gl}{\alpha}\otimes \circ \rightarrow \circ \qquad \rho_{\circ}:= -(\Id_{\circ}\otimes \mathrm{ev}_{\bullet})\circ c_{\bullet,\circ\circ}
\label{equation:dual_standard_rep_action_type_A}
\end{equation}
The $\Lul{gl}{\alpha}$ modules $(\bullet,\rho_{\bullet})$, $(\circ,\rho_{\circ})$ are sent to the standard and dual standard representation by $\ssF{n}$ when $\alpha=n\in \Z_+$. 

Let $\irep(\GL_{\alpha},R)^{\mathrm{Tw}}$ denote the same underlying pseudo-tensor category as $\irep(\GL_{\alpha},R)$, but with braiding
\begin{equation}
    c^{\mathrm{Tw}}_{w,w'}:=(-1)^{|w|\cdot |w'|}c_{w,w'}
\end{equation}
on words $w,w'$. Here $|w|$ denotes the length of the word $w$.
The twisted braiding appears in the following parity duality. 
\begin{prop}[\cite{St_Interpolating_Deligne}]
\label{prop:GL_duality}
    There is an equivalence of symmetric pseudo-tensor categories:
    \begin{equation}
\Pi_{\alpha}:\irep(\GL_{\alpha},R)\rightarrow \irep(\GL_{-\alpha},R)^{\mathrm{Tw}},
\end{equation}
\begin{equation}
\Pi_{\alpha}(\bullet):=\bullet,\qquad \Pi_{\alpha}(c_{w,w'}):= c^{\mathrm{Tw}}_{w,w'}=(-1)^{|w|\cdot |w'|} c_{w,w'}
    \end{equation}
Moreover, this equivalence identifies the canonical Lie algebra objects: 
\begin{equation}
    \Pi_{\alpha}(\Lul{gl}{\alpha})=\Lul{gl}{-\alpha},\qquad \Pi_{\alpha}(\Lul{sl}{\alpha})=\Lul{sl}{-\alpha}
\end{equation}
with the second equality requiring $\alpha\in R^{\times}$.
\end{prop}
We use the symbol $\Pi_{\alpha}$ to emphasize that this equivalence interpolates the parity functor from supermathematics. 

Lastly, and most importantly, we need the following `interpolation principle'.
\begin{prop}[Interpolation Principle for $\GL_{\alpha}$]
\label{prop:interpolation_principle_A}

    Fix finite words $w,w'$ composed of $\bullet,\circ$ and $R=\F=\C$. For sufficiently large $N$, the functor $\ssF{N}$  induces an isomorphism:
    \begin{equation}\ssF{N}:\Hom_{\underline{\GL}_{N}}(w,w')\rightarrow \mathrm{Hom}_{\GL_N}(\ssF{N}(w),\ssF{N}(w'))
    \end{equation}
    
\end{prop}
\begin{proof}
    This follows from mixed Schur--Weyl duality for $\GL_N$. 
\end{proof}
\subsubsection{Interpolating Types $\Or_n,\Sp_{2n}$}
Assume throughout this subsection that $\F$ does not have characteristic $2$. 

The interpolating category in the orthogonal setting is quite similar to the Type $A$ case, except we now have a symmetric isomorphism $\theta:V\rightarrow  V^*$. This corresponds to considering words only in the alphabet $\{\bullet\}$ instead of $\{\bullet,\circ\}$.
\begin{definition}[Interpolating Orthogonal Skeleton]
\label{definition:interpolating_orthogonal} 

    The category $\irep_0(\Or_{\alpha},R)$ is the analogous category constructed in  Definition~\ref{definition:interpolating_skeleton_type_A}, except words are now in the alphabet $\{\bullet\}$ and Condition $4$ in Definition~\ref{definition:word_walled_Brauer_diagram} is removed.

\end{definition}

As in Type $A$, we pass to the pseudo-Abelian envelope to obtain the interpolating category.
\begin{definition}[Deligne's Interpolating Category $\irep(\Or_{\alpha},R)$]
    The category $\irep(\Or_{\alpha},R)$ is the pseudo-Abelian envelope of $\irep_0(\Or_{\alpha},R)$. 
\end{definition}
We will use the following Type $B$ analogues of Proposition~\ref{prop:properties_type_A_interpolating}.
\begin{prop}
\label{prop:properties_type_B_interpolating}
The category $\irep(\Or_{\alpha},R)$ is a symmetric pseudo-tensor category over $R$. If $R=\F=\C$ and $\alpha=n\in\Z_{+}$, then there is a full braided pseudo-tensor functor:
\[
    \ssF{n}:\irep(\Or_{n},\C)\rightarrow \mathrm{Rep}(\Or_n)
\]
sending $\bullet$ to the standard representation $\C^n$. Moreover, if the corresponding Brauer algebras are semisimple at $\alpha$, then $\irep(\Or_{\alpha},\F)$ is semisimple. In particular, over $\C$, this holds for $\alpha\notin\Z$.
\end{prop} There is a canonical Lie algebra $\Lul{so}{\alpha}:=((\bullet\bullet,h_{\bullet,2}),[-,-])$. Recall from Equation~\eqref{eqn:symmetrizer_anti_symmetrizer} that $e_{\bullet,2},h_{\bullet,2}$ denote the symmetrizer/anti-symmetrizer respectively on the object $\bullet\bullet$. The following holds for $n\in \Z_+$:
\begin{equation}
    \Id\ox \theta:\ssF{n}(\Lul{so}{n})\xrightarrow{\cong} \fr{so}_n 
\end{equation}
\begin{equation}
    \ssF{n}(\cv_{\bullet})= (\Id\ox \theta^{-1})\circ \cv_{V}
    \label{eqn:theta_coevaluation}
\end{equation}
Similar to the Type $A$ case, $\bullet$ can be thought of as the standard representation of $\Lul{so}{\alpha}$ with action defined through:
\begin{equation}
    \rho^B_{\bullet}:\Lul{so}{\alpha}\otimes \bullet\rightarrow \bullet \qquad \rho^B_{\bullet}:= (\Id_{\bullet}\otimes \mathrm{ev}_{\bullet})\circ (h_{\bullet,2}\otimes \Id_{\bullet})
\label{equation:standard_rep_action_type_B}
\end{equation}
The $\Lul{so}{\alpha}$ module $(\bullet,\rho^B_{\bullet})$ is sent to the standard through $\ssF{n}$ when $\alpha=n\in \Z_+$. 
\begin{remark}
\label{remark:Pfaffian_killed}
    It is important to note that $\ssF{n}(\Lul{so}{n})\cong \fr{so}_n$ as an object internal to $\mathrm{Rep}(\Or_n)$. Consequently, when we are taking $\fr{so}_n$-invariants internal to $\mathrm{Rep}(\Or_n)$, we are really taking $\fr{so}_n$-invariants and $O_n$-invariants. This is relevant, for example, when looking at the center of $\fr{so}_n$ internal to $\mathrm{Rep}(\Or_{n})$ that is $\Hom_{\fr{so_n}}(1,U(\fr{so_n}))^{\Or_n}$. When $n$ is odd, this will be the usual center, but when $n=2m$ is even this will be $Z(\fr{so}_{2m})^{\Z/2}$. 
\end{remark}
Interpolating the symplectic category $\mathrm{Rep}(\Sp_{m},\F)$ is almost identical to interpolating $\mathrm{Rep}\ (\Or_{n},\F)$, except now $\theta:V\rightarrow V^*$ is an anti-symmetric isomorphism. 
\begin{definition}[Interpolating Symplectic Skeleton]
The category $\irep_0(\Sp_{\alpha},R)$ is defined almost identically to $\irep_0(\Or_{\alpha},R)$, except that composition is defined differently. Instead of stacking two diagrams and scaling by a factor of $\alpha^{\ell(Y,X)}$ one should scale by $(-\alpha)^{\ell(Y,X)}$. Additionally, the braiding is defined as in Definition~\ref{definition:interpolating_orthogonal} except that it is multiplied by $(-1)^{|w|\cdot|w'|}$, where $|w|$ is the length of the word $w$. 
\end{definition}

Put differently, $\irep_0(\Sp_{\alpha},R)$ is the symmetric tensor category $\irep_0(\Or_{-\alpha},R)$ with a twisted braiding.

\begin{definition}[Deligne's Interpolating Category $\irep(\Sp_{\alpha},R)$]
 The category $\irep(\Sp_{\alpha},R)$ is the pseudo-Abelian envelope of $\irep_0(\Sp_{\alpha},R)$. 
\end{definition}
 $\irep(\Sp_{\alpha},R)$ is a symmetric pseudo-tensor category over $R$ and when $R=\C$, $\alpha=2n\in \Z_+$ there is a full braided pseudo-tensor functor:
\begin{equation}
    \mathrm{ss}_{2n}:\irep(\Sp_{2n},\C)\rightarrow \mathrm{Rep}(\mathrm{Sp}_{2n})
\end{equation}
sending $\bullet$ to the standard representation. There is also a canonical Lie algebra $\Lul{sp}{\alpha}:=((\bullet\bullet,e_{\bullet,2}),[-,-])$, with the bracket defined in the same way. Additionally:
\begin{equation}\Id\ox \theta:\ssF{2n}(\Lul{sp}{2n})\xrightarrow{\cong} \fr{sp}_{2n}\end{equation}
\begin{equation}
    \ssF{2n}(\cv_{\bullet})= (\Id\ox \theta^{-1})\circ \cv_{V}
\end{equation}
The precise action of $\Lul{sp}{\alpha}$ on $\bullet$ is given through:
\begin{equation}
    \rho^C_{\bullet}:\Lul{sp}{\alpha}\otimes \bullet\rightarrow \bullet \qquad \rho^C_{\bullet}:= (\Id_{\bullet}\otimes \mathrm{ev}_{\bullet})\circ (e_{\bullet,2}\otimes \Id_{\bullet})
\label{equation:standard_rep_action_Type_C}
\end{equation}
The parity duality for the orthogonal/symplectic setting is essentially tautological:
\begin{prop}
    There is a braided pseudo-tensor equivalence:
    \begin{equation}
        \Pi_{\alpha}:\irep({\Sp}_{\alpha},R)\rightarrow\irep(\Or_{-\alpha},R)^{\mathrm{Tw}}
    \end{equation}
    \begin{equation}
        \Pi_{\alpha}(\bullet)=\bullet,\quad  \Pi_{\alpha}(c_{w,w'})=(-1)^{|w|\cdot |w'|}c_{w,w'}
    \end{equation}
Furthermore, there is an equality of Lie algebras:
\begin{equation}
    \Pi_{\alpha}(\Lul{sp}{\alpha})=\Lul{so}{-\alpha}
    \label{eqn:identificaiton_symplectic_neg_orthogonal}
\end{equation}
\end{prop}
See \cite[Section 2]{Pavel_Complex_Rank_2} for details. 
\begin{prop}[Interpolation Principle for $\Or_{\alpha}, \Sp_{\alpha}$]
\label{prop:interpolation_principle_BC}

    Let $w,w'$ be words in $\bullet$, and let $R=\F=\C$. For sufficiently large $N$, the functor $\ssF{N}$ induces an isomorphism: 
    \begin{equation}\ssF{N}:\Hom_{\underline{G}_N}(w,w')\rightarrow \mathrm{Hom}_{G_N}(\ssF{N}(w),\ssF{N}(w'))
    \end{equation}
    where $G_N$ is either $\Or_N$ or $\Sp_N$.
\end{prop}
\begin{proof}
    This follows from Schur--Weyl duality for the orthogonal and symplectic groups. 
\end{proof}

\subsection{Interpolating Categories Through Ultraproducts}
\label{subsection:Ultraproduct_interpolating_category}

Several of the main structural results in this paper are proved by comparing the internal categorical constructions with the ultraproducts of finite-rank objects.  We therefore recall an ultraproduct realization of $\irep(G_{\alpha},\C)$.  For the basics of ultraproducts, see \cite{ultraproducts_source}.  This perspective was first outlined by Deligne for transcendental rank in \cite{St_Interpolating_Deligne}, and later extended to arbitrary complex rank by Harman in \cite{harman2016delignecategorieslimitsrank}.

The relevant realization is the following theorem:

\begin{theorem}[{\cite[Theorem 1.1]{harman2016delignecategorieslimitsrank}}]
\label{thm:ultraproduct_deligne}

Let $(\alpha_n)_{n\geq1}$ be a sequence of positive integers with $\alpha_n \to \infty$, 
and let $(p_n)_{n\geq1}$ be a sequence of prime numbers with $p_n \to \infty$.  
Fix a non-principal ultrafilter $\mathcal U$ on $\mathbb N$, and choose an isomorphism
$
\phi:\prod_{\mathcal U} \overline{\F}_{p_n} \rightarrow \mathbb C$, denoting by 
$
\alpha \in \mathbb C
$
 the image of the sequence $\alpha_1 \in  \F_{p_1}, \alpha_2 \in \F_{p_2}, \dots$ 
under the isomorphism. For each $n\in \mathbb{N}$, denote the tautological representation of $G_{\alpha_n}$ over $\overline{\F}_{p_n}$ by $V_n$. 
Let $\cC_n:=\mathrm{Rep}(G_{\alpha_n},\overline{\F}_{p_n})$ and let $\widehat{\mathcal C}_{\mathcal U}$ denote the ultraproduct category 
$\prod_{\mathcal U} \mathcal C_n$. Set
$V \in \widehat{\mathcal C}_{\mathcal U}$ to be the object corresponding to the ultraproduct of the sequence $(V_n)_{n=1}^{\infty}$. If $\mathcal C_{\mathcal U} \subset \widehat{\mathcal C}_{\mathcal U}$ 
denotes the full  symmetric pseudo-tensor category generated by $V$ 
under duals, tensor products, direct sums, and direct summands, then $\mathcal C_{\mathcal U}$ is equivalent, as a symmetric pseudo-tensor category over $\mathbb C$ to $\irep(G_{\alpha},\C)$. 

\end{theorem}
\begin{example}

    Suppose that $n_0\in \Z_{\geq 0}$ and set $\alpha_n:=n_0+p_n$. Then $\irep(G_{n_0})\cong \cC_{\mathcal{U}}$ for some non-principal ultrafilter on $\mathbb{N}$.  More generally, for any algebraic $\alpha\in \overline{\mathbb{Q}}$ there will be infinitely many primes where the minimal polynomial of $\alpha$ has a root in $\F_p$.  By imitating the integer case one can realize $\irep(G_{\alpha})\cong \cC_{\mathcal{U}}$ through ultraproducts. This realization will be key in describing the center at algebraic ranks.
\end{example}
\begin{remark}
    To obtain transcendental $\alpha\in \C$, one needs to replace the finite fields with $\overline{\mathbb{Q}}$. In this case $\prod_{\mathcal{U}}\overline{\mathbb{Q}}\cong \C$, non-canonically, and Theorem \ref{thm:ultraproduct_deligne} will hold. 
\end{remark}
\begin{corollary} 
\label{corollary:canonical_Lie_algebra_ultraproduct}

    With the same notation as in Theorem~\ref{thm:ultraproduct_deligne} denote by $\fr{g}_{\alpha_n}$ the canonical Lie algebra. Under the identification given in Theorem \ref{thm:ultraproduct_deligne} we have:
    \begin{equation}
        \prod_{\mathcal U}\fr{g}_{\alpha_n}\cong \Lul{g}{\alpha}
    \end{equation}
\end{corollary}

\section{Categorical Vertex Algebras }

In this section, we recall the definition of a vertex algebra in a symmetric pseudo-tensor category over an $\mathbb{F}$-algebra $R$ of characteristic zero and the standard results one would expect to hold.  For general background on vertex algebras we refer to \cite{Vertex_Algebras_for_beginners, Axiomatic_approach_VOA}, and for vertex algebras over commutative rings to \cite{Vertex_Rings}. While this paper was being prepared, Zeng's work \cite{zeng2025largenvertexalgebras} appeared, developing a closely related framework for monoidal categories over a field. Where appropriate, we use this framework rather than reproducing standard arguments in full.

\begin{definition}[Vertex Algebra over $R$ in $\cC$]
\label{defn:cat_vertex_algebra}
Let $\cC$ be a symmetric pseudo-tensor category over $R$ of characteristic $0$. A vertex algebra in $\cC$ is a tuple $(V,Y,\ket{0},\mathsf{T})$ consisting of:
\begin{enumerate}
    \item Space of States: $V\in \IC$
    \item Vacuum Vector: A morphism $\ket{0}:\mathbf{1}_{\cC}\rightarrow V$
    \item Translation Map: A morphism $\mathsf{T}:V\rightarrow V$

    \item State-Field Correspondence: a map:
    \begin{equation}
    Y(-,z):=\sum_{n\in \Z}m_nz^{-n-1}\in \Hom_{\IC}(V\ox V,V)[[z,z^{-1}]]\end{equation}
\end{enumerate}
These satisfy the following conditions:
\begin{enumerate}
    \item Field Condition: For finite objects $X_1,X_2\in \cC$ and morphisms $\eta:X_1\rightarrow V,\xi:X_2\rightarrow V$, there exists $N_{X_1,X_2}\in \Z_{\geq 0}$ such that for $k>N_{X_1,X_2}$ we have:
    \begin{equation}
        m_k\circ (\eta\ox\xi)=0
    \end{equation}
    \item Vacuum Axioms: 
    \begin{equation}
        Y(-,z)\circ (\ket{0}\ox \Id_V)=\Id_V 
    \end{equation}
    \begin{equation}
        Y(-,z)\circ (\Id\ox \ket{0})\in \Hom_{\cC}(V\ox \mathbf{1}_{\cC},V)[[z]] 
    \end{equation}
    \begin{equation}  Y(-,z)\circ (\Id_V\ox \ket{0})|_{z=0}=\Id    \end{equation}
        \item Translation Axioms: 
        \begin{equation}
            \mathsf{T}\circ \ket{0}=0
        \end{equation}
        \begin{equation}
            \mathsf{T}\circ Y(-,z)-Y(-,z)\circ (\Id\ox \mathsf{T})=\pr_zY(-,z)
        \end{equation}
        \item Locality: For finite objects $X_1,X_2\in \cC$ and morphisms $\eta \in \Hom_{\IC}(X_1,V), \xi\in \Hom_{\IC}(X_2,V)$ there exists a non-negative integer $N_{X_1,X_2}\in \Z_{\geq 0}$ such that:
        \begin{equation}
        \begin{split}
            (z-w)^{N_{X_1,X_2}}(Y(-,z)\circ (\Id_V\ox Y(-,w))\circ (\eta\ox \xi\ox \Id_V)-\\ Y(-,w)\circ (\Id_V\ox Y(-,z))\circ (c_{V,V}\ox \Id)\circ (\eta\ox \xi\ox \Id_V))=0
        \end{split}
        \end{equation}
\end{enumerate}
\end{definition}

\begin{definition}[Finite Labeled Field]
\label{defn:field_labeled}

A field labeled by a finite object $X$ acting on $V$ is a collection of morphisms:
    \begin{equation}
        A_n:X\otimes V\rightarrow V, \quad \text{ for all }n\in \Z
    \end{equation}
    such that for any finite object $U\in \cC$ and morphism $\xi:U\rightarrow V$, there exists $N\in \Z_+$ such that for all $k>N$:
    \begin{equation}
        A_{k}\circ (\Id_X\ox \xi)=0
    \end{equation}
\end{definition}

\begin{definition}[Normal Ordered Product]

    Let:
    \begin{equation}
    A(z)=\sum_{n\in \Z}A_nz^{-n-1},\quad B(z)=\sum_{n\in \Z}B_nz^{-n-1}
    \end{equation}
    be two fields labeled by finite objects $X_1,X_2$ respectively. The normally ordered product $:A(z)B(w):$ is defined as:
    \begin{equation}
        :A(z)B(w): \ := A(z)_+\circ (\Id_{X_1}\otimes B(w))+B(w)\circ (\Id_{X_2}\ox A(z)_-)\circ (c_{X_1,X_2}\ox \Id_V)
    \end{equation}
    where $A(z)_+=\sum_{n<0}A_nz^{-n-1}$ and $A(z)_-=\sum_{n\geq 0}A_nz^{-n-1}$. 
\end{definition}
\begin{lemma}[Dong's Lemma {\cite[Lemmas 3.3, 3.4]{zeng2025largenvertexalgebras}}]

    If $A(z),B(z),C(z)$ are mutually local fields labeled by $X_1,X_2,X_3$, then $:A(z)B(z):$ is a field labeled by $X_1\otimes X_2$ and is local with $C(z)$.
\end{lemma}
\begin{prop}[Skew-Symmetry {\cite[Proposition 3.1]{zeng2025largenvertexalgebras}}]

Let $X_1,X_2$ be finite objects in $\cC$ with morphisms $\eta\in \Hom_{\IC}(X_1,V), \xi\in \Hom_{\IC}(X_2,V)$, then the following identity holds:
\begin{equation}
    Y(-,z)\circ (\eta\otimes \xi)=e^{zT}\circ Y(-,(-z))\circ c_{V,V}\circ (\eta\otimes \xi)
\end{equation}
\end{prop}
\begin{theorem}[Borcherds’ Formula {\cite[Theorem~3.2]{zeng2025largenvertexalgebras}}]
\label{theorem:brocherds_formula}

Let $X_1,X_2,X_3$ be finite objects in $\cC$ with morphisms:
\[
\eta : X_1 \to V 
\quad \xi : X_2 \to V
\quad \gamma : X_3 \to V
\]
Then for all integers $r,k,\ell$, we have:
\begin{equation}\label{eq:borcherds}
\begin{aligned}
  &\sum_{j\geq 0} \binom{r}{j}\,
    m_{r+k-j} \circ (m_{\ell+j}\otimes \Id_V)
    \circ (\eta\otimes \xi\otimes \gamma) \\
  &= \sum_{j\geq 0} \binom{\ell}{j} (-1)^j 
     \Bigl[ m_{r+\ell-j}\circ (\Id_V\otimes m_{k+j})
     - (-1)^{\ell}\,
       m_{\ell+k-j}\circ (\Id_V\otimes m_{r+j})
       \circ (c_{V,V}\otimes \Id_V)\Bigr]
     \circ (\eta\otimes \xi\otimes \gamma).
\end{aligned}
\end{equation}
\end{theorem}

\begin{corollary}[{\cite[Proposition 3.2]{zeng2025largenvertexalgebras}}]
\label{corollary:zeng_categorical_OPE}

    Let $X_1,X_2$ be finite objects in $\cC$ with morphisms $\eta:X_1\rightarrow V, \xi:X_2\rightarrow V$, then the following identity holds:
    \begin{equation}
        [Y(\eta,z),Y(\xi,w)]=\sum_{n=0}^{\infty}Y(m_n(\eta\ox \xi),w)\pr^{(n)}_w\delta(z-w)
    \end{equation}
\end{corollary}

\begin{definition}

    Let $(V,Y_V,\ket{0}_V,\mathsf{T}_V), (W,Y_W,\ket{0}_W,\mathsf{T}_W)$ be vertex algebras in a symmetric pseudo-tensor category $\cC$ over $R$ of characteristic $0$. A morphism of vertex algebras is a map $f:V\rightarrow W\in \IC$ such that:
    \begin{equation}
        f\circ Y_V=Y_W\circ (f\otimes f), \quad f\circ \ket{0}_V=\ket{0}_W, \quad f\circ \mathsf{T}_V =\mathsf{T}_W\circ f
    \end{equation}
    \label{defn:morphism_vertex_algebras}
    \end{definition}

\section{Vertex Algebras Associated to Interpolated Affine Lie Algebras}
\label{section:vertex_algebras_universal_affine_interpolated}

Besides free field constructions, vertex algebras associated to affine Lie algebras $\widehat{\fr{g}}$  provide one of
the most accessible sources of vertex algebras. The goal of this section is to carry that construction over to complex rank, for affine Lie algebras living in interpolating categories $\irep(G_{\alpha},R)$ where $R$ is characteristic $0$.

When $R=\C$, interpolating categories admit two complementary realizations: an internal categorical one, reviewed in Subsection \ref{subsection:Del_interpolating_category}, and an ultraproduct realization, reviewed in Subsection \ref{subsection:Ultraproduct_interpolating_category}.  We develop universal affine vertex algebras in both languages and then compare the results. More precisely, the section proceeds in three steps. First, we construct the relevant affine Lie algebras and vacuum modules internally in any symmetric pseudo-tensor category over an $\mathbb{F}$-algebra $R$ of characteristic $0$. Second, in the case of an interpolating category over $\C$, we construct the same objects through bounded ultraproducts of finite-rank vacuum modules in positive characteristic. Third, using these vacuum modules, we construct universal affine vertex algebras in interpolating categories through both realizations and show that the resulting categorical vertex algebras agree when $R=\C$.  The affine Lie algebra construction is mostly formal. The vacuum module is subtler, especially on the categorical side, and some technical details are therefore deferred to the appendix.

\subsection{The Categorical Method}
Throughout this subsection, $G_{\alpha}$ denotes one of $\GL_{\alpha}$, $\Or_{\alpha}$, or $\Sp_{\alpha}$. 

\subsubsection{Affine Lie Algebras in Pseudo-Tensor Categories}
\label{section:affine_lie_algebra_poly}

Let $\cC$ be a symmetric pseudo-tensor category over $R$ of characteristic $0$.  The first step is to mimic the usual passage from a Lie algebra to its loop algebra.  In the categorical setting this means forming the $\mathrm{Ind}$-object:
\begin{equation}
L\fr{g}:=\fr{g}\ox R[t,t^{-1}]
\end{equation}
This should be thought of as a direct sum of copies of $\fr{g}$ indexed by $\Z$, one for each loop degree. 

\begin{definition}[Affinization of a Categorical Lie Algebra]
Let $B:\fr{g}\ox \fr{g}\rightarrow \mathbf{1}_{\cC}$ be a symmetric $\fr{g}$-invariant bilinear form. The affinization of $\fr{g}$ with respect to $B$ is the $\IC$ Lie algebra $(\widehat{\fr{g}}_B,[-,-]|_{\widehat{\fr{g}}_B})$ where $\widehat{\fr{g}}_B:=L\mathfrak{g}\oplus \mathbf{1}_{\cC}$, and the Lie algebra structure is defined as:
\begin{align}
    [\fr{g}t^n,\fr{g}t^m]&:=[-,-]_{\fr{g}}t^{n+m}+n\delta_{m+n,0}B \\ 
[\mathbf{1}_{\cC},\fr{g}t^n] &:= 0, \quad  [\fr{g}t^n,\mathbf{1}_{\cC}]:=0, \quad  [\mathbf{1}_{\cC}, \quad\mathbf{1}_{\cC}]:=0
\end{align}
\label{defn:affinization}
\end{definition}
A straightforward calculation shows that this defines a Lie algebra in $\mathrm{Ind}(\cC)$.
\begin{remark}
More generally, the same definition applies if the central term:
\[n\delta_{m+n,0}B\]
is replaced by a Lie algebra $2$-cocycle:
\[\tilde{B}:L\fr{g}\otimes L\fr{g}\rightarrow \mathbf{1}_{\cC}.\]
By this we mean that $\tilde{B}$ satisfies:
\[\tilde{B}\circ (12)=-\tilde{B}\]
as morphisms $(L\fr{g})^{\otimes 2}\rightarrow \mathbf{1}_{\cC}$, and:
\[\tilde{B}\circ([-,-]_{L\fr{g}}\otimes \Id_{L\fr{g}})+ \tilde{B}\circ([-,-]_{L\fr{g}}\otimes \Id_{L\fr{g}})\circ (123) + \tilde{B}\circ([-,-]_{L\fr{g}}\otimes \Id_{L\fr{g}})\circ (132)
=0\]
as morphisms $(L\fr{g})^{\otimes 3}\rightarrow \mathbf{1}_{\cC}$, where the permutations act through the symmetric structure on $\cC$. In this case, one obtains a Lie algebra structure on
$L\fr{g}\oplus \mathbf{1}_{\cC}$ in $\IC$.
\end{remark}
 An example of a symmetric $\fr{g}$-invariant bilinear form that exists for all Lie algebras in a symmetric pseudo-tensor category is the \emph{Killing form}.
\begin{definition}
The Killing form of a Lie algebra $(\fr{g},[-,-])$ in $\cC$ is defined as:
\begin{equation}
    \mathrm{Kil}_{\fr{g}}=\mathrm{ev}_{\fr{g}^*}\circ ([-,-]\ox \Id)\circ (\Id\ox [-,-]\ox\Id )\circ (\Id\ox \mathrm{coev}_{\fr{g}})
\end{equation}
\end{definition}
In the case of interpolating categories, another example is given through the \emph{standard form}.
\begin{definition}
Consider the interpolating category $\irep(G_{\alpha},R)$ with canonical Lie algebra $\Lul{g}{\alpha}$.
The standard form is the symmetric invariant bilinear form associated to the standard representation $V=\bullet$:
\begin{equation}
    \mathrm{B}_{\Lul{g}{\alpha}}=\mathrm{ev}_{V^*}\circ (\rho_V\ox \Id)\circ (\Id\ox \rho_V\ox\Id )\circ (\Id\ox \mathrm{coev}_V)
\end{equation}
where $\rho_V:\fr{g}\ox \bullet\rightarrow \bullet$ is the module structure map given in Equations \eqref{equation:standard_rep_action_type_A},\eqref{equation:standard_rep_action_type_B}, \eqref{equation:standard_rep_action_Type_C}. 
\end{definition}

In the cases $\Lul{g}{\alpha}=\Lul{so}{\alpha}$ and $\Lul{g}{\alpha}=\Lul{sp}{\alpha}$, and also in the case
$\Lul{g}{\alpha}=\Lul{sl}{\alpha}$ when $\alpha\in R^{\times}$, there exists $h^{\vee}_{\Lul{g}{\alpha}}\in R$ such that:
\begin{equation}
    \mathrm{Kil}_{\Lul{g}{\alpha}}=h^{\vee}_{\Lul{g}{\alpha}}\cdot \mathrm{B}_{\Lul{g}{\alpha}}
\end{equation}
We refer to $h^{\vee}_{\Lul{g}{\alpha}}$ as the \emph{dual Coxeter number of} $\Lul{g}{\alpha}$ and $-h^{\vee}_{\Lul{g}{\alpha}}$ is the \emph{critical level}. By examining the classical case we see that the dual Coxeter numbers are:
\begin{equation}
      h^{\vee}_{\Lul{so}{\alpha}}=\alpha-2, \  h^{\vee}_{\Lul{sp}{\alpha}}=\alpha+2, \ h^{\vee}_{\Lul{sl}{\alpha}}=2\alpha
\end{equation}
We make the convention that:
\begin{equation}
    h^{\vee}_{\Lul{gl}{\alpha}}=2\alpha
\end{equation}
so that the critical level of $\Lul{gl}{\alpha}$ is $-2\alpha$. 
It is important to notice that the usual Coxeter number of $\fr{sp}_{2n}$ is usually $n+1$, but in this case we have $h^{\vee}_{\Lul{sp}{2n}}=2n+2$ due to our normalization of $B$.
Similarly, the dual Coxeter number of $\Lul{sl}{\alpha}$ is twice the usual one due to our normalization. 
\begin{remark}
    Alternatively, one can define the critical level through quadratic Lie algebras as explained in \cite{Pavel_Complex_Rank_2}.
\end{remark}

\begin{notation}
    For the sake of simplicity we will denote the affinization of $\Lul{sp}{\alpha},\Lul{so}{\alpha}$ with respect to their standard forms simply by $\Lulhat{sp}{\alpha}$ and $\Lulhat{so}{\alpha}$.
We denote the affinization of $\Lul{gl}{\alpha}$ with respect to $\frac{1}{2\alpha}\mathrm{Kil}_{\Lul{gl}
 {\alpha}}$ by $\Lulhat{gl}{\alpha}$.

 \label{notation:quadratic_form}
\end{notation}
The interpolated parity functor will interact with the symmetric invariant bilinear forms as follows:
\begin{equation}
   \Pi_{\alpha}(\mathrm{Kil}_{\Lul{sp}{\alpha}})=\mathrm{Kil}_{\Lul{so}{-\alpha}}, \  \Pi_{\alpha}(\mathrm{B}_{\Lul{sp}{\alpha}})=-\mathrm{B}_{\Lul{so}{-\alpha}}
\end{equation}
\begin{equation}
   \Pi_{\alpha}(\mathrm{Kil}_{\Lul{gl}{\alpha}})=\mathrm{Kil}_{\Lul{gl}{-\alpha}}, \  \Pi_{\alpha}(\mathrm{B}_{\Lul{gl}{\alpha}})=-\mathrm{B}_{\Lul{gl}{-\alpha}}
\end{equation}
This implies that:
\begin{equation}
   h^{\vee}_{\Lul{sp}{\alpha}}  =-h^{\vee}_{\Lul{so}{-\alpha}}, \ h^{\vee}_{\Lul{gl}{\alpha}}= -h^{\vee}_{\Lul{gl}{-\alpha}}
\end{equation}
\begin{remark}
    The interpolated parity functors induce isomorphisms of affine Lie algebras in the interpolating categories:
    \begin{equation}
        \Pi_{\alpha}(\Lulhat{sp}{\alpha})\cong \Lulhat{so}{-\alpha}
        \label{eqn:Lie_alg_parity_functor_iso_BC}
    \end{equation}
    \begin{equation}
   \Pi_{\alpha}(\Lulhat{gl}{\alpha})\cong \Lulhat{gl}{-\alpha}
        \label{eqn:Lie_alg_parity_functor_iso_A}
    \end{equation}
where the isomorphisms  multiply the tensor unit by $-1$.
\end{remark}

\subsubsection{Interpolating the Vacuum Module Categorically}

Having defined the affinization $\widehat{\fr{g}}$ of a Lie algebra $\fr{g}$ in a symmetric pseudo-tensor category $\cC$ we will now define its vacuum modules categorically.
For readability, we defer most proofs to the appendix. The corresponding results are stated both here and in the appendix, where their proofs are given.
As a first step we need the universal enveloping algebra of a Lie algebra in $\IC$.
Recall that $\mathcal{T}(\fr{g})$ denotes the tensor algebra of $\fr{g}$, and $S(\fr{g})$ the symmetric algebra of $\fr{g}$.
\begin{theorem}[\ref{theorem:categorical_PBW_appendix}]

    Let $\cC$ be a symmetric pseudo-tensor category over $R$ of characteristic $0$.
If $\fr{g}$ is a Lie algebra in $\IC$, then $S(\fr{g})\in \IC$ has an associative unital algebra structure $U(\fr{g}):=(S(\fr{g}),m_{\star})$ in $\IC$ making it the universal enveloping algebra.
\label{thm:universal_enveloping_algebra_categorical_construction}
\end{theorem}

\begin{remark}

    Notice, the algebra structure map $m_{\star}:S(\fr{g})\otimes S(\fr{g})\rightarrow S(\fr{g})$ induces a map of algebras:
    \begin{equation}
        \pi_{\fr{g}}:(\mathcal{T}(\fr{g}),m)\rightarrow (S(\fr{g}),m_{\star})
    \end{equation}
    In particular, as the symmetric algebra is defined as the object $S(\fr{g}):=(\mathcal{T}(\fr{g}),\mathrm{sym}_{\fr{g}})$ the endomorphism $\pi_{\fr{g}}\in \mathrm{End}_{\IC}(\mathcal{T}(\fr{g}))$ satisfies: 
    \begin{equation}
        \mathrm{sym}_{\fr{g}}\circ \pi_{\fr{g}}=\pi_{\fr{g}}
    \end{equation}
    Furthermore, one can also show that (Lemma \ref{lemma:appendix_symmetrizers_pi_properties}):
    \begin{equation}
        \pi_{\fr{g}}\circ \mathrm{sym}_{\fr{g}}=\mathrm{sym}_{\fr{g}}
    \end{equation}
    This implies that $\pi_{\fr{g}}$ is 
an idempotent as:
    \begin{align*}
        \pi_{\fr{g}}\circ \pi_{\fr{g}}=\pi_{\fr{g}}\circ (\mathrm{sym}_{\fr{g}}\circ \pi_{\fr{g}})=\mathrm{sym}_{\fr{g}}\circ \pi_{\fr{g}} =\pi_{\fr{g}} 
    \end{align*}
One should think of $\pi_{\fr{g}}$ as the categorical generalization of the PBW idempotent given in \cite{On_PBW}.
\end{remark}

\begin{prop}[\ref{theorem:representation_universal_enveloping_algebra}]
    The category of left/right modules of $\fr{g}$ in $\IC$ is equivalent to the category of left/right modules of $U(\fr{g})$ in $\IC$.
\label{prop:reps_lie_reps_universal}
\end{prop}

If $\fr{g}$ is a Lie algebra in $\IC$, we denote the category of left modules of $\fr{g}$ in $\IC$ by $\mathrm{Mod}_{\IC}(\fr{g})$.
\begin{corollary}

    Let $\cC,\cD$ be symmetric pseudo-tensor categories over $R$ of characteristic $0$, and $F:\cC\rightarrow \cD$ a braided pseudo-tensor functor. If $\fr{g}$ is a Lie algebra in $\IC$, the following holds canonically:
    \begin{equation}
        F(U(\fr{g}))\cong U(F(\fr{g}))
    \end{equation}
\end{corollary}
\begin{corollary}

    Let $\cC$ be a symmetric pseudo-tensor category over $R$ of characteristic $0$.
Let $f:R\rightarrow S$ be a morphism of commutative $\mathbb{F}$-algebras of characteristic $0$.
The universal enveloping algebra is invariant under change of base ring:
    \begin{equation}
        f^*(U(\fr{g}))= U(f^*(\fr{g}))
    \end{equation}
\end{corollary}

\begin{definition}[Lie Algebra Triple]
A Lie algebra triple in a symmetric pseudo-tensor category $\cC$ over $R$ of characteristic $0$ is a tuple $(\fr{g},\fr{n}_-,\fr{p})$ such that each entry is a Lie algebra in $\IC$, and there is a decomposition of $\IC$-objects:
\begin{equation}
    \fr{g}=\fr{n}_-\oplus \fr{p}
\end{equation}
such that each summand is a Lie subalgebra of $\fr{g}$.
\end{definition}
Notice a Lie algebra triple does not mean that $\fr{g}$ is isomorphic as a Lie algebra to the direct sum of two Lie algebras. Interaction between the summands is allowed. 
\begin{prop}[Generalized Verma Modules \ref{prop:generalized_verma_module_appendix}]
\label{prop:generalized_verma_module}

Let $\fr{g}$ be a Lie algebra in $\IC$.
For every Lie algebra triple $(\fr{g},\fr{n}_-,\fr{p})$
 there is a natural restriction functor:
    \begin{equation}
        \mathrm{Res}^{\fr{g}}_{\fr{p}}:\mathrm{Mod}_{\IC}({\fr{g}})\rightarrow \mathrm{Mod}_{\IC}(\fr{p})
    \end{equation}
    The restriction functor has a left adjoint:
    \begin{equation}
        \mathrm{Ind}^{\fr{g}}_{\fr{p}}:\mathrm{Mod}_{\IC}(\fr{p})\rightarrow\mathrm{Mod}_{\IC}({\fr{g}})
    \end{equation}
    such that if $\cD$ is another category satisfying the same assumptions and $F:\cC\rightarrow \cD$ is a braided pseudo-tensor functor, then:
    \begin{equation}
        F\circ \mathrm{Ind}_{\fr{p}}^{\fr{g}}=\mathrm{Ind}^{F(\fr{g})}_{F(\fr{p})}\circ F
    \end{equation}
Furthermore, these will commute with change of base ring.
\end{prop}

We developed the preceding results in order to apply them to the following special case.
\begin{definition}
\label{definition:vacuum_module}

    Let $\cC$ be a symmetric pseudo-tensor category over $R$ of characteristic $0$. Let $\fr{g}$ be a Lie algebra in $\cC$ and $\tilde{B}:L\fr{g}\otimes L\fr{g}\rightarrow \mathbf{1}_{\cC}$ a Lie  algebra $2$-cocycle in $\IC$. Then its affinization $(\widehat{\fr{g}}_{\tilde{B}},[-,-]|_{\widehat{\fr{g}}_{\tilde{B}}})$ has a natural Lie algebra triple given by:
    \begin{equation}
        \widehat{\fr{g}}_{\tilde{B}}=\fr{g}[t^{-1}]t^{-1}\oplus ((\fr{g}\oplus \mathbf{1}_{\cC})\oplus \fr{g}[t]t)
    \end{equation}
Namely, $\fr{n}_-=\fr{g}[t^{-1}]t^{-1}, \fr{p}=(\fr{g}\oplus \mathbf{1}_{\cC})\oplus \fr{g}[t]t$. If $X$ is a $\fr{g}$-module, then this can be extended to a $\fr{p}$-module by having $\mathbf{1}_{\cC}$ act as $k\in R$ and $\fr{g}[t]t$ act by zero.
Denote the corresponding induced module as:
\begin{equation}M(X,k):=(\mathrm{Ind}_{\fr{p}}^{\widehat{\fr{g}}_{\tilde{B}}}(X), \rho_{X,k})
\end{equation}
In the special case $X=\mathbf{1}_{\cC}$ where $\fr{g}$ acts trivially we define the \emph{vacuum module} as:
\begin{equation}
    V^k(\widehat{\fr{g}}_{\tilde{B}}):= (\mathrm{Ind}^{\widehat{\fr{g}}_{\tilde{B}}}_{\fr{p}}(\mathbf{1}_{\cC} ), \rho_k)
\end{equation}
\end{definition}

In Definition \ref{definition:vacuum_module} we include the $2$-cocycle as part of the data since it affects the resulting vertex algebra.
For example, in the case that $\fr{g}$ is Abelian one may get the free boson or free fermion depending on the $2$-cocycle.   When the $2$-cocycle  is clear from context, we suppress it from the notation and write simply $V^k(\fr{g})$. In particular, in the interpolating category $\irep(G_{\alpha},R)$, the vacuum module at level $k\in R$ associated with the affine Lie algebra $\Lulhat{g}{\alpha}$, formed using the bilinear form fixed in Notation~\ref{notation:quadratic_form}, will be denoted by
 $V^k(\Lul{g}{\alpha})$. 
Due to Proposition \ref{prop:generalized_verma_module}, the vacuum module will be preserved under braided pseudo-tensor functors and change of base ring.
\begin{corollary}

\label{corollary:functor_and_base_change}
    Let $\cC,\cD$ be symmetric pseudo-tensor categories over $R$ of characteristic $0$ and $F:\cC\rightarrow \cD$ a braided pseudo-tensor functor. 
If $\fr{g}, \tilde{B}$ are as in Definition \ref{definition:vacuum_module}, then there is a canonical isomorphism:
    \begin{equation}
        F(V^k(\widehat{\fr{g}}_{\tilde{B}}))\cong V^k(\widehat{F(\fr{g})}_{F(\tilde{B})})
    \end{equation}
Furthermore, if $\cC$ is a symmetric pseudo-tensor category over $R$ of characteristic $0$ and $f:R\rightarrow S$ is a morphism of commutative $\mathbb{F}$-algebras, then there is a canonical isomorphism:
\begin{equation}
    f^*(V^k(\widehat{\fr{g}}_{\tilde{B}}))\cong V^{f(k)}\bigl(\widehat{f^*(\fr{g})}_{f^*(\tilde{B})}\bigr)
\end{equation}
\end{corollary}
As a consequence of Equations \eqref{eqn:Lie_alg_parity_functor_iso_BC}, \eqref{eqn:Lie_alg_parity_functor_iso_A} we have an isomorphism of modules: 
\begin{equation}
\Pi_{\alpha}(V^k(\Lul{so}{\alpha}))\cong V^{-k}(\Lul{sp}{-\alpha})
\label{eqn:module_iso_level_change}
\end{equation}
\begin{equation}
\Pi_{\alpha}(V^k(\Lul{gl}{\alpha}))\cong V^{-k}(\Lul{gl}{-\alpha})
\label{eqn:module_iso_level_change_gl}
\end{equation}
\subsection{The Ultraproduct Method}

Throughout this subsection, $G_{\alpha}$ denotes one of $\GL_{\alpha}$, $\Or_{\alpha}$, or $\Sp_{\alpha}$.  Fix a sequence of primes $(p_n)_{n\geq1}$ with $p_n\to\infty$, a non-principal ultrafilter $\mathcal{U}$ on $\mathbb{N}$, and a sequence of positive integers $(\alpha_n)_{n=1}^{\infty}$ with $\alpha_n\to\infty$ whose ultraproduct corresponds to $\alpha$ under some identification $\prod_{\mathcal{U}}\overline{\mathbb{F}}_{p_n}\cong \C$.

\subsubsection{Affine Lie algebras through Ultraproducts}

By Corollary \ref{corollary:canonical_Lie_algebra_ultraproduct}, the ultraproduct of the canonical Lie algebras $\prod_{\mathcal{U}}\fr{g}_{\alpha_n}$ is isomorphic to $\Lul{g}{\alpha}$ under Theorem \ref{thm:ultraproduct_deligne}.  The subtle point is that loop algebras do not commute naively with ultraproducts: the unrestricted ultraproduct $\prod_{\mathcal{U}}L\fr{g}_{\alpha_n}$ contains terms of unbounded degree and is therefore much larger than $L\Lul{g}{\alpha}$.  The correct construction is obtained by imposing boundedness in the degree.  More generally, recall that the category of finite subsets of $\Z$, denoted $I_{\Z}$, is filtered under inclusion.

\begin{definition}[Restricted Ultraproduct]

Let $\cC_n$ be a family of categories indexed by $\mathbb{N}$, and denote the category formed by their ultraproduct as $\widehat{\mathcal{C}}_{\mathcal{U}}$.

If there is a collection of functors $F_n:I_{\Z}\rightarrow \cC_n$, this induces a functor:

\begin{equation}
    \prod_{\mathcal{U}}^{\Delta}F_n:I_{\Z}\rightarrow  \widehat{\mathcal{C}}_{\mathcal{U}}, \left(\prod_{\mathcal{U}}^{\Delta}F_n\right)(S):=\prod_{\mathcal{U}}F_n(S), \left(\prod_{\mathcal{U}}^{\Delta}F_n\right)(f):=\prod_{\mathcal{U}}F_n(f)
\end{equation}

For each $n\in \mathbb{N}$, denote the Ind-object corresponding to $F_n$ by $X_n^{F_n}$. The restricted ultraproduct of $(X_n^{F_n})_{n=1}^{\infty}$ is by definition the $\mathrm{Ind}(\widehat{\mathcal{C}}_{\mathcal{U}})$-object:
\begin{equation}
  \prod_{\mathcal{U}}^{\mathrm{B}}X_n^{F_n}:=\mathrm{colim}\left(\prod_{\mathcal{U}}^{\Delta}F_n\right)
\end{equation}
\end{definition}

\begin{lemma}

    Let $\cC_n:=\mathrm{Rep}(G_{\alpha_n},\overline{\mathbb{F}}_{p_n})$, consider $L\fr{g}_{\alpha_n}$ with the natural $\Z$-grading.
There is an isomorphism of Lie algebras in $\mathrm{Ind}(\irep(G_{\alpha},\C))$:
    \begin{equation}
         \prod_{\mathcal{U}}^{\mathrm{B}}L\fr{g}_{\alpha_ n}\cong L\Lul{g}{\alpha}
    \end{equation}
under the identification given in Theorem \ref{thm:ultraproduct_deligne}.

Furthermore:
\begin{equation}
    \mathrm{Kil}_{\Lul{g}{\alpha}}=\prod_{\mathcal{U}}\mathrm{Kil}_{\Lul{g}{\alpha_n}}
\end{equation}
and so:
\begin{equation}
    \prod_{\mathcal{U}}^B \widehat{\fr{g}}_{\alpha_n} \cong \Lulhat{g}{\alpha}
\end{equation}
\end{lemma}
\begin{proof}
 This follows from the fact that the ultraproduct commutes with finite direct sums.
\end{proof}
The next step is to prove that the bounded ultraproduct of universal enveloping algebras coincides with the categorical construction.
This requires we have control over the filtration given by the categorical universal enveloping algebra construction.
If $\fr{g}$ is a Lie algebra over $\mathbb{Q}$, then there are two different but isomorphic constructions of the universal enveloping algebra:

\begin{enumerate}

    \item The categorical construction in $\cC=\mathrm{Vect}_{\mathbb{Q}}^{\mathrm{fd}}$, denoted by:
    \begin{equation*}
        \fr{U}(\fr{g})=(S(\fr{g}),m_{\star},1)
    \end{equation*}
    with filtration $F^{\bullet}\fr{U}(\fr{g})$:
    \begin{equation*}
        F^i\fr{U}(\fr{g}):=\bigoplus_{k=0}^iS^k(\fr{g})
    \end{equation*}
    \item The quotient construction where:
    \begin{equation*}
        U(\fr{g}):=(T(\fr{g})/\langle [x,y]-(x\otimes y -y\otimes x)\rangle ,m,1)
    \end{equation*}
    with filtration $F^{\bullet}U(\fr{g})$ induced by the PBW theorem.
\end{enumerate}

The categorical construction over $R$ requires that $m!\in R^{\times}$ for all $m\in \Z_+$ as the multiplication is defined using symmetrizers.
We assumed $R$ contains a field of characteristic $0$ so this was not an issue, but to show our ultraproduct construction agrees with the categorical construction we need to contend with Lie algebras in positive characteristic. The solution is to notice that if $\mathbb{F}_{p}$ is a field of characteristic $p$, $\fr{g}_0$ a Lie algebra over $\mathbb{F}_p$ and $m<p$, then $m!$ is invertible and one may realize $S^m(\fr{g}_0)$ as a direct summand of $\fr{g}_0^{\otimes m}$ through the symmetrizer map.
If $\fr{g}_0$ comes from a Lie algebra with a $\Z$-form, say $\fr{g}_0=\fr{g}_{\Z}\otimes_{\Z}\mathbb{F}_p$, then showing the ultraproduct construction agrees with the categorical construction is equivalent to showing the filtration from the quotient construction $U(\fr{g}_0)$ comes from the filtration of $\fr{U}(\fr{g}_{\mathbb{Q}})$ after a suitable change of base.
In particular, if $\fr{g}$ has a $\Z$-form $\fr{g}_{\Z}$ one can verify that the multiplication on filtered components of $\fr{U}(\fr{g}_{\mathbb{Q}})$ for $r,s\in \Z_{\geq 0}$:
\begin{equation}
    m_{\star,(r,s)}:F^r\fr{U}(\fr{g}_{\mathbb{Q}})\otimes F^s\fr{U}(\fr{g}_{\mathbb{Q}})\rightarrow F^{r+s}\fr{U}(\fr{g}_{\mathbb{Q}})
\end{equation}
may be restricted to be defined over $\Z[\frac{1}{(r+s)!}]$.
Namely, there is a well-defined map:
\begin{equation}
    m_{\star,(r,s)}:F^r\fr{U}(\fr{g}_{\Z[\frac{1}{(r+s)!}]})\otimes F^s\fr{U}(\fr{g}_{\Z[\frac{1}{(r+s)!}]})\rightarrow F^{r+s}\fr{U}(\fr{g}_{\Z[\frac{1}{(r+s)!}]})
\end{equation}

The following Lemma shows that such a restriction agrees with the multiplication of $U(\fr{g}_{\Z[\frac{1}{m!}]})$.
\begin{lemma}
Let $\fr{g}_{\Z}$ be a Lie algebra over $\Z$. 
Consider $\fr{U}(\fr{g}_{\mathbb{Q}}):=(S(\fr{g}_{\mathbb{Z}})\otimes_{\Z}\mathbb{Q},m_{\star},1)$ as a $\mathbb{Q}$-algebra.
The symmetrization map induces an isomorphism of $\Z[\frac{1}{m!}]$-modules:
\[
\widetilde{\mathrm{sym}}_m:
F^mS(\fr{g}_{\Z})\otimes \Z\!\left[\frac{1}{m!}\right]
\longrightarrow
F^m U\!\left(\fr{g}_{\Z[\frac{1}{m!}]}\right)
\]
with the following property: for all $r,s\in \Z_{\geq 0}$ with $r+s\leq m$, we have:
\begin{equation}
    m|_{F^rU(\fr{g}_{\Z[\frac{1}{m!}]})\otimes F^sU(\fr{g}_{\Z[\frac{1}{m!}]})}\circ (\widetilde{\mathrm{sym}}_r\otimes \widetilde{\mathrm{sym}}_s)= \widetilde{\mathrm{sym}}_{r+s}\circ m_{\star,(r,s)}
\end{equation}
\label{lemma:filtration_iso}
\end{lemma}
\begin{proof}
By PBW we know there is a $\mathbb{Q}$-algebra isomorphism given by the symmetrizer map $\mathrm{sym}: \fr{U}(\fr{g}_{\mathbb{Q}})\rightarrow U(\fr{g}_{\mathbb{Q}})$.
Restricting $\mathrm{sym}$ to the subset $F^mS(\fr{g}_{\mathbb{Z}})\otimes \Z[1/m!]$ we see that we obtain a $\Z[1/m!]$-linear unipotent map with respect to the PBW filtration, and hence an isomorphism of $\Z[1/m!]$-modules.
In particular, we see it will be compatible with the prescribed multiplications as $\mathrm{sym}$ is an algebra isomorphism over $\mathbb{Q}$.
\end{proof}
 Denote the canonical Lie algebras of  $\mathrm{Rep}(G_{\alpha_n},\overline{\mathbb{F}}_{p_n})$ by $\fr{g}_n$ and fix a non-principal ultrafilter $\mathcal{U}$ on $\mathbb{N}$.
We can take the restricted ultraproducts of $U(\fr{g}_{\alpha_n})$ and the result will be the categorical universal enveloping algebra as the next result shows.
\begin{prop}
    Under the identification given in Theorem \ref{thm:ultraproduct_deligne}, we have: \[\prod_{\mathcal{U}}^BU(\fr{g}_{\alpha_n})\in \mathrm{Ind}(\irep(G_{\alpha}))\] 
    Furthermore, there is an isomorphism of associative unital algebras in $\mathrm{Ind}(\irep(G_{\alpha}))$:
    \begin{equation}
            \prod_{\mathcal{U}}^BU(\fr{g}_{\alpha_n})\cong U(\Lul{g}{\alpha})
    \end{equation}
\label{prop:ultraproduct_universal_Lie}
\end{prop}
\begin{proof}
If $p_n>m$, then $m!$ is invertible in $\mathbb{F}_p$ and $F^mU(\fr{g}_{\Z[1/m!]})\otimes_{\Z[1/m!]}\overline{\mathbb{F}_p}\cong F^mU(\fr{g}_{\overline{\mathbb{F}_p}})$.
This implies that for $p_n>m$ we have an isomorphism $F^m\fr{U}(\fr{g})\otimes_{\Z[1/m!]}\overline{\mathbb{F}_p}\cong F^mU(\fr{g}_{\ov{\mathbb{F}}_p})$ compatible with multiplication within the filtered part by Lemma \ref{lemma:filtration_iso}.
But evidently, this implies that:
\begin{equation}
    F^mU(\Lul{g}{\alpha})\cong \prod_{\mathcal{U}}F^mU(\fr{g}_{\alpha_n})
\end{equation}
as filtered associative unital algebras.
\end{proof}
A similar result holds for the canonical affine Kac--Moody algebras in an interpolating category.
\begin{prop}
\label{prop:ultraproduct_universal_Kac_Moody}
Assume that $\alpha\in \C^{\times}$. Under the identification given in Theorem \ref{thm:ultraproduct_deligne}, we have:
\[\prod_{\mathcal{U}}^B U(\widehat{\fr{g}}_{\alpha_n})\in\mathrm{Ind}(\irep(G_{\alpha}))\]
Furthermore, there is an isomorphism of associative unital algebras in $\mathrm{Ind}(\irep(G_{\alpha}))$:
\begin{equation}
    \prod_{\mathcal{U}}^B U(\widehat{\fr{g}}_{\alpha_n})
    \cong
    U(\Lulhat{g}{\alpha})
\end{equation}
\end{prop}
\begin{proof}
In this case, the argument proceeds as in Proposition \ref{prop:ultraproduct_universal_Lie}, mutatis mutandis.
Namely, one proves the affine analogue of Lemma \ref{lemma:filtration_iso}. The $\Z_{\geq 0}$-filtration is replaced by a $\Z_{\geq 0}\times \Z$-filtration, and the resulting map is an isomorphism over $\Z[\frac{1}{m!}][\frac{1}{h^{\vee}}]$ in Type $A$ and over $\Z[\frac{1}{m!}]$ in Types $B$ and $C$.
With this established, one observes that in Type $A$, under the identification $\C\cong \prod_{\mathcal{U}}\overline{\mathbb{F}}_{p_n}$, the elements $\alpha_n\in \mathbb{F}_{p_n}$ are nonzero for almost all $n$. Hence the dual Coxeter numbers $2\alpha_n$ are also nonzero in $\mathbb{F}_{p_n}$ for almost all $n$. This allows us to apply the affine analogue of Lemma \ref{lemma:filtration_iso} in all types and conclude the result.
\end{proof}
Notice that the isomorphisms in Propositions \ref{prop:ultraproduct_universal_Lie}, \ref{prop:ultraproduct_universal_Kac_Moody} will preserve the Lie algebra triple decomposition of $\widehat{\fr{g}}_B$ and so in particular this implies that:
\begin{prop}
    Under the identification given in Theorem \ref{thm:ultraproduct_deligne} there is an isomorphism of $\Lulhat{g}{\alpha}$-modules:
    \begin{equation}
        \prod^B_{\mathcal{U}}V^{k_n}(\widehat{\fr{g}}_{\alpha_n})\cong V^k(\Lulhat{g}{\alpha})   
    \end{equation}
    where $k$ is the number corresponding to the ultraproduct element $[(k_n)_{n\in \mathbb{N}}]$.
\label{prop:ultraproduct_vacuum_module}
\end{prop}
\subsection{Universal Affine Vertex Algebras}
With the categorical and ultraproduct constructions of affine Lie algebras in place, we now pass to the corresponding universal affine vertex algebras.  The first theorem gives the internal categorical construction, after which we compare it with the ultraproduct realization.

\begin{theorem}[Universal Affine Vertex Algebra in Symmetric Pseudo-Tensor Categories]

Let $\cC$ be a symmetric pseudo-tensor category over $R$ of characteristic $0$.
Let $\fr{g}$ be a Lie algebra in $\cC$ and $\tilde{B}:L\fr{g}\otimes L\fr{g}\rightarrow \mathbf{1}_{\cC}$ be a Lie algebra $2$-cocycle in $\IC$.
For $k\in R$ there exists a $\cC$-vertex algebra structure $V^k(\fr{g},{\tilde{B}}):=(V^k(\widehat{\fr{g}}_{\tilde{B}}),Y,\ket{0},\mathsf{T})$  on $V^k(\widehat{\fr{g}}_{\tilde{B}})$ such that the vacuum
$\ket{0}$ is the canonical morphism $\mathbf{1}_{\cC}\rightarrow V^k(\widehat{\fr{g}}_{\tilde{B}})$, the derivation is $\mathsf{T}=-\pr_t$ and:
    \begin{equation}
        Y(\fr{g}t^{-1},z)=\sum_{n\in \Z}\left(\rho_{k}|_{\fr{g}t^n}\right) z^{-n-1}
    \end{equation}
\label{theorem:vertex_algebra_structure_universal_affine}
\end{theorem}
\begin{proof}
The proof of this is essentially the same as the $\mathrm{Vect}_{\C}$-case, and so we omit the more straightforward details  such as the vacuum axioms and translation covariance.
Let the field labeled by $\mathbf{1}_{\cC}$ be the identity field, and for $\fr{g}t^{-1}$ define:
    \begin{equation}
        Y(\fr{g}t^{-1},z):=\sum_{n\in \Z}(\rho_k|_{\fr{g}t^n})z^{-n-1}
    \end{equation}
    where $\rho_k:\widehat{\fr{g}}_{\tilde{B}}\otimes V^k(\widehat{\fr{g}}_{\tilde{B}})\rightarrow V^k(\widehat{\fr{g}}_{\tilde{B}})$ is the module structure.
From the usual proof one can check that this is a field labeled by $\fr{g}t^{-1}$ in the sense of Definition \ref{defn:field_labeled} and is local in the sense of Definition \ref{defn:cat_vertex_algebra}.
Using $\pr^{(n)}_z:=\frac{1}{n!}\pr_z^n$ one obtains a field labeled by $\fr{g}t^{-n-1}$ through:
    \begin{equation}
        Y(\fr{g}t^{-n-1},z):=\pr^{(n)}_zY(\fr{g}t^{-1},z)
    \end{equation}
 and hence a field labeled by $\fr{g}[t^{-1}]t^{-1}$.
Using the categorical Dong's lemma one obtains a local field labeled by $(\fr{g}[t^{-1}]t^{-1})^{\ox n}$, and hence a field labeled by $\mathcal{T}(\fr{g}[t^{-1}]t^{-1})$.
Applying the symmetrization map gives the desired state-field correspondence.
This will be local as any finite object $f:X\rightarrow V^k(\widehat{\fr{g}}_{\tilde{B}})$ factors as a morphism into $\tilde{f}:X\rightarrow \bigoplus^{m}_{k=0}S^k(\fr{g}[t^{-1}]t^{-1})$. 
\end{proof}

\begin{remark}
Applying Theorem \ref{theorem:vertex_algebra_structure_universal_affine} to $\mathrm{Vect}_{\C}$ or $\mathrm{sVect}_{\C}$ yields the usual construction of 
universal affine (super)vertex algebras. 
\end{remark}

The following is immediate from the construction in Theorem \ref{theorem:vertex_algebra_structure_universal_affine} and Corollary \ref{corollary:functor_and_base_change}.
\begin{prop}

    Let $F:\cC\rightarrow \cD$ be a braided pseudo-tensor functor over $R$ of characteristic $0$, and $\fr{g},\tilde{B}$ as in Theorem \ref{theorem:vertex_algebra_structure_universal_affine}.
There is a canonical isomorphism of vertex algebras in $\cD$:
    \begin{equation}
        F(V^k(\widehat{\fr{g}},{\tilde{B}}))\cong V^k(\widehat{F(\fr{g})},F(\tilde{B}))
    \end{equation}
    If $f:R\rightarrow S$ is a morphism of commutative $\mathbb{F}$-algebras of characteristic $0$ and $\cC$ is a symmetric pseudo-tensor category over $R$, then there is a canonical isomorphism:
    \begin{equation}
        f^*(V^k(\widehat{\fr{g}},\tilde{B}))\cong V^{f(k)}(f^*(\widehat{\fr{g}}),f^*(\tilde{B}))
    \end{equation}

\label{prop:vertex_algebra_functor_compatibility}
\end{prop}
Just as for the affinization of certain Kac--Moody algebras, Proposition \ref{prop:vertex_algebra_functor_compatibility} allows us to interpolate certain universal 
affine vertex algebras. Recall that when the $2$-cocycle on $L\fr{g}$ is clear from context, we write the vacuum module $V^k(\widehat{\fr{g}}_{\tilde{B}})$ simply by $V^k(\fr{g})$. By abuse of notation, when the $2$-cocycle is clear from context, we will write the vertex algebra $V^k(\fr{g},{\tilde{B}})$ simply as $V^k(\fr{g})$. 

\begin{definition}[Interpolated Universal Affine Vertex Algebras]
Let $R$ be a $\C$-algebra.
The interpolated universal affine vertex algebra of Type $B$ at $\alpha\in R$ is the construction $V^k(\Lul{so}{\alpha})$ in the category $\irep(\Or_{\alpha},R)$.
The interpolated universal affine vertex algebra of Type $C$ at $\alpha\in R$ is the construction $V^k(\Lul{sp}{\alpha})$ in the category $\irep(\Sp_{\alpha},R)$.
In the Type $A$ case, assume that $\alpha\in R^{\times}$. The interpolated universal affine vertex algebra of Type $A$ at $\alpha$ is the construction $V^k(\Lul{gl}{\alpha})$ in the category $\irep(\GL_{\alpha},R)$.
\end{definition}
They are interpolated in the sense that when $R=\C, \alpha\in \mathbb{N}$ we have by Proposition \ref{prop:vertex_algebra_functor_compatibility} an isomorphism of $\C$-vertex algebras:
\begin{equation}
    \ssF{n}(V^k(\Lul{g}{n}))\cong V^k(\fr{g}_n)
\end{equation}
where $\Lul{g}{n}$ is the canonical Lie algebra in the corresponding interpolating category.
By Equations \eqref{eqn:module_iso_level_change}, \eqref{eqn:module_iso_level_change_gl} we see that the interpolated parity functor $\Pi_{\alpha}$ at $\alpha\in R$ induces a vertex algebra isomorphism in  $\irep(\Or_{-\alpha},R)$:
\begin{equation}
    \Pi_{\alpha}(V^k(\Lul{sp}{\alpha}))\cong V^{-k}(\Lul{so}{-\alpha})
    \label{eqn:VOA_isomorphism_Type_BC}
\end{equation}
and when $\alpha\in R^{\times}$ a vertex algebra isomorphism  in $\irep(\GL_{-\alpha},R)$:
\begin{equation}
      \Pi_{\alpha}(V^k(\Lul{gl}{\alpha}))\cong V^{-k}(\Lul{gl}{-\alpha})
    \label{eqn:VOA_isomorphism_Type_A}
\end{equation}

Alternatively, one can construct $V^k(\Lul{g}{\alpha})$ through ultraproducts which we now briefly explain.
By Proposition \ref{prop:ultraproduct_vacuum_module} we know that $\prod_{\mathcal{U}}^{B}V^{k_n}(\fr{g}_{\alpha_n})\cong V^k(\Lul{g}{\alpha})$ as $\Z_{\geq 0}$-filtered $\Lulhat{g}{\alpha}$-modules.
On the other hand the vertex algebra structure on either side is  given through the $\Lulhat{g}{\alpha}$-module structure and defined inductively using the PBW filtration.
Therefore, we may conclude that:
\begin{prop}

Assume that $\alpha\in \C^{\times}$. Under the identification in Theorem \ref{thm:ultraproduct_deligne} the $\Z_{\geq 0}$-filtered $\Lulhat{g}{\alpha}$-module $\prod_{\mathcal{U}}^{B}V^{k_n}(\fr{g}_{\alpha_n})$ has a categorical vertex algebra structure in $\irep(G_{\alpha},\C)$ induced by the $\overline{\mathbb{F}}_{p_n}$-vertex algebra structure of $V^{k_n}(\fr{g}_{\alpha_n})$.
Furthermore, there is an isomorphism of categorical vertex algebras:
\begin{equation}
    \prod_{\mathcal{U}}^{B}V^{k_n}(\fr{g}_{\alpha_n})\cong V^k(\Lul{g}{\alpha})
\end{equation}
induced by the isomorphism in Proposition \ref{prop:ultraproduct_vacuum_module}
\label{prop:ultraproduct_vacuum_module_vertex_algebra}
\end{prop}
\section{Interpolating the Center at the Critical Level}
\label{section:interpolating_center_critical_level}

As in the previous section, all rings $R$ are associative, commutative, unital algebras over a field $\F$ of characteristic $0$.  The goal of this section is to pass from the universal affine vertex algebras constructed in Section~\ref{section:vertex_algebras_universal_affine_interpolated} to their centers at the critical level, and then to identify explicit algebraically independent generators. The main point is that the finite-rank formulas for higher Segal--Sugawara vectors can be expressed in a form which depends polynomially on the rank, and hence admits a direct interpolation. We use some standard terminology from Poisson vertex algebras in this section. The formal definitions and conventions are recalled in Section~\ref{section:PVA_background}.

\subsection{The Center at the Critical Level and its Poisson Vertex Algebra Structure}

Throughout this subsection, we assume that $\cC$ is a symmetric pseudotensor category over $R$ of characteristic $0$. Additionally, we assume $\fr{g}$ is a Lie algebra in $\cC$ such that its loop algebra has a Lie algebra $2$-cocycle $\tilde{B}:L\fr{g}\otimes L\fr{g}\rightarrow \mathbf{1}_{\cC}$. We first define the center in this generality and then show that these centers naturally carry a Poisson vertex algebra structure compatible with functoriality and base change.

Motivated by the classical realization as $\fr{g}[t]$-invariants, we make the following definition.

\begin{definition}
    Let $k\in R$.
    The center of $V^k(\fr{g},\tilde{B})$ is defined as:
    \begin{equation}
        \fr{z}_k(\widehat{\fr{g}}_{\tilde{B}}):= \Hom_{\fr{g}[t]}(\mathbf{1}_{\cC},V^k(\fr{g},\tilde{B}))
    \end{equation}
\end{definition}

\begin{theorem}

For any $k\in R$, the center $\fr{z}_{k}(\widehat{\fr{g}}_{\tilde{B}})$ has the structure of a Poisson vertex algebra over $R$.
If $F:\cC\rightarrow \cD$ is a braided pseudo-tensor functor over $R$, then there is a canonical morphism of $R$-linear Poisson vertex algebras:
\begin{equation}
    \ov{F}_k:\fr{z}_k(\widehat{\fr{g}}_{\tilde{B}})\rightarrow \fr{z}_k(F(\widehat{\fr{g}})_{F(\tilde{B})})
\end{equation}
Furthermore, if $f:R\rightarrow S$ is a morphism of $\F$-algebras, then there is a canonical morphism:
\begin{equation}
    \overline{f}:f^*(\fr{z}_k(\widehat{\fr{g}}_{\tilde{B}}))\rightarrow  \fr{z}_{f(k)}(f^*(\widehat{\fr{g}})_{f^*(\tilde{B})})
\end{equation}
\label{theorem:invariant_ind_PVA_structure}
\end{theorem}

\begin{proof}
One can show that if $(V,Y,\mathsf{T},\ket{0})$ is a vertex algebra in $\cC$, then $\Hom_{\IC}(\mathbf{1}_{\cC},V)$ is a vertex algebra over $R$ by defining the state-field correspondence through the modes for $a,b\in \Hom_{\IC}(\mathbf{1}_{\cC},V)$ by:
     \begin{equation}
        a_{(n)}b:= m_n\circ (a\otimes_{\cC} b)\in\Hom_{\cC}(\mathbf{1}_{\cC},V)
    \end{equation}
    See \cite{zeng2025largenvertexalgebras} where this is also considered.
    To see that the center is a vertex sub-algebra, one applies Borcherds' formulas in Theorem \ref{theorem:brocherds_formula}. Observe that, for $\eta\in \fr{z}_k(\widehat{\fr{g}}_{\tilde{B}})$, by Corollary \ref{corollary:zeng_categorical_OPE}:
    \begin{equation}
        [Y(\eta,z),Y(\fr{g}t^{-1},w)]=0
    \end{equation}
    As $\fr{g}t^{-1}$ generates the spaces of states, the usual argument shows that, for any finite morphism $\xi:X\rightarrow V^k(\widehat{\fr{g}}_{\tilde{B}})$, $[Y(\eta,z),Y(\xi,w)]=0$.
    Hence, $\fr{z}_k(\widehat{\fr{g}}_{\tilde{B}})$ is a commutative vertex algebra.
    By extending $R$ to $R\otimes_{\F}(\F[\epsilon]/(\epsilon^2))$ we may apply the standard argument \cite[Proposition 16.24]{Frenkel_Ben_Zvi} from which evaluating $\epsilon=0$ gives a Poisson vertex algebra structure on the original center $\fr{z}_k(\widehat{\fr{g}}_{\tilde{B}})$.
    For the base change statement, we first verify that $f^*(\fr{z}_k(\widehat{\fr{g}}_{\tilde{B}}))$ lands in the desired target.
    By Proposition \ref{prop:vertex_algebra_functor_compatibility} we know that $f^*$ sends $\fr{z}_k(\widehat{\fr{g}}_{\tilde{B}})$ to $\mathrm{Hom}_{\cC}(\mathbf{1}_{\cC},V^{f(k)}(f^*(\widehat{\fr{g}})_{f^*(\tilde{B})}))\otimes_R S$.
    Furthermore, the action of $\fr{g}[t]$ will become the action of $f^*(\fr{g})[t]$ and so indeed we have a map $\ov{f}:f^*(\fr{z}_k(\widehat{\fr{g}}_{\tilde{B}}))\rightarrow  \fr{z}_{f(k)}(f^*(\widehat{\fr{g}})_{f^*(\tilde{B})})$.
    Next, it suffices to show this is a Poisson vertex algebra morphism.
    But this just follows from Proposition \ref{prop:vertex_algebra_functor_compatibility} and the fact that it is $\F$-linear so it preserves the ring of dual numbers.
\end{proof}

As we are only concerned with the center at the critical level, we adopt the following abbreviation:
\begin{equation}
    \fr{z}(\Lulhat{g}{\alpha}):=\fr{z}_{-h^{\vee}}(\Lulhat{g}{\alpha})
\end{equation}
Recall, by interpolating twice the usual trace form the critical numbers are:
\begin{equation}
      -h^{\vee}_{\Lul{so}{\alpha}}=-\alpha+2, \  -h^{\vee}_{\Lul{sp}{\alpha}}=-\alpha-2, \ -h^{\vee}_{\Lul{gl}{\alpha}}=-2\alpha
\end{equation}

\begin{notation}

Let $G_{\alpha}=\GL_{\alpha}$, $\Or_{\alpha}$, or $\Sp_{\alpha}$. In the case of the functor $\ssF{n}:\irep(G_n,\C)\rightarrow \mathrm{Rep}(G_n,\C)$ we denote the induced morphism:
\begin{equation}
    (\overline{\mathrm{ss}_n})_{-h^{\vee}_{\Lul{g}{n}}}:\fr{z}(\Lulhat{g}{n})\rightarrow \fr{z}(\widehat{\fr{g}}_n)
\end{equation}
by:
\begin{equation}
\ssV{n}:\fr{z}(\Lulhat{g}{n})\rightarrow \fr{z}(\widehat{\fr{g}}_n)
\end{equation}
\end{notation}

We call a Poisson vertex algebra with the same differential algebra structure but with $\lambda$-bracket multiplied by $-1$ the opposite Poisson vertex algebra. An anti-isomorphism of Poisson vertex algebras means an isomorphism of differential algebras which identifies the source with the opposite Poisson vertex algebra of the target. See Section~\ref{section:PVA_background}. With this convention, the parity functor induces the following duality.

\begin{prop}
For $\alpha\in R^{\times}$ the parity functor $\Pi_{\alpha}:\irep(\GL_{\alpha},R)\rightarrow \irep(\GL_{-\alpha},R)^{Tw}$ induces an anti-isomorphism of Poisson Vertex Algebras over $R$:
\begin{equation}
\Pi_{\alpha}:\fr{z}(\Lulhat{gl}{\alpha})\rightarrow \fr{z}(\Lulhat{gl}{-\alpha})
\end{equation}
Similarly, for $\alpha\in R$ the parity functor $\Pi_{\alpha}:\irep(\Sp_{\alpha},R)\rightarrow \irep(\Or_{-\alpha},R)^{Tw}$ induces an anti-isomorphism of Poisson vertex algebras over $R$:
\begin{equation}
    \Pi_{\alpha}:\fr{z}(\Lulhat{sp}{\alpha})\rightarrow \fr{z}(\Lulhat{so}{-\alpha})
\end{equation}
\label{prop:interpolated_parity_functor_anti}
\end{prop}

\begin{proof}
It follows from Equations \eqref{eqn:VOA_isomorphism_Type_BC}, \eqref{eqn:VOA_isomorphism_Type_A}, that we have an isomorphism of differential algebras.
To see that it reverses the $\lambda$-bracket, notice that we are changing the level by a sign.
Therefore, deforming by the level will have the $\lambda$-bracket switch its sign.
\end{proof}

We have shown through categorical means that there is a natural Poisson vertex algebra structure on $\fr{z}(\Lulhat{g}{\alpha})$.
Equivalently, one can define the center for $\alpha\in \C$ through ultraproducts.

\begin{prop}

     Let $G_{\alpha}=\GL_{\alpha}$, $\Or_{\alpha}$, or  $\Sp_{\alpha}$ with Lie algebras $\fr{g}_{\alpha}$. Assume in the $\GL$ case that $\alpha$ is non-zero, and let $(\alpha_n)_{n\in \Z_+}$ be a sequence as in Theorem \ref{thm:ultraproduct_deligne}.
     The $\overline{\mathbb{F}}_{p}$-Poisson vertex algebra structure of $\fr{z}^{k_n}(\widehat{\fr{g}}_{\alpha_n})$ induces a $\C$-Poisson vertex algebra structure on $\prod_{\mathcal{U}}^{B}\fr{z}^{k_n}(\widehat{\fr{g}}_{\alpha_n})$ and furthermore the isomorphism in Proposition \ref{prop:ultraproduct_vacuum_module_vertex_algebra} induces an isomorphism of $\C$-Poisson vertex algebras:
     \begin{equation}
         \prod_{\mathcal{U}}^{B}\fr{z}^{k_n}(\widehat{\fr{g}}_{\alpha_n})\cong \fr{z}^k(\Lulhat{g}{\alpha})
     \end{equation}
     where $\alpha,k$ are the images of $[(k_n)],[(\alpha_n)]$ under an isomorphism $\prod_{\mathcal{U}}\overline{\mathbb{F}}_{p_n}\cong \C$.
     \label{prop:ultraproduct_center_vertex_algebra}
\end{prop}

\begin{proof}
    This follows from the observation that the isomorphism in Proposition \ref{prop:ultraproduct_vacuum_module_vertex_algebra} can be upgraded to an isomorphism over $\C[\epsilon]/(\epsilon^2)$ by base change to $\Z[\epsilon]/(\epsilon^2)$.
\end{proof}

To construct explicit generators, we separate the problem into two steps: the finite-rank input and the complex-rank interpolation. The next subsection recalls Molev's formulas in the classical setting. Subsection~\ref{section:interpolating_molevs_invariants} then explains how these formulas can be interpolated in the categories $\irep(G_{\alpha},R)$.

\subsection{Segal--Sugawara Vectors of Classical Lie Algebras}
\label{subsection:molevs_invariants_in_finite_rank}

To obtain explicit generators for the interpolated centers, we first recall the finite-rank Segal--Sugawara families constructed by Molev \cite{Molev_Type_B_C_Invaraints}.  These formulas provide the model that will later be transported to complex rank.

\begin{definition}[Complete Set of Segal--Sugawara Vectors]
Let $\fr{g}$ be a reductive Lie algebra of rank $n$, and recall that $\fr{z}(\widehat{\fr{g}})$ has a natural derivation $\mathsf{T}$ acting on it.
A complete set of Segal--Sugawara vectors is a set of homogeneous vectors $S_1,\cdots, S_n\in \fr{z}(\widehat{\fr{g}})$ such that the set $\{\mathsf{T}^iS_k: 1\leq k\leq n, i\in \Z_{\geq 0}\}$ is algebraically independent and generates $\fr{z}(\widehat{\fr{g}})$ as a $\C$-algebra. Equivalently, $S_1,\ldots,S_n$ freely generate $\fr{z}(\widehat{\fr{g}})$ as a differential algebra.
\end{definition}

The construction to be interpolated in Type $A$ is defined by taking a sum of symmetrized $\mu$-minors or symmetrized $\mu$-permanents over all partitions of $m\in \Z_{\geq 0}$. More precisely, denote a partition $\mu_1\geq \cdots \geq \mu_{\ell}>0$ of $m\in \Z_{\geq 0}$ by a tuple $\mu:=(\mu_1,\cdots,\mu_{\ell})$. The natural number $\ell$ is called the length of $\mu$ and is denoted by $\ell(\mu)$. When the partition is fixed, we often write simply $\ell=\ell(\mu)$. The symmetrized $\mu$-minors and symmetrized $\mu$-permanents from \cite{Molev_Type_B_C_Invaraints} in Type $A$ are:
\begin{align}
    D^A(\mu)&:=\frac{1}{\ell!}\sum_{\sigma\in S_{\ell}}\sum_{i_1,\cdots, i_{\ell}=1}^n \mathrm{sgn}(\sigma)\cdot E_{i_{\sigma(1)},i_1}t^{-\mu_1}\cdots E_{i_{\sigma(\ell)},i_{\ell}}t^{-\mu_{\ell}}\in U(\fr{gl}_n[t^{-1}]t^{-1})
    \label{eqn:D_immanant_type_A}
\\
    P^{A}(\mu)&:=\frac{1}{\ell!}\sum_{\sigma\in S_{\ell}}\sum_{i_1,\cdots, i_{\ell}=1}^n E_{i_{\sigma(1)},i_1}t^{-\mu_1}\cdots E_{i_{\sigma(\ell)},i_{\ell}}t^{-\mu_{\ell}}\in U(\fr{gl}_n[t^{-1}]t^{-1})
    \label{eqn:P_immanant_type_A}
\end{align}
where $E_{i,j}$ denotes the standard elementary matrices.

To present Molev's construction in a more uniform manner, we introduce the following notation.
Let $R$ be a commutative $\mathbb{Q}$-algebra and let $x\in R$. For $n\in \Z_{+}$ and $\ell\leq m\in \Z_+$ set:
\begin{equation}
    (x)_n:=\prod_{k=0}^{n-1}(x-k) \qquad  \binom{x}{n}:=\frac{(x)_n}{n!} \qquad   Q_{m,\ell}(x):=\frac{\ell!}{m!}\prod_{k=\ell}^{m-1}(x-k)
\label{eqn:binomical_coefficient_over_ring}
\end{equation}
with the convention that $Q_{m,m}(x)=1$. Notice that if $\binom{x}{\ell}$ is invertible, then $Q_{m,\ell}(x)=\binom{x}{m}\binom{x}{\ell}^{-1}$.

If $\mu\vdash m$, denote the number of permutations in $S_m$ of cycle type $\mu$ by $c_{\mu,m}$.
We may now recall Molev's construction in Type $A$: 

\begin{theorem}[\cite{Molev_Type_B_C_Invaraints}, Theorem 2.1]

For $1\leq m\leq n$ the elements:
\begin{equation}
    \phi_{m,n}:=\sum_{\mu\vdash m} c_{\mu,m}\cdot Q_{m,\ell}(n)\cdot D^A(\mu)
    \label{eqn:molev_finite_rank_type_a_invariant_anti}
\end{equation}
\begin{equation}
    \psi_{m,n}:= \sum_{\mu\vdash m} c_{\mu,m}\cdot (-1)^{\ell+m}Q_{m,\ell}(-n)\cdot P^{A}(\mu)
    \label{eqn:molev_finite_rank_type_a_invariant_sym}
\end{equation}
lie in the center $\fr{z}(\widehat{\fr{gl}}_n)$.
Furthermore, $\{\phi_{1,n},\cdots, \phi_{n,n}\}, \{\psi_{1,n},\cdots, \psi_{n,n}\}$ are complete sets of Segal--Sugawara vectors.

\label{theorem:finite_rank_Molev_type_A_invariant}
\end{theorem}

\begin{proof}
We only need to check that:
    \begin{equation}
        (-1)^{\ell+m}Q_{m,\ell}(-n)=\binom{n+m-1}{m}\binom{n+\ell-1}{\ell}^{-1}, \ \ \   Q_{m,\ell}(n)=\binom{n}{m}\binom{n}{\ell}^{-1}
    \end{equation} 
    so that Equations \eqref{eqn:molev_finite_rank_type_a_invariant_anti} and \eqref{eqn:molev_finite_rank_type_a_invariant_sym} agree with the formulation given in the proof of \cite[Theorem 2.1]{Molev_Type_B_C_Invaraints}.
    This follows from the basic properties of falling and raising factorials.
\end{proof}

To describe the effect of the interpolated parity functor on the Type $A$ generators, we need to record the effect the Cartan involution has in the finite-rank case. 
Recall that the Cartan involution of $\fr{gl}_n$ is the Lie algebra involution defined $E_{i,j}\mapsto -E_{j,i}$. This naturally induces an involution on the center:
\begin{equation}
    \nu_n:\fr{z}(\widehat{\fr{gl}}_n)\rightarrow \fr{z}(\widehat{\fr{gl}}_n)
\end{equation}
which we also denote by $\nu_n$.
One checks, see \cite{Molev} for details, that its action on Molev's higher Segal--Sugawara vectors is: 
\begin{align}
    \nu_n(\phi_{m,n}) &= \sum_{\mu\vdash m }c_{\mu,m}\cdot Q_{m,\ell}(n)(-1)^{\ell}D^A(\mu) \label{equation:cartan_fintie_rank_calc_1} \\
    \nu_n(\psi_{m,n}) &= \sum_{\mu\vdash m }c_{\mu,m}\cdot Q_{m,\ell}(-n)(-1)^{m}P^A(\mu) \label{equation:cartan_fintie_rank_calc_2}
\end{align}

The constructions to be interpolated in Types $B$ and $C$ are similar to the Type $A$ constructions. In these cases, however, there is only one complete set of Segal--Sugawara vectors in each type. The Type $B$ family is formed from sums of symmetrized $\mu$-permanents, while the Type $C$ family is formed from sums of symmetrized $\mu$-minors. This parallel will later be explained by the interpolated parity functor relating Types $B$ and $C$.

We define modified generators $F_{ij}$ that incorporate the structure of the underlying symmetric or anti-symmetric bilinear form associated with the orthogonal and symplectic Lie algebras, respectively.
In the symmetric setting, we define the corresponding isomorphism on the standard basis $e_i$ and dual basis $e_i^*$ by:
 \begin{equation}
 \theta:\C^{2n+1}\rightarrow (\C^{2n+1})^*, \ \theta(e_i):=e_{2n-i+2}^*
 \label{eqn:symmetric_isomorphism}
 \end{equation} 
In the skew-symmetric setting, define the corresponding isomorphism by:  
 \begin{equation}
     \theta:\C^{2n}\rightarrow (\C^{2n})^*, \ \theta(e_i):=\begin{cases}
         e_{2n-i+1}^* & 1\leq i\leq n \\
         -e_{2n-i+1}^* & n<i\leq 2n
     \end{cases}
     \label{eqn:anti_symmetric_isomorphism}
 \end{equation}
 In both cases, we set: 
\begin{equation}
F_{ij}:= E_{i,j}-(\theta^{-1}\otimes \theta)(e_j^*\otimes e_i)
\end{equation}
Notice that we are using the identification $\End(V)\cong V\otimes V^*$ for the relevant standard representation $V$. 

The symmetrized $\mu$-permanent in Type $B$ and symmetrized $\mu$-minor in Type $C$ are respectively:

\begin{align}
    P^{B}(\mu) &:= \frac{1}{\ell!}\sum_{\sigma\in S_{\ell}}\sum_{i_1,\cdots, i_{\ell}=1}^{2n+1} F_{i_{\sigma(1)},i_1}t^{-\mu_1}\cdots F_{i_{\sigma(\ell)},i_{\ell}}t^{-\mu_{\ell}} \in U(\fr{so}_{2n+1}[t^{-1}]t^{-1})
    \label{eqn:P_immanant_type_B}
\\
    D^{C}(\mu) &:= \frac{1}{\ell!}\sum_{\sigma\in S_{\ell}}\sum_{i_1,\cdots, i_{\ell}=1}^{2n} \mathrm{sgn}(\sigma)\cdot F_{i_{\sigma(1)},i_1}t^{-\mu_1}\cdots F_{i_{\sigma(\ell)},i_{\ell}}t^{-\mu_{\ell}} \in U(\fr{sp}_{2n}[t^{-1}]t^{-1})
    \label{eqn:D_immanant_type_C}
\end{align}
The higher Segal--Sugawara vectors in Types $B$ and $C$ are:

\begin{theorem}[\cite{Molev_Type_B_C_Invaraints}, Theorem 2.3]
For even $m$ with $2\leq m\leq 2n$, the element:
\begin{equation}
\phi_{m,2n+1}^{B}:=\sum_{\substack{\mu\vdash m \\ \ell(\mu) \text{ even}}}c_{\mu,m}\cdot Q_{m,\ell}(-2n)\cdot P^B(\mu)
\label{eqn:molev_finite_rank_type_b_invariant}
\end{equation}
belongs to the center $\fr{z}(\widehat{\fr{so}}_{2n+1})$.
For the same values of $m$, the element: 
\begin{equation}
\phi_{m,2n}^{C}:=\sum_{\substack{\mu \vdash m, \\ \ell(\mu)\text{ even}}} c_{\mu,m}\cdot Q_{m,\ell}(2n+1)\cdot D^C(\mu)
\label{eqn:molev_finite_rank_type_c_invariant}
\end{equation}
belongs to the center $\fr{z}(\widehat{\fr{sp}}_{2n})$.
Furthermore, the sets $\{\phi_{2,2n+1}^B,\cdots, \phi^B_{2n,2n+1}\}$ and $\{\phi^C_{2,2n},\cdots, \phi_{2n,2n}^C\}$ are complete sets of Segal--Sugawara vectors in their respective centers.
\end{theorem}

\begin{proof}
    We only need to check that for even $\ell,m$: 
    \begin{equation}
        Q_{m,\ell}(2n+1)=\binom{2n+1}{m}\binom{2n+1}{\ell}^{-1}, \ \  Q_{m,\ell}(-2n)=\binom{(2n+1)+m-2}{m}\binom{(2n+1)+\ell-2}{\ell}^{-1}
    \end{equation}
    so that Equations \eqref{eqn:molev_finite_rank_type_c_invariant} and \eqref{eqn:molev_finite_rank_type_b_invariant} agree with the formulations given in the proof of \cite[Theorem 2.3]{Molev_Type_B_C_Invaraints}.
    This follows from the observation that for even $r$:
    \begin{equation}
        \binom{T}{r} = \binom{-T+r-1}{r}
    \end{equation}
\end{proof}

\subsection{Interpolating Segal--Sugawara Vectors}
\label{section:interpolating_molevs_invariants}
  We now transport the finite-rank constructions from the previous subsection into interpolating categories and use them to build explicit generators of the interpolated Feigin--Frenkel center. Let $G_{\alpha}$ denote $\GL_{\alpha}$, $\Or_{\alpha}$, or $\Sp_{\alpha}$.

\begin{definition}[Complete Set of Interpolated Segal--Sugawara Vectors]
Consider the interpolating category $\irep(G_{\alpha},R)$ with corresponding Lie algebra $\Lul{g}{\alpha}$.
A complete set of interpolated Segal--Sugawara vectors is a subset $\{S_{m,\alpha}:m\in \Z_+\}\subset \fr{z}(\Lulhat{g}{\alpha})$ such that the set $\{\mathsf{T}^iS_{m,\alpha}:i\in \Z_{\geq 0}, m\in \Z_+\}$, where $\mathsf{T}$ denotes the translation operator, is algebraically independent over $R$ and generates $\fr{z}(\Lulhat{g}{\alpha})$ as an $R$-algebra. In other words, $\{S_{m,\alpha}:m\in \Z_+\}$ generates $\fr{z}(\Lulhat{g}{\alpha})$ as a differential $R$-algebra.
\end{definition}

In the last subsection, we recalled that Molev constructs higher Segal--Sugawara vectors by first introducing symmetrized $\mu$-minors and $\mu$-permanents, and then summing these expressions over all partitions of $m$. The first step in interpolating this construction is therefore to rewrite the symmetrized $\mu$-minors and $\mu$-permanents categorically. 

We begin with the simplest example. Consider the quadratic Casimir element of $\fr{gl}_n$:
\begin{equation}
    \sum_{i,j=1}^n E_{i,j}\otimes E_{j,i}
\end{equation}
This can be constructed categorically. Let $V:=\C^n$,  and denote the standard basis of $V$ by $e_1,\cdots, e_n$ and the dual basis of $V^*$ by $e_1^*,\cdots, e^*_n$. Under the identification $\mathrm{End}(V)\cong V\otimes V^*$ we have $E_{i,j}=e_i\otimes e_j^*$. Hence:
\begin{equation}
    (\cv_V\otimes \cv_V)(1)=\sum_{i_1,i_2=1}^n e_{i_1}\otimes e_{i_1}^*\otimes e_{i_2}\otimes e_{i_2}^*
\end{equation}
Applying the braiding that switches the first and third tensor factors gives:
\begin{equation}
    \sum_{i_1,i_2=1}^ne_{i_2}\otimes e_{i_1}^*\otimes e_{i_1}\otimes e_{i_2}^*=\sum_{i_1,i_2=1}^nE_{i_2,i_1}\otimes E_{i_1,i_2}
\end{equation}
Thus the usual summation over matrix indices is obtained from coevaluations, braidings, and the multiplication in the universal enveloping algebra. 

More generally, for $\ell\in \Z_+$, define an embedding $S_{\ell}\hookrightarrow S_{2\ell}$ by sending $\sigma\in S_{\ell}$ to $\widetilde{\sigma}\in S_{2\ell}$, where:
\begin{equation}
     \tilde{\sigma}(k):=\begin{cases} 2\sigma(\frac{k+1}{2})-1 & k \text{ is odd}\\ k & k\text{ is even}
    \end{cases}
\end{equation}
 This embedding lets $\sigma\in S_{\ell}$ act on $(V\otimes V^*)^{\otimes \ell}$ by permuting the $V$-factors while leaving the $V^*$-factors fixed. With this convention, one checks that:
 \begin{equation}
     (\tilde{\sigma}\circ (\mathrm{coev}_V)^{\otimes \ell})(1)=\sum_{i_1,\cdots, i_{\ell}=1}^{n} E_{i_{\sigma(1)},i_1}\otimes \cdots \otimes E_{i_{\sigma(\ell)},i_{\ell}}
 \end{equation}
If $\pi_{\fr{gl}_n}:T(\fr{gl}_n)\rightarrow U(\fr{gl}_n)$ denotes the multiplication map induced by the associative algebra $U(\fr{gl}_n)$ we see that:
 \begin{equation}
     \pi_{\fr{gl}_n}\circ (\tilde{\sigma}\circ (\mathrm{coev}_V)^{\otimes \ell})(1)=\sum_{i_1,\cdots, i_{\ell}=1}^{n} E_{i_{\sigma(1)},i_1}\cdots E_{i_{\sigma(\ell)},i_{\ell}}
 \end{equation}

To obtain Molev's symmetrized $\mu$-minors and permanents we now only need to insert loop degrees, which can be accomplished by adding some extra notation. For $r\in \Z_{\geq 0}$ let:
\begin{equation}
    \cv_{\bullet}[-r]\in  \Hom(1,(\bullet\otimes \bullet^*)t^{-r})
    \end{equation}
denote the coevaluation morphism placed in loop degree $-r$. For a tuple of positive integers $\mu=(\mu_1,\ldots, \mu_{\ell})$ set:
\begin{equation}
(\cv_{\bullet})^{\otimes \mu}:= (\cv_{\bullet}[-\mu_1])\otimes \cdots \otimes (\cv_{\bullet}[-\mu_{\ell}])
\end{equation}

The next lemma allows us to obtain the symmetrized $\mu$-minors and permanents categorically. For readability, we denote the semisimplification functor at $n\in \Z_{+}$ for the corresponding classical group $G$ by $\mathrm{ss}_n^G$, and $\pi_{\fr{g}}:T(\fr{g})\rightarrow U(\fr{g})$ for the multiplication map induced by the associative algebra $U(\fr{g})$.
\begin{lemma}
    Let $\mu\vdash m$, recall that $\ell:=\ell(\mu)$, and let $\sigma\in S_{\ell}$.
    The following holds in the respective representation categories:
    \begin{equation}
       \pi_{\fr{gl}_n}\circ \ssF{n}^{\mathrm{GL}}(\tilde{\sigma}\circ (\cv_{\bullet})^{\otimes \mu})=\sum_{i_1,\cdots, i_{\ell}=1}^n E_{i_{\sigma(1)},i_1}t^{-\mu_1}\cdots E_{i_{\sigma(\ell)},i_{\ell}}t^{-\mu_{\ell}} \tag{Type A}
    \end{equation}
      \begin{equation}
          \pi_{\fr{so}_{2n+1}}\circ (\Id\otimes \theta)^{\otimes \ell}\circ \ssF{2n+1}^{\mathrm{O}}((h_{\bullet,2})^{\otimes \ell}\circ\tilde{\sigma}\circ ( \cv_{\bullet})^{\otimes \mu})=\sum_{i_1,\cdots, i_{\ell}=1}^{2n+1} F_{i_{\sigma(1)},i_1}t^{-\mu_1}\cdots F_{i_{\sigma(\ell)},i_{\ell}}t^{-\mu_{\ell}} \tag{Type B}
    \end{equation}
    \begin{equation}
        \pi_{\fr{sp}_{2n}}\circ (\Id\otimes \theta)^{\otimes \ell}\circ \ssF{2n}^{\mathrm{Sp}}((e_{\bullet,2})^{\otimes \ell}\circ \tilde{\sigma}\circ ( \cv_{\bullet})^{\otimes \mu})=\sum_{i_1,\cdots, i_{\ell}=1}^{2n} F_{i_{\sigma(1)},i_1}t^{-\mu_1}\cdots F_{i_{\sigma(\ell)},i_{\ell}}t^{-\mu_{\ell}} \tag{Type C}
    \end{equation}
    Here $\theta:V\rightarrow V^*$ denotes the symmetric or skew-symmetric isomorphism defined in Equations \eqref{eqn:symmetric_isomorphism} and \eqref{eqn:anti_symmetric_isomorphism}, respectively.
\end{lemma}

\begin{proof}
Type A follows from the discussion preceding the lemma.
Expanding the left-hand side of the Type $B$ case, and recalling Equation \ref{eqn:theta_coevaluation}, we obtain:
\begin{align*}
   &(\Id\otimes \theta)^{\otimes \ell}\circ (h_{V,2})^{\otimes \ell}\circ (\Id\otimes \theta^{-1})^{\otimes \ell}\circ  \tilde{\sigma}\circ \left(\sum_{i_1,\cdots, i_\ell=1}^{2n+1} e_{i_1}\otimes e_{i_1}^*t^{-\mu_1}\otimes \cdots \otimes e_{i_{\ell}}\otimes e_{i_\ell}^*t^{-\mu_\ell}\right) \\ 
   &= (\Id-c_{V^*,V}\circ\theta\otimes \theta^{-1})^{\otimes\ell}\circ \tilde{\sigma}\circ \left(\sum_{i_1,\cdots, i_\ell=1}^{2n+1} e_{i_1}\otimes e_{i_1}^*t^{-\mu_1}\otimes \cdots \otimes e_{i_\ell}\otimes e_{i_\ell}^*t^{-\mu_\ell}\right) \\
   &= (\Id-c_{V^*,V}\circ\theta\otimes \theta^{-1})^{\otimes\ell}\circ \left(\sum_{i_1,\cdots, i_\ell=1}^{2n+1}E_{i_{\sigma(1)},i_{1}}[-\mu_{1}]\otimes \cdots \otimes E_{i_{\sigma(\ell)},i_{\ell}}[-\mu_{\ell}]\right) 
\end{align*}
Now notice that:
\begin{align}
    (\Id-c_{V^*,V}\circ\theta\otimes \theta^{-1})\circ E_{i,j} &= E_{i,j}-((\theta^{-1}\otimes \theta)\circ c_{V,V^*})(e_i\otimes e_j^*)\\
    &= E_{i,j}-(\theta^{-1}\otimes \theta)(e_{j}^*\otimes e_i) = F_{i,j}
\end{align}
The Type $C$ argument follows in the same way, except a sign is introduced due to the braiding.
See \cite[Section 1.5]{Molev} for more details.
\end{proof}
Denote the following elements of $\C[S_{\ell}]$ by:
 \begin{equation}
    A^{(\ell)}:=\frac{1}{\ell!}\sum_{\sigma\in S_{\ell}}\mathrm{sgn}(\sigma)\cdot \sigma \qquad H^{(\ell)}:=\frac{1}{\ell!}\sum_{\sigma\in S_{\ell}}\sigma
 \end{equation}
 We extend $\sigma\mapsto \widetilde{\sigma}$ linearly to 
$\C[S_{\ell}]\rightarrow \C[S_{2\ell}]$, and denote the images of $A^{(\ell)}$ and $H^{(\ell)}$ by $\widetilde{A^{(\ell)}}$ and
$\widetilde{H^{(\ell)}}$.
\begin{corollary}
    The following equalities hold in their respective representation categories for $n\in \Z_+$:
    \begin{equation}
        D^A(\mu)=\pi_{\fr{gl}_n}\circ \ssF{n}^{\mathrm{GL}}(\widetilde{A^{(\ell)}}\circ (\cv_{\bullet})^{\otimes \mu})
    \end{equation}
    \begin{equation}
        P^A(\mu)=\pi_{\fr{gl}_n}\circ \ssF{n}^{\mathrm{GL}}(\widetilde{H^{(\ell)}}\circ (\cv_{\bullet})^{\otimes \mu})
    \end{equation}
    \begin{equation}
        P^B(\mu)=  \pi_{\fr{so}_{2n+1}}\circ (\Id\otimes \theta)^{\otimes \ell}\circ \ssF{2n+1}^{\mathrm{O}}((h_{\bullet,2})^{\otimes \ell}\circ\widetilde{H^{(\ell)}}\circ ( \cv_{\bullet})^{\otimes \mu})
    \end{equation}
    \begin{equation}
        D^C(\mu)=\pi_{\fr{sp}_{2n}}\circ (\Id\otimes \theta)^{\otimes \ell}\circ \ssF{2n}^{\mathrm{Sp}}((e_{\bullet,2})^{\otimes \ell}\circ \widetilde{A^{(\ell)}}\circ ( \cv_{\bullet})^{\otimes \mu})
    \end{equation}
\end{corollary}

The preceding corollary shows that Molev's symmetrized minors and permanents are obtained from categorical morphisms in interpolating categories after applying the finite-rank semisimplification functors. We therefore define the interpolated higher Segal--Sugawara vectors in Type $A$ by:
\begin{equation}
    \phi_{m,T}:=\pi\left(\sum_{\mu\vdash m} Q_{m,\ell}(T)\cdot c_{\mu,m}\cdot \widetilde{A^{(\ell)}}\circ (\cv_{\bullet})^{\otimes \mu}\right)
\end{equation}
\begin{equation}
    \psi_{m,T}:=\pi\left(\sum_{\mu\vdash m }(-1)^{\ell+m}\cdot Q_{m,\ell}(-T)\cdot c_{\mu,m}\cdot \widetilde{H^{(\ell)}}\circ (\cv_{\bullet})^{\otimes \mu}\right)
\end{equation}
Here $\pi:\mathcal{T}(\Lul{gl}{T}[t^{-1}]t^{-1})\rightarrow U(\Lul{gl}{T}[t^{-1}]t^{-1})$ is the multiplication map induced by the universal enveloping algebra structure. For Type $A$ we set $R=\C[T,T^{-1}]$.

\begin{theorem}[Interpolated Segal--Sugawara Vectors: Type $A$]
\label{theorem:interpolate_segal_sugawara_type_A} 
The Type $A$ families $\{\phi_{m,T}: m\in \Z_+\}$ and $\{\psi_{m,T}: m\in \Z_+\}$ satisfy:  
\begin{enumerate}
    \item \textbf{Specialization:} For $\alpha\in \C^{\times}$ let $\mathrm{Ev}_{T=\alpha}:f_{\alpha}^*(\irep(\GL_{T},R))\rightarrow \irep(\GL_{\alpha},\C)$ denote the change of base functor induced by the evaluation morphism $f_{\alpha}:R\rightarrow \C$, where $p(T,T^{-1})\mapsto p(\alpha,\alpha^{-1})$. 
    By defining:
    \begin{equation} \phi_{m,T=\alpha}:=\mathrm{Ev}_{T=\alpha}(\phi_{m,T})\qquad \psi_{m,T=\alpha}:=\mathrm{Ev}_{T=\alpha}(\psi_{m,T})\end{equation}
    the specialized sets $\{\phi_{m,T=\alpha}\mid m\in \Z_+\}$ and $\{\psi_{m,T=\alpha}\mid m\in \Z_+\}$ are complete sets of interpolated Segal--Sugawara vectors of $\fr{z}(\Lulhat{gl}{\alpha})$. In particular, $\mathrm{Ev}_{T=\alpha}(\fr{z}(\Lulhat{gl}{T}))\cong \fr{z}(\Lulhat{gl}{\alpha})$ as $\C$-Poisson vertex algebras.
    \item \textbf{Completeness:} The sets $\{\phi_{m,T}:m\in \Z_+\}$ and
    $\{\psi_{m,T}:m\in \Z_+\}$ form complete sets of interpolated
    Segal--Sugawara vectors for $\fr{z}(\Lulhat{gl}{T})$.
    \item \textbf{Compatibility with Finite Rank:} Under the mapping $\ssV{n}$, the specialized generators correspond to Molev's classical invariants:
    \begin{equation}
    \label{eqn:prop_type_A_ss_n_Gen_1}
        \ssV{n}(\phi_{i,T=n})=
        \begin{cases}
            \phi_{i,n} & \text{ if } 1\leq i\leq n \\
            0 & \text{ otherwise}
        \end{cases}
    \end{equation}
    and for $1\leq i\leq n$:
    \begin{equation}
    \label{eqn:prop_type_A_ss_n_Gen_2}
         \ssV{n}(\psi_{i,T=n}) = \psi_{i,n}
    \end{equation}
    \item \textbf{Asymptotic $\lambda$-brackets:} Let $\ssIV{n}:\fr{z}(\widehat{\fr{gl}}_n)\rightarrow \fr{z}(\Lulhat{gl}{n})$ be the section of differential algebras given by $\phi_{i,n}\mapsto \phi_{i,T=n}$ for $1\leq i\leq n$. For fixed positive integers $i_0,j_0$ and sufficiently large $n$, the $\lambda$-bracket respects the section:
    \begin{equation}
        \ssIV{n}(\{ \phi_{i_0,n} \ {}_{\lambda} \ \phi_{j_0,n} \}) = \{ \phi_{i_0,T=n} \ {}_{\lambda} \ \phi_{j_0,T=n} \}
    \end{equation}
\end{enumerate}
\end{theorem}

The proof of the completeness and specialization statements has three steps. First, we specialize to sufficiently large positive integer ranks and apply the Interpolation Principle to show that the families $\{\phi_{m,T}:m\in \Z_+\}$ and $\{\psi_{m,T}:m\in \Z_+\}$ are $\Lul{gl}{T}[t]$-invariant and freely generate differential subalgebras over $R$. Second, we use the ultraproduct realization of $\fr{z}(\Lulhat{gl}{\alpha})$, to show that for every $\alpha\in \C^\times$ the specialized families $\{\phi_{m,T=\alpha}:m\in \Z_+\}$ and $\{\psi_{m,T=\alpha}:m\in \Z_+\}$ are complete sets of Segal--Sugawara vectors in $\fr{z}(\Lulhat{gl}{\alpha})$. Third, we use an exhaustive filtration of $\fr{z}(\Lulhat{gl}{T})$, together with finite generation and Nakayama's lemma, to lift this generation statement from all specializations to the generic $R$-family. This proves that $\{\phi_{m,T}:m\in \Z_+\}$ and $\{\psi_{m,T}:m\in \Z_+\}$ form complete sets of interpolated Segal--Sugawara vectors.

To carry out these steps, we need two technical lemmas. The first concerns two filtrations on the center. For Step 1, we need an exhaustive filtration to which the Interpolation Principle can be applied. This is provided by the PBW filtration
\(F^{\bullet}\fr{z}(\Lulhat{gl}{T})\). For Step 3, we need a filtration whose filtered pieces are objects of the interpolating category and are finite over the base ring. This is provided by the loop filtration
\(\mathscr{F}^{\bullet}\fr{z}(\Lulhat{gl}{T})\). We explain both constructions below.

The second lemma identifies the relevant filtered pieces in the ultraproduct realization of the center. It is the main input for Step 2, where we transfer finite-rank generation statements to the specializations at complex rank.
 
To that end, recall in  $\irep(\GL_{T},R)$ the universal affine vertex algebra $V^k(\Lul{gl}{T})$ is isomorphic as an object to $S(\Lul{gl}{T}[t^{-1}]t^{-1})$.
This realization carries two filtrations needed below.

The first is the PBW filtration, defined for $m\in\Z_{\geq 0}$ by:
\begin{equation}
    F^mV^k(\Lul{gl}{T})
    :=
    \bigoplus_{a=0}^m
    S^a(\Lul{gl}{T}[t^{-1}]t^{-1}).
\end{equation}
The second is the loop filtration, defined for $m\in\Z_{\geq 0}$ by:
\begin{equation}
    \mathscr{F}^mV^k(\Lul{gl}{T})
    :=
    \bigoplus_{\substack{(r_s)_{s\geq 0}\\ \sum_{s\geq 0}(s+1)r_s\leq m}}
    \bigotimes_{s\geq 0} S^{r_s}(\Lul{gl}{T}t^{-s-1}).
\end{equation}
Only finitely many factors in the tensor product are nontrivial in each summand. These filtrations induce filtrations on the center:
\begin{equation}
    F^m\fr{z}(\Lulhat{gl}{T})
    :=
    \Hom_{\Lul{gl}{T}[t]}
    \left(
        \mathbf{1}_{\irep(\GL_T,R)},
        F^mV^{-h^{\vee}}(\Lul{gl}{T})
    \right)
\end{equation}
and
\begin{equation}
    \mathscr{F}^m\fr{z}(\Lulhat{gl}{T})
    :=
    \Hom_{\Lul{gl}{T}[t]}
    \left(
        \mathbf{1}_{\irep(\GL_T,R)},
        \mathscr{F}^mV^{-h^{\vee}}(\Lul{gl}{T})
    \right)
\end{equation}
Similarly, we can define these filtrations on the specializations $\fr{z}(\Lulhat{gl}{\alpha})$ for $\alpha \in \C^{\times}$. 
\begin{lemma}
\label{lemma:center_filtrations_exhaustive}
The induced filtrations $F^{\bullet}\fr{z}(\Lulhat{gl}{T})$ and $\mathscr{F}^{\bullet}\fr{z}(\Lulhat{gl}{T})$ are exhaustive. Moreover, for every $m\in \Z_{\geq 0}$, the loop-filtered piece $\mathscr{F}^m\fr{z}(\Lulhat{gl}{T})$ is a finitely generated $R$-module. Similarly, for every $\alpha\in \C^\times$, the corresponding filtrations $F^{\bullet}\fr{z}(\Lulhat{gl}{\alpha})$ and $\mathscr{F}^{\bullet}\fr{z}(\Lulhat{gl}{\alpha})$ are exhaustive.
\end{lemma}
\begin{proof}
Suppose that $f\in \fr{z}(\Lulhat{gl}{T})$. Since $f$ is a morphism from the monoidal unit into $V^{-h^{\vee}}(\Lul{gl}{T})$, it factors through a finite subobject of $V^{-h^{\vee}}(\Lul{gl}{T})$. Since the filtrations $F^{\bullet}V^{-h^{\vee}}(\Lul{gl}{T})$ and $\mathscr{F}^{\bullet}V^{-h^{\vee}}(\Lul{gl}{T})$ are exhaustive, this finite subobject is contained in some filtered piece. Hence $f$ lies in $F^m\fr{z}(\Lulhat{gl}{T})$ and in $\mathscr{F}^{m'}\fr{z}(\Lulhat{gl}{T})$ for some $m,m'\in \Z_{\geq 0}$. This proves that the induced filtrations on the center are exhaustive. The same argument applies after specialization to $\alpha\in \C^{\times}$.

For the finite generation claim, note that $\mathscr{F}^m\fr{z}(\Lulhat{gl}{T})$ is an $R$-submodule of $\Hom_{\irep(\GL_T)} \left(\mathbf{1}_{\irep(\GL_T)},\mathscr{F}^mV^{-h^{\vee}}(\Lul{gl}{T})\right)$
which is a free $R$-module of finite rank. Since $R=\C[T,T^{-1}]$ is a principal ideal domain, every submodule of a finite-rank free $R$-module is finitely generated. Therefore $\mathscr{F}^m\fr{z}(\Lulhat{gl}{T})$ is finitely generated over $R$.
\end{proof}
To utilize the ultraproduct realization of $\fr{z}(\Lulhat{gl}{n_0})$ for $n_0\in \Z_+$, we prove an asymptotic version of \cite[Theorem 1.1]{Modular_Center}.
Here we use the notation of \cite{Modular_Center} where, if $R$ is a ring, then $V_R:=V^{-h^{\vee}}((\fr{gl}_{n_0+p})_{\Z})\otimes_{\Z} R$.
That is, we have a $\Z$-form for the universal affine vertex algebra at the critical level.
Furthermore, let $P_1,\cdots, P_{n_0+p}\in \C[\fr{gl}^{*}_{n_0+p}]$ denote homogeneous generators and $P_{1,-1},\cdots, P_{n_0+p,-1}$ the lifted generators in $S(\fr{gl}_{n_0+p}[t^{-1}])^{\fr{gl}_{n_0+p}[t]}$.
We refer the reader to \cite{Modular_Center} for details.

\begin{lemma}
    \label{lemma:large_p_filtration_lemma}
    Fix $m_0\in \Z_{\geq 0}$, and let $\mathbf{k}$ be a field of characteristic $p$ such that $p$ is a very good prime for $\mathfrak{gl}_{n_0+p}$.
    Let $Q=\mathbb{Z}[\frac{1}{m_0!}]$. Suppose there is a homomorphism $Q\rightarrow \mathbf{k}$ and that there exist elements $S_{i,-1}\in  V_Q\cap \fr{z}(V_{\C})$ such that $\mathrm{gr}(S_{i,-1})=P_{i,-1}$ for $1\leq i\leq m_0$. If the image of $S_{i,-1}$ under $V_{Q}\rightarrow V_{\mathbf{k}}$ is also denoted $S_{i,-1}$, then the following is true:
    \begin{equation}
        F^{m_0} \fr{z}(V_{\mathbf{k}})\cong F^{m_0}\mathbf{k}[\mathsf{T}^jS_{i,-1}: 1\leq i\leq m_0, j\in \Z_{\geq 0}]
    \end{equation}
\end{lemma}

\begin{proof}
    This follows by first noting that endowing $V^{-h^{\vee}}(\fr{gl}_{n_0+p})$ with the PBW grading yields $\mathrm{gr} (V^{-h^{\vee}}(\fr{gl}_{n_0+p}))\cong S(\fr{gl}_{n_0+p}[t^{-1}])$ as $\mathfrak{gl}_{n_0+p}[t]$-algebras. Furthermore, by \cite[Theorem 1.1]{Modular_Center} we know for very good primes that if $P_{1},\cdots, P_{n_0+p}\in \C[\fr{gl}_{n_0+p}^*]$ are homogeneous generators, then they lift to homogeneous generators $P_{1,-1},\cdots, P_{n_0+p,-1}$ of $S(\fr{gl}_{n_0+p}[t^{-1}])^{\fr{gl}_{n_0+p}[t]}$ as a $\mathbf{k}[\mathsf{T}]$-algebra. 
    In particular, this implies that if $p$ is very large, then the first $m_0$-graded components of $S(\fr{gl}_{n_0+p}[t^{-1}])^{\fr{gl}_{n_0+p}[t]}$ are generated by $P_{1,-1},\cdots, P_{m_0,-1}$, assuming $P_{i,-1}$ has degree $i$.
    On the other hand, notice we have:
    \begin{equation}
        \mathrm{gr}(\fr{z}(\fr{gl}_{n_0+p}))\subset \mathrm{gr}(V^{-h^{\vee}}(\fr{gl}_{n_0+p}))^{\fr{gl}_{n_0+p}[t]}
    \end{equation}
    This hypothesis implies that the first $m_0$-graded components of $\mathrm{gr}(\fr{z}(\fr{gl}_{n_0+p}))$ are generated by $P_{1,-1},\cdots, P_{m_0,-1}$.
    Therefore, the $m_0$-filtered component of the center is generated by their lifts $S_{1,-1},\cdots, S_{m_0,-1}$.
\end{proof}

In particular, we may apply Lemma \ref{lemma:large_p_filtration_lemma} to Molev's generators $\phi_{1,n_0+p},\cdots, \phi_{m,n_0+p}$ and $\psi_{1,n_0+p},\cdots, \psi_{m,n_0+p}$ when $p$ is sufficiently larger than $m$. 

\begin{proof}[Proof of Theorem~\ref{theorem:interpolate_segal_sugawara_type_A}]
We will present the proof for $\{\phi_{m,T}\mid m\in \Z_{+}\}$ as the other case follows similarly. First, we must show that each $\phi_{m,T}$ is $\Lul{gl}{T}[t]$-invariant and that $\{\mathsf{T}^i\phi_{m,T}\mid m\in \Z_{+}, i\in \Z_{\geq 0}\}$ is algebraically independent over $R$.

To that end, fix $m_0\in \Z_+$. For $r\gg2m_0$, Lemma \ref{prop:interpolation_principle_A} implies that $\ssV{r}|_{F^{m_0}}:F^{m_0}\fr{z}(\Lulhat{gl}{r})\cong F^{m_0}\fr{z}(\widehat{\fr{gl}}_{r})$ as $\C$-vector spaces. 
By Corollary \ref{corollary:functor_and_base_change}, the functor $\ssF{r}$ sends $\Lul{gl}{r}[t]$-invariants to $\fr{gl}_{r}[t]$-invariants. This implies 
for $r\gg2(m_0+1)$, $\mathrm{ss}_{r}(f)\in F^{m_0}\fr{z}(\widehat{\fr{gl}}_{r})$ if and only if $f\in F^{m_0}\fr{z}(\Lulhat{gl}{r})$. 
To see this note if $(M,\rho_M),(M',\rho_{M'})$ are $\fr{gl}_{r}[t]$-modules, then $f\in \Hom_{\fr{gl}_{r}[t]}(M,M')$ is equivalent to the following equation holding:
\begin{equation}
    f\circ \rho_M= \rho_{M'}\circ (\Id\otimes f)
\end{equation}
As $r\gg2(m_0+1)$, the Interpolation Principle implies that:
\begin{equation}
    \ssF{r}(\rho_{V^{-2r}(\Lul{gl}{r})}\circ (\Id\otimes f))=\rho_{V^{-2r}(\fr{gl}_{r})}\circ (\Id\otimes \mathrm{ss}_{r}(f))=0\Leftrightarrow \rho_{V^{-2r}(\Lul{gl}{r})}\circ (\Id\otimes f)=0
\end{equation}
With this in mind, notice that $\ssF{r}(\phi_{m_0,T=r})=\phi_{m_0,r}\in \fr{z}(\widehat{\fr{gl}}_{r})$, and so $\phi_{m_0,T=r}$ is a $\Lul{gl}{r}[t]$-invariant. 
As $\Lul{gl}{T}[t]$-invariance is a polynomial condition in $T,T^{-1}$, this implies that $\phi_{m_0,T}$ is $\Lul{gl}{T}[t]$-invariant.
Algebraic independence of the set $\{\phi_{m,T}:m\in \Z_+\}$ follows from similar reasoning.

Next, we need to show that $\{\phi_{m,T}:m\in \Z_+\}$ is a complete set of Segal--Sugawara vectors over $R$. 
We do this by using ultraproducts to show that for each $\alpha\in \C^{\times}$, $\{\phi_{m,T=\alpha}:m\in \Z_+\}$ is a complete set of Segal--Sugawara vectors over $\C$, and then concluding from commutative algebra that $\{\phi_{m,T}:m\in \Z_+\}$ is a complete set of Segal--Sugawara vectors over $R$.

Now recall that Proposition \ref{prop:ultraproduct_center_vertex_algebra} gives an isomorphism of vertex algebras over $\C$:
\begin{equation}
    \prod_{\mathcal{U}}^{B}\fr{z}(\widehat{\fr{gl}}_{n_0+p})
    \cong
    \fr{z}(\Lulhat{gl}{n_0}).
\end{equation}
This isomorphism is compatible with the PBW filtrations. Moreover, for each fixed $i$, the ultraproduct of Molev's generators in positive characteristic is identified with the specialized interpolated Segal--Sugawara vector $\phi_{i,T=n_0}$. Indeed, both sides are obtained from the same finite categorical formula, and ultraproducts preserve the categorical operations appearing in that formula.

It follows from Lemma \ref{lemma:large_p_filtration_lemma} that, as $\C$-vector spaces:
\begin{equation}
    F^m\fr{z}(\Lulhat{gl}{n_0})
    \cong 
    F^{m}\C[\mathsf{T}^j\phi_{i,T=n_0}:1\leq i\leq m,\ j\in \Z_{\geq 0}].
\end{equation}
By Lemma \ref{lemma:center_filtrations_exhaustive} the PBW filtration is exhaustive. This implies the elements $\{\phi_{m,T=n_0}:m\in \Z_{+}\}$ generate $\fr{z}(\Lulhat{gl}{n_0})$ as a differential algebra. Furthermore, since these vectors arise as the ultraproduct of free differential generators they generate a free differential subalgebra. Hence they form a complete set of Segal--Sugawara vectors at $T=n_0$.

The same argument applies to an arbitrary $\alpha\in \C^\times$, after choosing
an appropriate sequence $(\alpha_n)_{n=1}^{\infty}$ whose ultraproduct rank is $\alpha$.

By Lemma \ref{lemma:center_filtrations_exhaustive}, the $R$-module  $\mathscr{F}^m\fr{z}(\Lulhat{gl}{T})$ is finitely generated. Consider the inclusion:
\begin{equation}
    \iota_{m}:\mathscr{F}^m R[\mathsf{T}^j\phi_{i,T}:i\in \Z_+,j\in \Z_{\geq 0}]\rightarrow \mathscr{F}^m\fr{z}(\Lulhat{gl}{T})
\end{equation}
For every $\alpha\in \C^{\times}$, the specialization of $\iota_m$ at $\fr{m}_{\alpha}:=(T-\alpha)$ is an isomorphism, by the completeness at the specialization $T=\alpha$. As both $R$-modules are finitely generated, Nakayama's lemma implies that $(\iota_m)_{\fr{m}_{\alpha}}$ is an isomorphism of $R_{\fr{m}_{\alpha}}$-modules for every $\alpha\in \C^{\times}$. Therefore, $\iota_m$ is an isomorphism of $R$-modules and $\{\phi_{m,T}:m\in \Z_+\}$ is a complete set of interpolated Segal--Sugawara vectors over $R$. 

Equations~\eqref{eqn:prop_type_A_ss_n_Gen_1}, \eqref{eqn:prop_type_A_ss_n_Gen_2} follow from evaluation.

The last statement follows from boundedness of the $\lambda$-bracket with respect to the PBW filtration. Indeed, for fixed $i_0,j_0$, the bracket $\{\phi_{i_0,T}{}_{\lambda}\phi_{j_0,T}\}$ lies in a fixed PBW-filtered piece,
for instance in degree at most $i_0+j_0$. For $n$ sufficiently large relative to this filtered degree, the interpolation principle implies that $\ssV{n}$ is an isomorphism on the relevant filtered piece. Since $\ssV{n}$ is compatible with the affine vertex algebra structure, it intertwines the corresponding $\lambda$-brackets. Thus the desired identity may be checked after applying $\ssV{n}$, where it is precisely the finite-rank statement.
\end{proof}

\begin{remark}
    It is interesting to note that $\alpha=0$ is excluded as we are unable to use \cite[Theorem 1.1]{Modular_Center} in the ultrafilter lemma since $p$ is a bad prime for $\fr{gl}_p$ in characteristic $p$.
\end{remark}

These two sets of interpolated Segal--Sugawara vectors are related to one another through the parity functor and the interpolated Cartan involution, which we will now explain.
Recall from the previous subsection that, for each $n\in \Z_+$, we have the Cartan involution:
\begin{equation}
    \nu_n:\fr{z}(\widehat{\fr{gl}}_n)\rightarrow \fr{z}(\widehat{\fr{gl}}_n)
\end{equation}
Using the categorical description of $\phi_{m,n}, \psi_{m,n}$, with the explicit calculations given in Equations~\eqref{equation:cartan_fintie_rank_calc_1}, \eqref{equation:cartan_fintie_rank_calc_2} we see that $\nu_n(\phi_{m,n})$ and $\nu_n(\psi_{m,n})$ can be described as the automorphism sending $\mathrm{coev}_{\bullet}\mapsto -\mathrm{coev}_{\bullet}$.
This naturally leads to the interpolated Cartan involution:

\begin{prop}[Interpolated Cartan Involution]

The following morphism is a well-defined Poisson vertex algebra automorphism of $\fr{z}(\Lulhat{gl}{T})$: 
\begin{equation}
\nu_T(\phi_{m,T}):=\pi\left(\sum_{\mu\vdash m,\ } (-1)^{\ell}Q_{m,\ell}(T)\cdot c_{\mu,m}\cdot \widetilde{A^{(\ell)}}\circ (\cv_{\bullet})^{\otimes \mu}\right)
\end{equation}
Furthermore, we have:
\begin{equation}
\nu_T(\psi_{m,T})=\pi\left(\sum_{\mu\vdash m,\ } (-1)^{m}Q_{m,\ell}(-T)\cdot c_{\mu,m}\cdot \widetilde{H^{(\ell)}}\circ (\cv_{\bullet})^{\otimes \mu}\right)
\end{equation}
Similarly, one can evaluate $T=\alpha$, for a non-zero $\alpha\in \C^{\times}$, and obtain a Poisson automorphism $\nu_{\alpha}:\fr{z}(\Lulhat{gl}{\alpha})\rightarrow\fr{z}(\Lulhat{gl}{\alpha})$.
\end{prop}

\begin{proof}
    To see that this is well defined, notice that $\nu_T(\phi_{i,T})\in \Hom_{\irep(\GL_{T})}(\mathbf{1}_{\cC},V^{-h}(\Lul{gl}{T}))$. It suffices to verify that it is $\Lul{gl}{T}[t]$-invariant.
    This follows from applying the interpolation principle. Similarly, the Poisson automorphism property and the action of $\nu_T(\psi_{i,T})$ follow from the interpolation principle.
\end{proof}

The following corollary relates the two complete sets of Segal--Sugawara vectors constructed by Molev in \cite[Theorem 2.1]{Molev_Type_B_C_Invaraints} to the parity functor $\Pi_{\alpha}:\irep(\GL_{\alpha},\F)\rightarrow \irep(\GL_{-\alpha},\F)^{Tw}$. Its proof follows from an easy calculation and so is omitted. 

\begin{corollary}
    The Poisson anti-isomorphism $\Pi_{\alpha}\circ \nu_{\alpha}:\fr{z}(\Lulhat{gl}{\alpha})\rightarrow \fr{z}(\Lulhat{gl}{-\alpha})$ satisfies the following:
    \begin{equation}
        (\Pi_{\alpha}\circ \nu_{\alpha})(\phi_{i,\alpha})=(-1)^{i}\psi_{i,-\alpha} \ \ \ (\Pi_{\alpha}\circ \nu_{\alpha})(\psi_{i,\alpha})=(-1)^i \phi_{i,-\alpha}
    \end{equation}
    \label{cor:anti_isomorphism_centers_generators}
\end{corollary}

As we will see later in the paper, the appearance of the sign $(-1)^i$ is closely related to the duality between complete and elementary symmetric functions.

The interpolated Segal--Sugawara vectors in $\irep(\Or_{T},\C[T])$, $\irep(\Sp_T,\C[T])$ are defined for $m\in \Z_+$ respectively by:
\begin{equation}
    \phi_{2m,T}^{B}:=\pi\left(\sum_{\substack{\mu\vdash 2m,\ \\ \ell \text{ even}}} Q_{2m,\ell}(-T+1)\cdot c_{\mu,2m}\cdot (h_{\bullet,2})^{\otimes \ell}\circ \widetilde{H^{(\ell)}}\circ (\cv_{\bullet})^{\otimes \mu}\right)
\end{equation}
\begin{equation}
    \phi_{2m,T}^{C}:=\pi\left(\sum_{\substack{\mu\vdash 2m,\ \\ \ell \text{ even}}}Q_{2m,\ell}(T+1)\cdot c_{\mu,2m}\cdot (e_{\bullet,2})^{\otimes \ell}\circ \widetilde{A^{(\ell)}}\circ (\cv_{\bullet})^{\otimes \mu}\right)
\end{equation}
where the superscript denotes the corresponding type of the interpolating category and $\pi:\mathcal{T}(\Lul{g}{T}[t^{-1}]t^{-1})\rightarrow U(\Lul{g}{T}[t^{-1}]t^{-1})$ is the canonical map from the tensor algebra to the universal enveloping algebra.

\begin{theorem}[Interpolated Segal--Sugawara Vectors: Type $B/C$]
\label{thm:interpolated_segal_sugawara_Type_B_C}

Let $G_{\alpha}=\Or_{\alpha}$ or $\Sp_{\alpha}$ and set $R=\C[T]$.
The Type $B$ family $\{\phi^B_{2m,T}: m\in \Z_+\}$ and Type $C$ family $\{\phi^C_{2m,T}: m\in \Z_+\}$ satisfy:  

\begin{enumerate}
    \item \textbf{Specialization:} For any $\alpha\in \C$ let $\mathrm{Ev}_{T=\alpha}:f_{\alpha}^*(\irep(G_{T},R))\rightarrow \irep(G_{\alpha},\C)$ denote the change of base functor induced by the evaluation morphism $f_{\alpha}:R\rightarrow \C$, where $p(T)\mapsto p(\alpha)$. By defining: \begin{equation}\phi_{2m,T=\alpha}^B:=\mathrm{Ev}_{T=\alpha}(\phi_{2m,T}^B)\qquad \phi_{2m,T=\alpha}^C:=\mathrm{Ev}_{T=\alpha}(\phi_{2m,T}^C)\end{equation} the sets $\{\phi^B_{2m,T=\alpha}:m\in \Z_+\}$ and $\{\phi^C_{2m,T=\alpha}:m\in \Z_+\}$ form complete sets of interpolated Segal--Sugawara vectors of $\fr{z}(\Lulhat{so}{\alpha})$ and $\fr{z}(\Lulhat{sp}{\alpha})$ respectively. In particular, $\mathrm{Ev}_{T=\alpha}(\fr{z}(\Lulhat{so}{T}))\cong \fr{z}(\Lulhat{so}{\alpha})$ and $\mathrm{Ev}_{T=\alpha}(\fr{z}(\Lulhat{sp}{T}))\cong \fr{z}(\Lulhat{sp}{\alpha})$
    as $\C$-Poisson vertex algebras.
     \item \textbf{Completeness:} The set $\{\phi^B_{2m,T}:m\in \Z_+\}$ is a complete set of interpolated Segal--Sugawara vectors for $\fr{z}(\Lulhat{so}{T})$ over $R$. The set $\{\phi^C_{2m,T}:m\in \Z_+\}$ is a complete set of interpolated Segal--Sugawara vectors for $\fr{z}(\Lulhat{sp}{T})$ over $R$. 
    \item \textbf{Compatibility with Finite Rank:} Under mapping back to finite rank, they satisfy:
    \begin{equation}
    \ssV{2n+1}(\phi^{B}_{2i,T=2n+1}) = \phi_{2i,2n+1}^{B} \quad \text{for } 1 \leq i \leq n.
\end{equation}
    \begin{equation}
        \ssV{2n}(\phi^{C}_{2i,T=2n})=
        \begin{cases}
            \phi_{2i,2n}^{C} & \text{ if } 1\leq i\leq n\\
            0 & \text{ otherwise}
        \end{cases}
    \end{equation}

    \item \textbf{Asymptotic $\lambda$-brackets:} Let $\Lul{g}{N}=\Lul{so}{N}$ or $\Lul{sp}{N}$, where $N=2n+1$ in the orthogonal case, and $N=2n$ in the symplectic case. Denote by $\ssIV{N}:\fr{z}(\widehat{\fr{g}}_N)\rightarrow \fr{z}(\Lulhat{g}{N})$ the section given by $\phi_{2i}\mapsto \phi_{2i,T=N}$ for $1\leq i\leq n$. For fixed positive integers $i_0,j_0$, there is a sufficiently large $N$ such that the following holds for the $\lambda$-bracket:
    \begin{equation}
        \ssIV{N}(\{ \phi_{2i_0,N} \ {}_{\lambda} \ \phi_{2j_0,N} \})=\{ \phi_{2i_0,T=N} \ {}_{\lambda} \ \phi_{2j_0,T=N} \}
    \end{equation}
\end{enumerate}
\end{theorem}

\begin{proof}
The proof is identical to the Type $A$ argument after replacing the finite-rank input by Molev's Type $B/C$ generators. The filtration lemma and ultraproduct lemma have the same proofs, with $\GL_T$ replaced by $\Or_T$ or $\Sp_T$ and with the corresponding interpolation principle. In particular, Types $B/C$ have no bad primes greater than $2$, which allows us to apply the ultraproduct construction for any $\alpha\in \C$.  
\end{proof}
The following Corollary relates the two complete sets of Segal--Sugawara vectors constructed by Molev in \cite[Theorem 2.1]{Molev_Type_B_C_Invaraints} to the parity functor $\Pi_{\alpha}:\irep(\Sp_{\alpha},\F)\rightarrow \irep(\Or_{-\alpha},\F)^{Tw}$.

\begin{corollary}
    The Poisson anti-isomorphism $\Pi_{\alpha}:\fr{z}(\Lulhat{sp}{\alpha})\rightarrow \fr{z}(\Lulhat{so}{-\alpha})$ satisfies the following:
    \begin{equation}
        \Pi_{\alpha}(\phi_{2i,\alpha}^C)=\phi_{2i,-\alpha}^B \ \
    \end{equation}
    \label{corollary:anti_isomorphism_generators}
\end{corollary}

\section{Applications to Lie Superalgebras and the Interpolated Gaudin Model}
\label{section:applications}

One of the most interesting applications of representation theory in complex rank is that it enables classical constructions to be transported to super representation theory via specific functors. In this section, we show that the interpolated higher Segal--Sugawara vectors constructed in Section \ref{section:interpolating_center_critical_level} recover substantial parts of the super Lie algebra story. More precisely, we show that the main theorem of \cite{Molev_Segal_Sugawara_Vectors_super_osp} on higher Segal--Sugawara vectors for $\fr{osp}_{M|2n}$ follows as a straightforward corollary. The analogous construction for $\fr{gl}_{m|n}$ from \cite{Molev_Segal_Sugawara_Vectors_super_gl} is recovered in the same way. 

Furthermore, it is well-known that the higher Hamiltonians of the Gaudin model, the so-called Bethe algebra, can be constructed using the Feigin--Frenkel center \cite{Gaudin_FF_center_first_ref} by evaluating at certain points.  While the authors of \cite{Rybnikov} construct the interpolated Bethe algebra by indirectly interpolating certain generators, we show that this algebra can be recovered directly from the interpolated Feigin--Frenkel center. 
\subsection{Application to Lie Superalgebras}

The connection between interpolating categories and super-representation theory rests on the following universal property and the corresponding Schur--Weyl dualities. Throughout this subsection, we assume non-zero superdimension in the general linear case, i.e., $m \neq n$.

\begin{theorem}[Universal Property of Interpolating Categories \cite{St_Interpolating_Deligne}] 
    Let $\cD$ be a symmetric pseudo-tensor category over $R$ of characteristic $0$.  If $V\in \cD$ is an object of dimension $t\in R$, then there exists a braided pseudo-tensor functor:
    \begin{equation}
        F:\irep(\GL_{t},R)\rightarrow \cD \qquad F(\bullet):= V
    \end{equation}
    Furthermore, if $V$ has a symmetric non-degenerate bilinear form, then there exists a braided pseudo-tensor functor: 
    \begin{equation}
        K:\irep(\Or_t,R)\rightarrow \cD \qquad K(\bullet):=V
    \end{equation}
    with a similar statement for symplectic non-degenerate bilinear forms. 
\end{theorem}

\begin{theorem}[Mixed Schur--Weyl Duality in the Super Case \cite{Mixed_Schur_WEyl_Duality_Super_gl,Schur_Weyl_Duality_Orthosymplectic}]
In the general linear case, let $V:=\C^{m|n}$ and $W:=V^*$. Then there is a surjective $\C$-algebra morphism: 
\begin{equation}
    B_{r,s}(m-n)\rightarrow \mathrm{End}_{\fr{gl}_{m|n}}(V^{\otimes r}\otimes W^{\otimes s})
\end{equation}
which is an isomorphism if $|m-n|\geq r+s$. 

In the orthosymplectic case, let $V:=\C^{M|2n}$. Then there is a surjective algebra morphism: 
\begin{equation}
    B_{r}(M-2n)\rightarrow \mathrm{End}_{\fr{osp}_{M|2n}}(V^{\otimes r})^{\mathrm{OSp}_{M|2n}}
\end{equation}
which is an isomorphism if $|M-2n|\geq r$. 
\end{theorem}
\begin{remark}
  We take $\mathrm{OSp}_{M|2n}$-invariants as we would like to work in $\mathrm{Rep}_{\mathrm{OSp}_{M|2n}}(\fr{osp}_{M|2n})$, namely the category of $\fr{osp}_{M|2n}$-modules internal to $\mathrm{Rep}(\mathrm{OSp}_{M|2n})$. Similar to the non-super case this is equivalent to working $\Z/2$-equivariantly when $M\in 2\Z_+$, see Remark \ref{remark:Pfaffian_killed}. 
\end{remark}
Combining these facts, we obtain the following proposition: 
\begin{prop}
    There exist full, braided pseudo-tensor functors over $\C$:
    \begin{equation} 
        F_{m|n}: \irep(\GL_{m-n},\C)\rightarrow \mathrm{Rep}_{\mathrm{GL}_{m|n}}(\fr{gl}_{m|n})
    \end{equation}
    \begin{equation}
        K_{M|2n}: \irep (\Or_{M-2n},\C)\rightarrow \mathrm{Rep}_{\mathrm{OSp}_{M|2n}}(\fr{osp}_{M|2n})
    \end{equation}
where $\mathrm{Rep}_G(\mathfrak{g})$ denotes the representations of $\mathfrak{g}$ internal to the category $\mathrm{Rep}(G)$. 
\end{prop}

As the functors $F_{m|n}$ and $K_{M|2n}$ are
braided pseudo-tensor functors, they preserve all of the
algebraic structures constructed in
Section~\ref{section:vertex_algebras_universal_affine_interpolated}:
Lie algebras, affinizations, universal enveloping
algebras, and vacuum modules. Composing $K_{M|2n}$
with the isomorphism induced by the standard
supersymmetric bilinear form
$\theta : \C^{M|2n} \to (\C^{M|2n})^*$ (see \cite[Section 2]{Molev_Segal_Sugawara_Vectors_super_osp} for details),
we obtain identifications:
\begin{equation}
    F_{m|n}(\Lul{gl}{m-n}) = \fr{gl}_{m|n},
    \qquad
    \Id\otimes \theta:K_{M|2n}(\Lul{so}{M-2n})
    \cong \fr{osp}_{M|2n},
\end{equation}
\begin{equation}
    F_{m|n}(\Lulhat{gl}{m-n}) = \widehat{\fr{gl}}_{m|n},
    \qquad
    K_{M|2n}(\Lulhat{so}{M-2n})
    \cong \widehat{\fr{osp}}_{M|2n},
\end{equation}
\begin{equation}
    F_{m|n}(V^k(\Lul{gl}{m-n})) = V^k(\fr{gl}_{m|n}),
    \qquad
   \tilde{\theta}: K_{M|2n}(V^k(\Lul{so}{M-2n}))
    \cong V^k(\fr{osp}_{M|2n}).
\end{equation}
Here $\tilde{\theta}$ is the isomorphism of supervertex algebras induced by $\theta$.
In the super setting, the natural morphism spaces are internal Hom's, which include both even and odd maps. Consequently, the Feigin--Frenkel center defined through external Hom captures only the even part. For this reason we make the following distinction.
\begin{definition}
Let $\fr{g}$ be a Lie superalgebra. The \emph{even
Feigin--Frenkel center} is defined as:
\begin{equation}
\fr{z}(\widehat{\fr{g}})_0
:= \mathrm{Hom}_{\fr{g}[t]}
   (\mathbf{1}, V^{-h^{\vee}}(\fr{g}))
\end{equation}
\end{definition}
Applying Theorem \ref{theorem:invariant_ind_PVA_structure} with $F_{m|n}$, $K_{M|2n}$ to the interpolated higher Segal--Sugawara vectors constructed in Theorems \ref{theorem:interpolate_segal_sugawara_type_A} and \ref{thm:interpolated_segal_sugawara_Type_B_C} we obtain: 
\begin{theorem}
\label{theorem:segal_sugawara_super_case}
There is a morphism of $\C$-Poisson vertex algebras: 
\begin{equation}
    \fr{z}_k(\Lulhat{gl}{m-n})\rightarrow \fr{z}_k(\widehat{\fr{gl}}_{m|n})_0
\end{equation}
In particular, at the critical level, we obtain higher Segal--Sugawara vectors in $\fr{z}(\widehat{\fr{gl}}_{m|n})$ for $r\in \Z_+$:  
\begin{align}
    F_{m|n}(\phi_{r,T=m-n})&=\sum_{\mu\vdash r,}Q_{r,\ell}(m-n)c_{\mu,r} \sum_{\sigma\in S_{\ell}} \mathrm{sgn}(\sigma) \sum_{i_1,\cdots, i_{\ell}=1}^{m+n} (-1)^{K(\underline{i},\sigma)}E_{i_{\sigma(1)},i_1}t^{-\mu_1}\cdots E_{i_{\sigma(\ell)},i_{\ell}}t^{-\mu_{\ell}}
   \\
    F_{m|n}(\psi_{r,T=m-n})&=\sum_{\mu\vdash r,}Q_{r,\ell}(m-n)c_{\mu,r} \sum_{\sigma\in S_{\ell}}   \sum_{i_1,\cdots, i_{\ell}=1}^{m+n}(-1)^{K(\underline{i},\sigma)}E_{i_{\sigma(1)},i_1}t^{-\mu_1}\cdots E_{i_{\sigma(\ell)},i_{\ell}}t^{-\mu_{\ell}}
\end{align}
where $(-1)^{K(\underline{i},\sigma)}$ is the sign obtained by applying the braiding $\tilde{\sigma}$ to $(e_{i_1}\otimes e_{i_1}^*)\otimes \cdots \otimes (e_{i_{\ell}}\otimes e_{i_{\ell}}^*)$ multiplied by $(-1)^{\overline{i}_1+\cdots +\overline{i}_{\ell}}$.  Up to rescaling, these are the first two sets of Segal--Sugawara vectors given in \cite[Corollary 3.3]{Molev_Segal_Sugawara_Vectors_super_gl}. 
\end{theorem}

\begin{theorem}
\label{theorem:Segal_Sugawara_Super_osp}
There is a morphism of $\C$-Poisson vertex algebras: 
\begin{equation}
    \fr{z}_k(\Lulhat{so}{M-2n})\rightarrow \fr{z}_k(\widehat{\fr{osp}}_{M|2n})_0
\end{equation}
In particular, at the critical level, we obtain higher Segal--Sugawara vectors in $\fr{z}(\widehat{\fr{osp}}_{M|2n})$ for even $r\in \Z_+$: 
\begin{equation}
    (\tilde{\theta}\circ K_{M|2n})(\phi^B_{r,T=M-2n})=\sum_{\substack{\mu\vdash r\\  \ even}}c_{\mu,r}\cdot Q_{r,\ell}(2n-M+1) \sum_{\sigma\in S_{\ell}}\sum_{i_1,\cdots, i_{\ell}=1}^{M+2n} (-1)^{K(\underline{i},\sigma)}F_{i_{\sigma(1)},i_1}t^{-\mu_1}\cdots F_{i_{\sigma(\ell)},i_{\ell}}t^{-\mu_{\ell}} 
\end{equation}
where $F_{i,j}=(1-c^{\mathrm{sVect}}_{V^*,V}\circ (\theta\otimes \theta^{-1}))(E_{i,j})$, and $K(\underline{i},\sigma)$ is as in Theorem \ref{theorem:segal_sugawara_super_case}.  Up to rescaling, this is the set of Segal--Sugawara vectors given in \cite[Equation 2.6]{Molev_Segal_Sugawara_Vectors_super_osp}. 
\end{theorem}

\begin{proof}
    Note that the polynomial coefficients in \cite[Equation 2.6]{Molev_Segal_Sugawara_Vectors_super_osp} will be the same as ours as $\ell,r$ are both even. 
\end{proof}

In \cite{Molev_Segal_Sugawara_Vectors_super_osp} the authors conjecture that when $M$ is odd, the Feigin--Frenkel center will be generated as a differential algebra by the set of Segal--Sugawara vectors in Theorem \ref{theorem:Segal_Sugawara_Super_osp}. One can prove an asymptotic version of this conjecture for the even part $\fr{z}(\widehat{\fr{g}})_0$ using mixed Schur--Weyl duality: 

\begin{theorem}[Asymptotic Description of Super Feigin--Frenkel Center]
\label{theorem:asymptotic_description}
Fix $m_0\in \Z_{+}$. In each case, we consider $\fr{z}(\widehat{\fr{g}})$ with the PBW filtration denoted by $F^{\bullet}\fr{z}(\widehat{\fr{g}})$. For sufficiently large $m-n$ the following holds:
\begin{equation}F^{m_0}(\fr{z}(\widehat{\fr{gl}}_{m|n}))_0= F^{m_0}\C[\mathsf{T}][F_{m|n}(\phi_{1,T=m-n}),\cdots,F_{m|n}(\phi_{m-n,T=m-n}) ]
\end{equation}
For sufficiently large $M-2n$, assuming $M$ is odd, the following holds: 
\begin{equation}F^{m_0}\fr{z}(\widehat{\fr{osp}}_{M|2n})_0= F^{m_0}\C[\mathsf{T}][(\tilde{\theta}\circ K_{M|2n})(\phi_{2,T=M-2n}^B),(\tilde{\theta}\circ K_{M|2n})(\phi_{4,T=M-2n}^B),\cdots]
\end{equation}
\end{theorem}
\begin{proof}
    This follows from the interpolation principles (Propositions \ref{prop:interpolation_principle_A}, \ref{prop:interpolation_principle_BC}) combined with the descriptions given in Theorems \ref{theorem:interpolate_segal_sugawara_type_A}, \ref{thm:interpolated_segal_sugawara_Type_B_C}. 
\end{proof}
\subsection{The Interpolated Gaudin Model}

In this subsection we show that the interpolated Bethe algebra of \cite{Rybnikov} can be recovered directly from the interpolated Feigin--Frenkel center, exactly as in finite rank. We briefly recall the setup from \cite{Rybnikov}. The quantum Gaudin model for $\fr{gl}_n$ is a quantum spin chain that takes $m$ representations of $\fr{gl}_n$ with highest weights $\mu_1,\dots, \mu_m$ and assigns to distinct points $z_1,\dots, z_m\in \C$ sites $V_{\mu_1},\dots, V_{\mu_m}$, where $V_{\mu_a}$ denotes the irreducible finite-dimensional $\fr{gl}_n$ module with highest weight $\mu_a$. The total Hilbert space is $V(\underline{\mu}):=V_{\mu_1}\otimes \cdots \otimes V_{\mu_m}$, and the commuting Hamiltonians are, for $a=1,\cdots,m$:
\begin{equation}
    H_a:=\sum_{\substack{b=1\\ b\neq a}}^m\frac{\Omega^{ab}}{z_a-z_b}\in U(\fr{gl}_n)^{\otimes m}
\end{equation}
where:
\begin{equation}
    \Omega=\sum_{i,j=1}^nE_{i,j}\otimes E_{j,i}
\end{equation}
is the quadratic Casimir, and for $\Omega^{ab}$, the first tensor factor acts on the $a$-th component and the second tensor factor acts on the $b$-th component.  These Hamiltonians commute with the diagonal action of $\fr{gl}_n$ on $V(\underline{\mu})$, and so naturally they preserve the subspace of singular vectors of $V(\underline{\mu})$, denoted $V(\underline{\mu})^{\mathrm{sing}}$. The Bethe algebra is the centralizer of the images of the $H_a$ in $\mathrm{End}(V(\underline{\mu})^{\mathrm{sing}})$ and is denoted by $\mathcal{B}_{n,\underline{\mu}}^{\mathrm{sing}}$. 

Equivalently, there is a \emph{universal Bethe algebra}:
\[
\mathcal{B}_{n,m}\subset U(\fr{gl}_n)^{\otimes m}
\]
whose action on $V(\underline{\mu})$ is $\fr{gl}_n$-equivariant. Let:
\[
\mathcal{B}_{n,\underline{\mu}}\subset \mathrm{End}_{\fr{gl}_n}(V(\underline{\mu}))
\]
denote the image of this action. Since the action is $\fr{gl}_n$-equivariant, $\mathcal{B}_{n,\underline{\mu}}$ preserves $V(\underline{\mu})^{\mathrm{sing}}$, and the restriction of this image to singular vectors is isomorphic to $\mathcal{B}_{n,\underline{\mu}}^{\mathrm{sing}}$. Interestingly, one can construct the universal Bethe algebra from $\fr{z}(\widehat{\fr{gl}}_n)$ \cite{Gaudin_FF_center_first_ref}, which by the above discussion means all Bethe algebras may be recovered. 

In particular, once one fixes a set of generators of the universal Bethe algebra, the Bethe algebra $\mathcal{B}_{n,\underline{\mu}}^{\mathrm{sing}}$ is determined by the images of these generators in:
\[
\End_{\fr{gl}_n}(V(\underline{\mu})).
\]
The authors of \cite{Rybnikov} use the fact that the universal Bethe algebra admits generators which can be described diagrammatically in $\mathrm{Rep}(\GL_n,\C)$ and hence interpolated to $\irep(\GL_{\alpha},\C)$. They then take the images of these generators in:
\[
\End_{\irep(\GL_{\alpha},\C)}(V_1\otimes\cdots\otimes V_m)
\]
for indecomposable \(V_i\in \irep(\GL_{\alpha},\C)\). The interpolated Bethe algebra is defined to be the subalgebra generated by these images. More concretely, $\mathcal{B}_{n,m}$ has generators $S^{(r)}_{k\ell j}\in \fr{gl}_n^{\otimes m}$ indexed by certain natural numbers $(k,\ell,j,r)$ that do not depend on $n$. These generators can be written as certain linear combinations of morphisms $\mathbf{1}_{\mathrm{Rep}(\GL_n)}\rightarrow \fr{gl}_n^{\otimes m}$ given by coevaluations and braidings, similar to the strategy in Subsection \ref{section:interpolating_molevs_invariants}. They then note that these coefficients are independent of $n$ and so they naturally generalize to the context of $\irep(\GL_{\alpha},\C)$ and are denoted by:
\begin{equation}
    \mathcal{S}_{k\ell j}^{(r)}\in \mathrm{Hom}_{\irep(\GL_{\alpha},\C)}(\mathbf{1}_{\irep(\GL_{\alpha},\C)},(\Lul{gl}{\alpha})^{\otimes m})
\end{equation}
for the same tuples of natural numbers $(k,\ell,j,r)$. This motivates the following. 
\begin{definition}[Interpolated Universal Bethe Algebra]

Let $\alpha\in \C$. The interpolated universal Bethe algebra is the $\C$-algebra generated by the elements $\mathcal{S}_{k\ell j}^{(r)}$:
\begin{equation}
    \underline{\mathcal{B}}_{\ \alpha,m}:=\C [\mathcal{S}_{k\ell j}^{(r)}] \subset \mathrm{Hom}_{\irep(\GL_{\alpha},\C)}(\mathbf{1}_{\irep(\GL_{\alpha},\C)},(U(\Lul{gl}{\alpha}))^{\otimes m})
\end{equation}
\end{definition}
  \begin{remark}
    The definition above is a slight reformulation of the construction in \cite{Rybnikov}. There, the authors construct the interpolated generators and use their images to define Bethe algebras associated to objects of \(\irep(\GL_{\alpha},\C)\). Since we have constructed universal enveloping algebras in our setting, we may instead regard the same generators as generating a universal Bethe algebra before passing to any particular representation.
\end{remark}

The interpolated universal Bethe algebra $\underline{\mathcal{B}}_{\ \alpha,m}$ can be obtained from the interpolated Feigin--Frenkel center $\fr{z}(\Lulhat{gl}{\alpha})$ just as in the classical case.  Namely, one chooses distinct complex numbers $u,z_1,\cdots, z_m$ and defines evaluation morphisms: 
\begin{equation}
    \mathrm{ev}_{z}^t:U(\Lul{gl}{\alpha}[t^{-1}]t^{-1})\rightarrow U(\Lul{gl}{\alpha})\qquad X t^n\mapsto X z^n
\end{equation}
For the $a$-th tensor slot, set: 
\begin{equation}(\mathrm{ev}_{z_a-u})^{(a)}:=\Id\otimes\cdots\otimes \mathrm{ev}_{z_a-u}^t\otimes \cdots \otimes  \Id
\end{equation} 
Using these evaluation morphisms, one defines the map:
\begin{equation}
    \tilde{\Psi}_m:\Lul{gl}{\alpha}[t^{-1}]t^{-1}\rightarrow U(\Lul{gl}{\alpha})^{\otimes m} \qquad \tilde{\Psi}_m:=\left(\sum_{a=1}^m(\mathrm{ev}_{z_a-u})^{(a)}\right)\circ \Delta
\end{equation}
where $\Delta:\Lul{gl}{\alpha}[t^{-1}]t^{-1}\rightarrow (\Lul{gl}{\alpha}[t^{-1}]t^{-1})^{\otimes m}$ denotes the iterated coproduct map. 
The map $\tilde{\Psi}_m$ extends multiplicatively to the universal enveloping algebra, giving:
\begin{equation}
    \Psi_m:U(\Lul{gl}{\alpha}[t^{-1}]t^{-1})\rightarrow U(\Lul{gl}{\alpha})^{\otimes m}
\end{equation}
Composing with the involution $\zeta: U(\Lul{gl}{\alpha}[t^{-1}]t^{-1})\rightarrow U(\Lul{gl}{\alpha}[t^{-1}]t^{-1})$ given by $X t^n \mapsto -X t^n$, we obtain the desired map on the object underlying the vacuum module: 
\begin{equation}
    \tilde{\Phi}_m:U(\Lul{gl}{\alpha}[t^{-1}]t^{-1})\rightarrow  U(\Lul{gl}{\alpha})^{\otimes m} \quad \tilde{\Phi}_m:= \Psi_m\circ \zeta
\end{equation}
 In the finite-rank case, one recovers the universal Bethe algebra by restricting $\tilde{\Phi}_m$ to the Feigin--Frenkel center, giving:
\begin{equation}
    \Phi_m:\fr{z}(\widehat{\fr{gl}}_n)\rightarrow U(\fr{gl}_{n})^{\otimes m} \quad \Phi_m(\fr{z}(\widehat{\fr{gl}}_n))=\mathcal{B}_{n,m}
\end{equation}
See \cite[Section 14.1]{Molev} for more details. In our setting, we have natural inclusions:
\begin{equation}
    \fr{z}(\Lulhat{gl}{\alpha})\hookrightarrow \mathrm{Hom}_{\irep(\GL_{\alpha},\C)}(\mathbf{1}_{\irep(\GL_{\alpha},\C)},V^{-h^{\vee}}(\Lul{gl}{\alpha}))\cong  \mathrm{Hom}_{\irep(\GL_{\alpha},\C)}(\mathbf{1}_{\irep(\GL_{\alpha},\C)},U(\Lul{gl}{\alpha}[t^{-1}]t^{-1}))
\end{equation}
Therefore, there is a natural map:
\begin{equation}
    \Phi_{m}: \fr{z}(\Lulhat{gl}{\alpha})\rightarrow  \mathrm{Hom}_{\irep(\GL_{\alpha},\C)}(\mathbf{1},U(\Lul{gl}{\alpha})^{\otimes m})
\end{equation}

\begin{theorem}
\label{theorem:bethe_algebra_interpolated_construction}

$\Phi_m$ is a homomorphism of $\C$-algebras, and its image $\Phi_m(\fr{z}(\Lulhat{gl}{\alpha}))$ is equal to the interpolated universal Bethe algebra $\underline{\mathcal{B}}_{\ \alpha,m}$.
\end{theorem}
\begin{proof}
This follows from the interpolation principle applied to the finite-rank construction recalled above.
\end{proof}

\part{Adler--Gelfand--Dickey Algebras at Complex Rank}

\section{Drinfeld--Sokolov Reduction of \texorpdfstring{$\fr{gl}_{\lambda}$}{gl-lambda} and \texorpdfstring{$\fr{po}_{\lambda}$}{po-lambda}}
\label{section:DS_complex_motivation}

We now turn to the Adler--Gelfand--Dickey side of the story.  This section reviews the finite-rank and complex-rank forms of Drinfeld--Sokolov reduction that motivate the definition of the Poisson vertex algebras $\cW(\fr{gl}_{\lambda})$ and $\cW(\fr{po}_{\lambda})$.  The discussion is primarily expository. We follow \cite{Khesin_Malikov_DS_reduction,Drinfeld_Sokolov}.
\color{black}
\subsection{The Second Gelfand--Dickey Bracket and Drinfeld--Sokolov Reduction}
We begin with the classical finite-rank picture, as it serves as a template for all later constructions.  Consider the manifold $M$ of monic order $n$ differential operators on the circle:
\begin{equation}
M=\{\pr^{n}+a_1(x)\pr^{n-1}+\cdots+ a_{n}(x)\mid a_i(x)\in C^{\infty}(\mathbb{R}/\mathbb{Z},\C)\}
\end{equation}
For our algebraic setting, we replace the functions $a_i(x)$ by differential
indeterminates and work over the differential polynomial algebra:
\begin{equation}
\widetilde A
:=
\C[u_i^{(r)}\mid 1\leq i\leq n,\ r\geq 0],
\qquad
\partial u_i^{(r)}=u_i^{(r+1)}.
\end{equation}
The second Gelfand--Dickey bracket is the corresponding Poisson structure on $M$. Equivalently, a Hamiltonian mapping $A^{(L)}:T_L^*M\rightarrow T_LM$ for every $L\in M$. The tangent space at $L$ is defined as differential operators of the form $K:=v_{1}\pr^{n-1}+\cdots +v_n$. The cotangent space can be identified, through the residue pairing, with pseudodifferential operators of the form:
\begin{equation}
    X:=\pr^{-1}\circ v_1+\cdots +\pr^{-n}\circ v_n
\end{equation}
For a pseudodifferential operator $P\in \widetilde{A}((\pr^{-1}))$ the residue  is defined as: 
\begin{equation}
    \mathrm{res}(P)= \pr^{-1}  \text{ coefficient of }P
\end{equation}
The residue pairing is then given by $\langle K,X\rangle :=\mathrm{res}( K\circ X)$. With this pairing in mind, the second Adler--Gelfand--Dickey bracket arises from the Adler-map $A^{(L)}$ defined as: 
\begin{equation}
    A^{(L)}(X):=(LX)_+L-L(XL)_+
\end{equation}
Here $(-)_+$ denotes taking the differential operator part of a pseudodifferential operator.

In their seminal work \cite{Drinfeld_Sokolov}, Drinfeld and Sokolov showed that the second Adler--Gelfand--Dickey bracket can be realized as a special case of Hamiltonian reduction on the affine Kac--Moody algebra $\widehat{\fr{gl}}_n$. Specifically, one considers a reductive Lie algebra with decomposition:
\begin{equation}
\fr{g}=\fr{n}_-\oplus \fr{h}\oplus \fr{n}_+
\end{equation}
determined by a principal nilpotent $f$. The relevant manifold $\widetilde{M}$ consists of operators of the form:
\begin{equation}
    \fr{L}=\pr_x+q+f
\end{equation}
where $q:S^1\rightarrow \fr{b}:=\fr{n}_+\oplus \fr{h}$. The gauge action on $\widetilde{\mathcal{M}}$ is given by operators $\mathrm{e}^{\mathrm{ad} S}, S:S^1\rightarrow \fr{n}_+$. Any two operators equivalent under this action are called gauge equivalent. Drinfeld--Sokolov reduction is obtained by passing to gauge-invariant functions, or equivalently to the quotient by this gauge action. Drinfeld and Sokolov proved that $\mathcal{M}$ has a natural Poisson structure. For $\fr{g}=\fr{gl}_n$, they show that the Poisson structure on $\mathcal{M}$ is naturally equivalent to $M$ with the second Adler--Gelfand--Dickey bracket.

Both of these constructions were later recast in the language of Poisson vertex algebras \cite{ADG_SKV,DS_Kac_DeSole} 
where they were shown to be isomorphic as Poisson vertex algebras via the Miura map \cite[Theorem 2.23]{ADG_SKV}. 

\subsection{The Complex-Rank Setting}

The complex-rank analogue is obtained by replacing finite-order differential operators with pseudodifferential symbols carrying a complex degree.  This is the setting in which the Khesin--Malikov picture naturally lives.

As mentioned above, a \emph{pseudodifferential operator} over $\widetilde{A}$ is an element of $\widetilde{A}((\pr^{-1}))$.  A \emph{pseudodifferential operator of complex rank} is then an expression of the form:
\begin{equation}
   \left(1+\sum_{i=-\infty}^{-1} u_i(x)\pr^i\right)\circ \pr^t \qquad t\in \C
\end{equation}
The collection of all such operators $\Psi\mathrm{DO}$ is the Lie group of pseudodifferential symbols introduced by Khesin and Zakharevich in \cite{Khesin_Zakharevich}.  Similar to how one can define a Poisson structure on monic differential operators of order $n$ through the second Adler--Gelfand--Dickey bracket, they define a generalized second Adler--Gelfand--Dickey bracket on $\Psi \mathrm{DO}$, endowing the space with a Poisson-Lie group structure. 

The generalized second Gelfand--Dickey bracket can alternatively be realized as Drinfeld--Sokolov reduction on Feigin's Lie algebras of complex rank. In \cite{Feigin_gl_lambda}, Feigin's Lie algebra $\fr{gl}_{\lambda}$ is defined as:
 \begin{equation}
     \fr{gl}_{\lambda}:=U(\fr{sl}_2)/\langle \Delta-(\lambda^2-1)/{2}\rangle 
 \end{equation} where $\Delta$ is the usual Casimir element of $\fr{sl}_2$. In particular, $\fr{gl}_{\lambda}=\fr{gl}_{-\lambda}$. For $\lambda=n$, there exists a natural quotient $\fr{gl}_{\lambda=n}\rightarrow \fr{gl}_n$, and therefore $\fr{gl}_{\lambda}$ can be considered as an interpolation of the Type $A$ case. Notice that this is reminiscent of the functor $\ssF{n}:\irep(\GL_{n},\C)\rightarrow \mathrm{Rep}(\GL_n,\C)$. Similarly, for the symplectic and orthogonal case Feigin introduced the Lie algebra $\fr{po}_{\lambda}\subset \fr{gl}_{\lambda}$ defined by:
\begin{equation}
    \fr{po}_{\lambda}:=\{x\mid  x\in \fr{gl}_{\lambda}, x^*=x\}
\end{equation}
where $*$ is a certain involution of $\fr{gl}_{\lambda}$. These algebras admit natural surjections for $n\in \Z_+$ to $\fr{sp}_{2n},\fr{so}_{2n+1}$ when $\lambda=2n$ and $2n+1$ respectively. As in the $\fr{gl}_{\lambda}$ case, we have $\fr{po}_{\lambda}=\fr{po}_{-\lambda}$. In particular, there is a natural surjection from $\fr{po}_{-2n}\rightarrow \fr{sp}_{2n}$. Compare this with the identification of $\irep(\Sp_{\alpha},\C)$ with $\irep(\Or_{-\alpha},\C)$ up to a twist of the braiding. Khesin and Malikov showed that one can apply Drinfeld--Sokolov reduction to $\fr{gl}_{\lambda}$, obtaining a Poisson structure on $\Psi DS_{\lambda}$ that coincides with the generalized second Adler--Gelfand--Dickey bracket \cite[Theorem 3.8]{Khesin_Malikov_DS_reduction}.  Similarly, applying Drinfeld--Sokolov reduction to $\fr{po}_{\lambda}$ yields a Poisson structure on the space of \emph{self-adjoint} pseudodifferential operators of complex rank, which also agrees with the restricted second Adler--Gelfand--Dickey bracket. 

The role of Part II is to supply the Poisson vertex algebra form of this
complex-rank Adler--Gelfand--Dickey theory. The constructions of Khesin--Zakharevich and Khesin--Malikov produce the corresponding Poisson-geometric structures on pseudodifferential symbols, but the comparison with the Feigin--Frenkel center requires a more algebraic object: a differential algebra of coefficients equipped with an explicit $\lambda$-bracket.  We construct this structure by extending the Adler map formalism of \cite{ADG_SKV} to pseudodifferential symbols of complex degree and to the relevant localizations of \(\C[T]\). This produces the Poisson vertex algebras
\[
    \cW(\fr{gl}_{\lambda})
    \qquad \text{and} \qquad
    \cW(\fr{po}_{\lambda})
\] which are the complex-rank Adler presentations used throughout Part II.
It should be noted that De Sole, Kac, and Valeri \cite[Section 2.3]{ADG_SKV} described the Poisson vertex algebra structure on pseudodifferential operators of rank $N\in \Z_+$ arising from a generic Adler map. Furthermore, Feigin, Rybnikov, and Uvarov \cite{Rybnikov} showed that, in the supercase, the pseudodifferential operators one should consider are ratios of certain differential operators (see \cite[Section 7]{Rybnikov}). A similar story may emerge in the context of super Feigin--Frenkel duality at the critical level.

\section{Poisson Vertex Algebra Background}
\label{section:PVA_background}

Throughout this section, $R$ denotes a commutative $\mathbb{F}$-algebra of characteristic $0$.  The purpose of the section is to develop a Poisson vertex algebra formalism that can accommodate complex powers of $\pr$.  The first step is therefore to introduce $\alpha$-shifted distributions and the corresponding algebra of shifted pseudodifferential symbols.
\subsection{$\alpha$-Shifted Distributions}

The shifted formalism begins by allowing distributions with exponents translated by an arbitrary parameter $\alpha\in R$.

\begin{definition}[$\alpha$-Shifted Distribution]
Let $\cA$ be an $R$-module, $\alpha\in R$. An $\cA$-valued \emph{$\alpha$-shifted distribution} in one variable $z$ is a formal series:
\begin{equation}
    a(z)=\sum_{m\in \Z}a_mz^{\alpha+m} \qquad a_m\in \cA
\end{equation}
In two variables it is:
\begin{equation}
    a(z,w):=\sum_{m,n\in \Z}a_{mn}z^{\alpha+m}w^{\alpha+n}
\end{equation}
\end{definition}
We call a $0$-shifted distribution just an $\cA$-valued distribution. Recall, the usual $\delta$-distribution is defined as:
\begin{equation}
    \delta(z-w):=\sum_{n\in \Z}\frac{z^{n}}{w^{n+1}}
\end{equation}
We define the $\alpha$-shifted delta distribution as:
\begin{equation}
    \delta_{\alpha}(z-w):=\sum_{n\in \Z}\frac{z^{n+\alpha}}{w^{n+1+\alpha}}
\end{equation}
The $\alpha$-shifted delta distribution has the following easily checked properties.
\begin{lemma}
\begin{enumerate}[label=(\alph*)]
\item For $\alpha\in R$ we have: 
\begin{equation}
\delta_{\alpha}(z-w)=\delta_{-\alpha}(w-z)
\label{eqn:delta_symmetry_switch_variables}
\end{equation}
\item If $a(z)$ is an $\alpha$-shifted distribution, then:
\begin{equation}
    a(z)\delta_{-\alpha}(z-w)=a(w)\delta_0(z-w)
    \label{eqn:delta_distribution_identity_change_variables}
\end{equation}
\end{enumerate}
\end{lemma}
The following notation is reserved only for $0$-shifted distributions:
\begin{equation}
    a(z)_+:=\sum_{n\in \Z_{\geq 0}}a_nz^n
\end{equation}
\begin{equation}
    a(z)_-=\sum_{n\in \Z_{<0}}a_nz^n
\end{equation}
We denote the power series expansion for $|z|>|w|$ by:
\begin{equation}
    i_z(z-w)^{-1}:=\sum_{n\in \Z_{\geq 0}}\frac{w^{n}}{z^{n+1}}
\end{equation}
Recall the identity:
\begin{equation}
    \delta(z-w)=i_z(z-w)^{-1}-i_w(z-w)^{-1}
\end{equation}
The residue of an $\cA$-valued distribution is defined as the $z^{-1}$-coefficient:
\begin{equation}
    \mathrm{Res}_{z}(\sum_{n\in \Z}a_nz^{-n-1}):=a_0
\end{equation}
\subsection{The Algebra of $\alpha$-Shifted Pseudodifferential Symbols}

We now package the shifted distributions into an algebra of pseudodifferential symbols, which is the natural home for the complex-rank Adler--Gelfand--Dickey formalism.

\begin{definition}[$R$-Linear Differential Algebra]
    An $R$-linear differential algebra $\cA$ is a unital commutative associative algebra over $R$ with $R$-linear derivation $\pr$.
\end{definition}

For $\alpha\in R$ and an $R$-linear differential algebra $\cA$, denote the $R$-module of all $\alpha$-shifted Laurent series with coefficients in $\cA$ by $\cA_{\alpha}((z))$. Recall that Equation \eqref{eqn:binomical_coefficient_over_ring} allows us to define binomial coefficients over a ring. 

\begin{definition}[Product of Shifted Pseudodifferential Operators]
\label{definition:product_pseudodifferential_operators}
Suppose that $\alpha_1,\alpha_2\in R$, the product of shifted pseudodifferential symbols $\circ:\cA_{\alpha_1}((\pr^{-1}))\times \cA_{\alpha_2}((\pr^{-1}))\rightarrow\cA_{\alpha_1+\alpha_2}((\pr^{-1}))$ is induced by:
\begin{equation}
    \pr^{\alpha_1+n}\circ (A\pr^{\alpha_2+m})=\sum_{k\in \Z_{\geq 0}}\binom{\alpha_1+n}{k}A^{(k)}\pr^{(\alpha_1+\alpha_2)+m+n-k}
\end{equation}
\end{definition}
$\alpha$-shifted pseudodifferential symbols serve as analogues of complex power pseudodifferential symbols mentioned in Section \ref{section:DS_complex_motivation}. The following notation will be needed later:
\begin{equation}
    \cA^{\alpha}((\pr^{-1})):= \pr^{\alpha}\circ \cA((\pr^{-1}))
    \label{eqn:right_alpha_shifted_notation}
\end{equation}

\begin{definition}[$\alpha$-Shifted Order/Degree]
    Let $A(\pr)=\sum_{k=-N}^{\infty}A_k\pr^{\alpha-k}$ be such that $A_{-N}\neq 0$. The \emph{$\alpha$-shifted order/degree of $A(\pr)$} is defined to be $N$. We denote this by $\mathrm{ord}_{\alpha}(A(\pr))$
\end{definition}
Notice, the order depends on our choice of $\alpha$. For example:
\begin{equation}
    \mathrm{ord}_{3}(\pr^5)=2 \ \ \mathrm{ord}_5(\pr^5)=0
\end{equation}
We denote the space of $\alpha$-shifted pseudodifferential symbols of $\alpha$-shifted degree at most $N\in \Z$ by:
\begin{equation}
    (\cA_{\alpha}((\pr^{-1})))_N:=\{A(\pr)\in \cA_{\alpha}((\pr^{-1})): \mathrm{ord}_{\alpha}(A(\pr))\leq N\}
    \label{eqn:bounded_degree_notation}
\end{equation}
If the largest $\alpha$-shifted degree term has coefficient $1$ we say that it is monic.

\begin{definition}[$\alpha$-Shifted Adjoint]
\label{defn:delta_shifted_adjoint}
    Let $A(\pr)=\sum_{k\in \Z}A_k\pr^{\alpha+k}\in \cA_{\alpha}((\pr^{-1}))$. The $\alpha$-shifted adjoint of $A(\pr)$ is defined to be:
\begin{equation}
    A(\pr)^*=\sum_{k\in \Z}(-1)^k\pr^{\alpha+k}\circ A_k
\end{equation}
\end{definition}
One can show that $A(\pr)^{**}=A(\pr)$. Observe that this differs from the usual definition of the adjoint, see \cite[1.2]{ADG_SKV} for example. The trade-off is it allows us to uniformly deal with symplectic \emph{and} orthogonal cases. 

\begin{definition}[Symbol of an $\alpha$-Shifted Pseudodifferential Operator]
  For $A(\pr)=\sum_{k\in \Z}A_k\pr^{\alpha-k}$, the symbol of $A(\pr)$ is the Laurent series $A(z)=\sum_{k\in \Z}A_kz^{\alpha-k}$. Denote this bijection by $\mathrm{Symb}_z:\cA_{\alpha}((\pr^{-1}))\rightarrow \cA_{\alpha}((z^{-1}))$.
\end{definition}
As in the classical case, this mapping is not an algebra homomorphism but satisfies:
\begin{lemma}

Let $A(\pr)\in \cA_{\alpha_1}((\pr^{-1})), B(\pr)\in \cA_{\alpha_2}((\pr^{-1}))$:
\begin{equation}
    \mathrm{Symb}(A(\pr)\circ B(\pr))=A(z+\pr)B(z)
\end{equation}
Here the right-hand side we expand in positive powers of $\pr$:
\begin{equation}
    (z+\pr)^{\alpha+m}:=\sum_{t\in \Z_{\geq 0}}\binom{\alpha+m}{t}\pr^{t}z^{\alpha+m-t}
\end{equation}
 In particular, $\pr$ will act on the coefficients of $B(z)$ as a derivation of $\cA$. 
 \label{lemma:symbol_of_product}
\end{lemma}
\begin{proof}
    By linearity it suffices to show that for $f\in \cA, m,n\in \Z$:
    \begin{equation}
        \mathrm{Symb}(\pr^{\alpha_1+m}\circ (f\pr^{\alpha_{2}+n}))=(z+\pr)^{\alpha_1+m}f z^{\alpha_2+n}
    \end{equation}
Expanding the left-hand side we have:
\begin{equation}
    \mathrm{Symb}\left(\sum_{k\in\Z_{\geq 0}}\binom{\alpha_1+m}{k}\pr^k(f)\pr^{\alpha_1+\alpha_2+n+m-k}\right)=\sum_{k\in \Z_{\geq 0}}\binom{\alpha_1+m}{k}\pr^k(f)z^{\alpha_1+\alpha_2+n+m-k}
\end{equation}
Expanding the other side we obtain:
\begin{equation}
    \sum_{k\in \Z_{\geq 0}}\binom{\alpha_1+m}{k}\pr^k(f)z^{\alpha_1+\alpha_2+n+m-k}
\end{equation}
Hence, they are equal.
\end{proof}
Throughout the rest of this paper, we will always expand in terms of the positive power of $\pr$. Now note that for $f,g\in \cA$:
\begin{equation}
    (\pr^{m}\circ g)\cdot f=\pr^{m}\cdot (gf)
\end{equation}
Here $\cdot$ denotes the action of $\pr^{m}$ on $f$ as an element of $\cA$. 
We will need the following identity. 
\begin{lemma}
\label{lemma:delta_shift_generating_fct}
    Let $L=\sum_{n\in \Z}a_n \pr^{n+\alpha}\in \cA_{\alpha}(\pr)$, then, for $f\in \cA$, the following equality holds:
    \begin{equation}
        \delta_{\alpha}(w-\pr)\circ (f\cdot L(\pr))= \delta_0(w-\pr)\circ (L(w+\pr)^*\cdot f)
    \end{equation}
\end{lemma}
\begin{proof}
It suffices to check they are the same after applying $\mathrm{Symb}_z$. By taking the symbol of the left-hand side and repeatedly applying Equations~\eqref{eqn:delta_symmetry_switch_variables},~\eqref{eqn:delta_distribution_identity_change_variables} to the left-hand side we obtain: 
\begin{align}
    \delta_{\alpha}(w-(z+\pr))\cdot (f\cdot L(z)) &= \delta_{\alpha}(w-(z+\pr))\cdot \left(\sum_{n\in \Z}(f\cdot a_n)z^{n+\alpha}\right)\\ 
    &= \delta_{\alpha}(w-(z+\pr))\cdot\left(\sum_{n\in \Z}z^{n+\alpha}(a_n f)\right)\\
    &= \delta_{-\alpha}((z+\pr)-w)\cdot\left(\sum_{n\in \Z}z^{n+\alpha}(a_n f)\right)\\
    &= \delta_{-\alpha}((z+\pr)-w)\cdot\left(\sum_{n\in \Z}((z+\pr)-\pr)^{n+\alpha}(a_n f)\right)\\
    &= \delta_{0}(w-(z+\pr))\cdot \left(\sum_{n\in \Z}(w-\pr)^{n+\alpha}\cdot ( a_n f)\right) \label{eqn:lemma_alpha_distribution_part_1}\\
    &=\delta_0(w-(z+\pr))\cdot (L(w+\pr)^*\cdot f)
\end{align} 
 Observe that $L(w+\pr)^*\cdot f\in \cA_{\alpha}((w^{-1}))$ and so in particular it is not itself a Laurent series in $\pr$. Taking the symbol of the right hand we see that indeed the equation holds. 

\end{proof}

For the sake of clarity we also define base change for pseudodifferential symbols:
\begin{definition}
    Let $R,S$ be commutative $\mathbb{F}$-algebras, and $f:R\rightarrow S$ a morphism. For $\alpha\in R, \beta\in S$ the morphism $f$ induces a map of $R$-modules:
    \begin{equation}
       f^*: \mathcal{A}_{\alpha}((\pr^{-1}))\rightarrow f(\mathcal{A})_{f(\alpha)}((\pr^{-1})) \ \ \ A\pr^{\alpha+m}\mapsto f(A)\pr^{f(\alpha)+m}
    \end{equation}
where $\mathcal{A}_{f(\alpha)}((\pr^{-1}))$ is given $R$-module structure through $f$. 
\label{definition:base_change_PS_symbol}
\end{definition}

\subsection{$R$-Linear Poisson Vertex Algebras and the Master Theorem}
\label{subsection:PVA_Background}

With the algebra of shifted pseudodifferential symbols in place, we now recall the $R$-linear notion of a Poisson vertex algebra and the master theorem used later to produce $\lambda$-brackets from Adler operators.

\begin{definition}[$R$-Linear $\lambda$-bracket]
An $R$-linear $\lambda$-bracket on an $R$-linear differential algebra $\cV$ is an $R$-linear map $\{\cdot_{\lambda}\cdot\}:\cV\otimes_R\cV\rightarrow \cV[\lambda]$ satisfying for all $f,g,h\in \cV$:
\begin{enumerate}
    \item Sesquilinearity: \begin{equation}
        \{\pr f_{\lambda}g\}=-\lambda\{f_{\lambda}g\}, \ \ \{f_{\lambda}\pr g\}=(\lambda+\pr)\{f_{\lambda}g\}
    \end{equation}
    \item Leibniz Rules: \begin{equation}
        \{f_{\lambda}gh\}=\{f_{\lambda}g\}h+\{f_{\lambda}h\}g, \ \{fh_{\lambda}g\}=\{f_{\lambda+\pr}g\}_{\rightarrow}h+\{h_{\lambda+\pr}g\}_{\rightarrow}f
    \end{equation}
\end{enumerate}
\end{definition}
Following \cite{ADG_SKV} we use the notation:
\begin{notation}
    If $\{f_{\lambda}g\}=\sum_{n\in \Z_{\geq 0}}\lambda^nc_n$, then:
\begin{equation}
    \{f_{\lambda+\pr}g\}_{\rightarrow}h:=\sum_{n\in \Z_{\geq 0}}c_n(\lambda+\pr)^nh
\end{equation}
\end{notation}
\begin{definition}[$R$-Linear Poisson vertex algebra (PVA)]
    An $R$-linear Poisson vertex algebra is an $R$-linear differential algebra $\cV$ endowed with an $R$-linear $\lambda$-bracket $\{\cdot_{\lambda}\cdot\}$, satisfying skew-symmetry and the Jacobi identity for all $f,g,h\in \cV$. 
    \begin{enumerate}
        \item Skew-Symmetry: $\{g_{\lambda}f\}=-\{f_{-\lambda-\pr}g\}$
        \item Jacobi-Identity: $\{f_{\lambda}\{g_{\mu}h\}\}-\{g_{\mu}\{f_{\lambda}h\}\}=\{\{f_{\lambda}g\}_{\lambda+\mu}h\}$
    \end{enumerate}
\end{definition}

\begin{definition}[Free Differential Generators]
Let $(\cV,\pr,\{\cdot_{\lambda}\cdot\})$ be an $R$-linear Poisson vertex algebra with derivation $\pr$. Let $I$ be a set such that for every $i\in I$ we have $u_i\in \cV$. The set $\{u_i:i\in I\}$ is a set of free differential generators if they are algebraically independent over $R[\pr]$, and the $R[\pr]$-algebra generated by them is all of $\cV$. 
\label{definition:differential_gens}
\end{definition}
\begin{definition}[Poisson vertex algebra Ideal]
Let $(\cV, \{\cdot_{\lambda}\cdot\})$ be an $R$-linear PVA. A Poisson vertex algebra ideal $I$ is an ideal of the differential algebra $\cV$ such that $\{I_{\lambda}\cV\}\subset I[\lambda]$
\end{definition}
\begin{definition}[Quotient Poisson vertex algebra]
Let $(\cV, \{\cdot_{\lambda}\cdot\})$ be an $R$-linear PVA and $I$ a PVA ideal. The quotient is the differential algebra $\cV/I$ with $\lambda$-bracket induced by the $\lambda$-bracket on $\cV$. 
\end{definition}
\begin{definition}[R-Linear PVA Filtration]
\label{definition:R_linear_filtration}
A $\Z_{\geq0}$-filtration on an $R$-linear PVA  $(\cV, \{\cdot_{\lambda}\cdot\})$  is a collection of $R$-submodules $F^i\cV$ such that $\cV=\bigcup_{i=0}^{\infty}F^i\cV$, and for all $i,j\in \Z_{\geq 0}$:
\begin{equation}
    F^i\cV\cdot F^j\cV\subset F^{i+j}\cV, \ \ \ \pr(F^i\cV)\subset F^{i+1}\cV \ \ \ (F^i\cV)_{(k)}(F^j\cV)\subset F^{i+j-k-1}\cV
\end{equation}
\end{definition}

\begin{definition}[Opposite Poisson vertex algebra]
Let $(\cV,\{\cdot_{\lambda}\cdot\})$ be an $R$-linear PVA where $\cV$ denotes the underlying $R$-differential algebra. The opposite PVA is defined as $(\cV,\{\cdot_{\lambda}\cdot\})^{op}:=(\cV,-\{\cdot_{\lambda}\cdot\})$ where $-\{\cdot_{\lambda}\cdot\}$ is the opposite bracket.  
\end{definition}

A Poisson anti-isomorphism is an isomorphism from a Poisson vertex algebra to an opposite Poisson vertex algebra.

\begin{definition}[Change of Base of Poisson vertex algebras]
Let $R$ and $S$ be commutative $\mathbb{F}$-algebras, with a morphism $f:R\rightarrow S$. Suppose that $(\cV,\{\cdot_{\lambda}\cdot\})$ is an $R$-linear PVA. The base change along $f$ is defined as:
\begin{equation}
    f^*(\cV,\{\cdot_{\lambda}\cdot\}):=(\cV\otimes_RS,\{\cdot_{\lambda}\cdot\}\otimes_R S)
\end{equation}
and $f^*(\cV,\{\cdot_{\lambda}\cdot\})$ will be an $S$-linear PVA.
\end{definition}
\begin{remark}
    This gives a functor from the category of $R$-linear PVAs to the category of $S$-linear PVAs. 
\end{remark}
\begin{notation}
If $R=\mathbb{F}[T]$ or $\mathbb{F}[T,T^{-1}]$ and $\alpha\in \mathbb{F}$ or $\mathbb{F}^{\times}$, we denote the change of base functor induced by evaluation at $\alpha$ by $\mathrm{Ev}_{T=\alpha}$.  
\end{notation}
$R$-linear PVAs of particular relevance will arise from free commutative differential algebras on a set $I$. More precisely, let $I$ be a set and let $\{u_m\}_{m\in I}$ be a collection of symbols. The differential algebra of differential polynomials over $R$ is defined by:
\begin{equation}
P_I^R = R[\pr][u_m : m \in I].
\end{equation}

Denoting $u_m^{(i)} := \pr^i u_m$, the algebra $P_I^R$ is the polynomial algebra over $R$ generated by the symbols $\{u_m^{(i)} : m \in I,\, i \in \mathbb{Z}_{\geq0}\}$.

We will primarily be interested in endowing $P_I^R$ with a PVA structure in the following two cases: $R=\mathbb{F}$, or $R$ is a localization of $\mathbb{F}[T]$. In the first case, there is a result that allows one to construct Poisson vertex algebra structures on $P_I^R$, known as the Master Theorem \cite[Theorem 1.15]{Master_Thm_Ref}. There is a straightforward extension of the Master Theorem that produces an $R$-linear Poisson vertex algebra structures in the second case.

\begin{theorem}[{\cite[Theorem 1.15]{Master_Thm_Ref}}]

Assume now that $R$ is an integral domain. Let $\cV=P_I^{R}$, and  $H=(H_{nm}(\lambda))\in \mathrm{Mat}_{I\times I}\cV[\lambda]$. 
\begin{enumerate}
    \item There is a unique $R$-linear $\lambda$-bracket $\{\cdot_{\lambda}\cdot\}_H$ on $\cV$, such that $\{u_n{}_{\lambda}u_m\}_H=H_{mn}(\lambda)$ for every $m,n\in I$, and it is given by the following Master Formula for $f,g\in \cV$:
    \begin{equation}
        \{f_{\lambda}g\}_H:=\sum_{\substack{m,n\in I\\ i,j\in \Z_{\geq 0}}}\frac{\pr g}{\pr u_n^{(i)}}(\lambda+\pr)^iH_{nm}(\lambda+\pr)(-\lambda-\pr)^j\frac{\pr f}{\pr u_m^{(j)
        }}
        \label{eqn:master_formula}
    \end{equation}
    \item The $\lambda$-bracket on $\cV$ is skew-symmetric if and only if $H$ is skew adjoint. That is:
    \begin{equation}
        \{u_n{}_{\lambda}u_m\}_H=-\{u_m{}_{-\lambda-\pr}u_n\}_H
    \end{equation}
    \item Assuming the above holds, the $\lambda$-bracket satisfies the Jacobi identity, thus making $\cV$ a $R$-PVA, provided the Jacobi identity holds on any triple of generators $(n,m,k\in I)$:
    \begin{equation}
        \{u_n{}_{\lambda}\{u_m{}_{\mu}u_k\}_H\}_H-\{u_{m}{}_{\mu}\{u_{n}{}_{\lambda}u_k\}_H\}_H=\{\{u_{n}{}_{\lambda}u_m\}_{H_{\lambda+\mu}}u_k\}_H
    \end{equation}
\end{enumerate}

\label{theorem:Master_Theorem}
\end{theorem}
\begin{proof}
    Let $S:=\mathrm{Frac}(R)$ and extend $\cV_{S}:=\cV\ox S=P^{S}_I$. Notice that $H\ox 1\in \mathrm{Mat}_{I\times I}\cV_S[\lambda]$ and will satisfy the axioms of the Master theorem for fields. Hence, there is a Poisson structure on $\cV_{S}$ induced by $H\ox 1$. But $\cV$ is naturally a subset of $\cV_{S}$ and it is clear from the formula given in Equation \eqref{eqn:master_formula} that it is closed under the $\{\cdot_{\lambda}\cdot\}$-bracket. Therefore, $\cV$ is an $R$-linear Poisson vertex algebra. 
\end{proof}
\color{black}
\section{Universal Adler Maps}
\label{section:universal_adler_maps}

As mentioned in Section \ref{section:DS_complex_motivation}, De Sole, Kac, and Valeri, in \cite{ADG_SKV}, translated the second Adler--Gelfand--Dickey bracket into the formalism of Poisson vertex algebras over $\C$. Their construction starts from the Adler map $A^{(L)}$, viewed as a linear map on pseudodifferential operators, and extracts from it a matrix of $\lambda$-brackets, from which the Master formula determines the Poisson vertex algebra structure. 

In this section, we carry out the analogous translation for the generalized second Adler--Gelfand--Dickey bracket of Khesin and Zakharevich. The main point is that the Adler map formalism can be developed for pseudodifferential symbols whose exponents lie in a commutative $\C$-algebra. Applying this over a localization $R$ of $\C[T]$ gives a generic-rank construction, and specialization of $T$ then recovers the corresponding complex-rank Poisson vertex algebras. These specializations are denoted by $\cW(\fr{gl}_{T=\alpha})$ and $\cW(\fr{po}_{T=\alpha})$.

Concretely, in Subsection \ref{subsection:adler_map_technical_results}, we use the formalism developed in Section \ref{section:PVA_background} to extend the Adler map construction in \cite{ADG_SKV} to this $R$-linear setting. The resulting construction is then used to define the $R$-Poisson vertex algebra $\cW(\fr{gl}_T)$ and its specializations in Subsection \ref{subsection:adler_presentation_gl_T}. Similarly, in Subsection \ref{subsection:adler_presentation_po_T}, we define the orthogonal/symplectic variant $\cW(\fr{po}_T)$ and its specializations.

Throughout the rest of the paper, we refer to these Adler--Gelfand--Dickey realizations as Adler presentations.

\subsection{The Adler Map and Technical Results}
\label{subsection:adler_map_technical_results}
Let $R$ be a commutative $\C$-algebra, and $\alpha\in R$. Recall the notation introduced in Equations \eqref{eqn:right_alpha_shifted_notation}, \eqref{eqn:bounded_degree_notation}.
 
The arguments in this subsection closely parallel \cite[Section 2]{ADG_SKV}, so we emphasize only the points that must be modified for interpolation.
\begin{definition}[The Adler Map]
Let $L\in \cA_{\alpha}((\pr^{-1}))$ be a pseudodifferential operator of $\alpha$-shifted degree $0$. The Adler map associated to $L$ is the $R$-linear map:
\begin{equation}
    A^{(L)}:\cA^{-\alpha}((\pr^{-1}))\rightarrow \cA_{\alpha}((\pr^{-1})) 
\end{equation} 
defined through the formula: 
\begin{equation}
    A^{(L)}(F):=(LF)_+L-L(FL)_+=L(FL)_--(LF)_-L
    \label{equation:Generic_Adler_Map}
\end{equation}
\end{definition}

As in \cite{ADG_SKV} the Adler map, when altered appropriately for our setting, is used to derive the matrix of $\lambda$-brackets. Namely, the second part of Equation \eqref{equation:Generic_Adler_Map} implies that: \[A^{(L)}(F)~\in \left(\cA_{\alpha}((\pr^{-1}))\right)_{-1}\]
Additionally, if $F\in \left(\cA^{-\alpha}((\pr^{-1}))\right)_{-1}$, then $(FL)_+=(LF)_+=0\Rightarrow A^{(L)}(F)=0$. Hence, $A^{(L)}$ induces a map:
\begin{equation}
    A^{(L)}:\frac{\cA^{-{\alpha}}((\pr^{-1}))}{\left(\cA^{-\alpha}((\pr^{-1}))\right)_{-1}}\rightarrow \left(\cA_{\alpha}((\pr^{-1}))\right)_{-1}
\end{equation} 
Next notice that:
\begin{equation}
\cA^{-\alpha}((\pr^{-1}))= \left(R[\pr,\pr^{-1}]\pr^{-\alpha}\right)\circ\cA+\left(\cA^{-\alpha}((\pr^{-1}))\right)_{-1}
\end{equation}
\begin{equation}
 \left(\left(R[\pr,\pr^{-1}]\pr^{-\alpha}\right)\circ\cA\right)\cap\left(\cA^{-\alpha}((\pr^{-1}))\right)_{-1}= (R[\pr^{-1}]\pr^{-\alpha})\circ\cA
\end{equation}
Therefore, by the second isomorphism theorem we have that:
\begin{equation}
    \frac{\cA^{-\alpha}((\pr^{-1}))}{\left(\cA^{-\alpha}((\pr^{-1}))\right)_{-1}}\cong \frac{\left(R[\pr,\pr^{-1}]\pr^{-\alpha}\right)\circ\cA}{(R[\pr^{-1}]\pr^{-\alpha})\circ\cA}
\end{equation}
If $I_{\alpha}:=\{-\alpha,-\alpha+1,\cdots \}$, then there is an isomorphism:

\begin{equation}
    \frac{\left(R[\pr,\pr^{-1}]\pr^{-\alpha}\right)\circ\cA}{(R[\pr^{-1}]\pr^{-\alpha})\circ\cA}\xrightarrow{\cong}\cA^{\oplus I_{\alpha}}, \ \ \ \sum_{n=0}^M\pr^{-\alpha+n}\circ F_{-\alpha+n}\mapsto (F_{i})_{i\in I_{\alpha}}
\end{equation}
Similarly, we also have an isomorphism:
\begin{equation}
    (\cA_{\alpha}((\pr^{-1})))_{-1}\xrightarrow{\cong}\cA^{ I_{\alpha}}, \ \ \ \sum_{n=0}^{\infty}P_{-\alpha+n}\pr^{\alpha-1-n}\mapsto (P_i)_{i\in I_{\alpha}}
\end{equation}
Combining these isomorphisms, we see that the Adler map $A^{(L)}$ induces a matrix:
\begin{equation}
H^{(L)}:\cA^{\oplus I_{\alpha}}\rightarrow \cA^{ I_{\alpha}} 
\end{equation}
defined for $f\in \cA$ by:
\begin{equation}
    H_{ij}^{(L)}(\pr)f:=\mathrm{Res}_{\pr}( A^{(L)}(\pr^j\circ f)\circ \pr^i) \quad \text{ for } i,j\in I_{\alpha}
\end{equation}
Notice this parallels the description given in Section \ref{section:DS_complex_motivation}, where the Hamiltonian mapping is induced by the residue pairing. 
We would like to show that $H^{(L)}(\lambda)$ satisfies the hypothesis of Theorem \ref{theorem:Master_Theorem}. To that end one needs a description of the generating function:
\begin{equation}
    H^L(\pr)(z,w):=\sum_{i,j\in I_{\alpha}} H_{ij}^{(L)}(\pr)z^{-i-1}w^{-j-1}=\sum_{i,j\in \Z-\alpha}H_{ij}^{(L)}(\pr)z^{-i-1}w^{-j-1}
\end{equation}
where $H_{ij}(\pr)=0$ if $i$ or $j$ lies in $-\alpha+\Z_{<0}$. 
\begin{lemma}
We have:
\begin{equation}
   H^{(L)}(\pr)(z,w)=L(w)i_w(w-(z+\pr))^{-1}\circ L(z)- L(z+\pr)i_w(w-(z+\pr))^{-1}
\circ L(w+\pr)^*\end{equation}

\label{lemma:matrix_induced_adler_description}
\end{lemma}
\begin{proof}
By definition we have for $f\in \cA$:
\begin{align}
H^{(L)}(\partial)(z,w)f
&=
\sum_{i,j\in I_{\alpha}}
\operatorname{Res}_{\partial}
\left(
A^{(L)}(\partial^{j}\circ f)\partial^i
\right)
z^{-i-1}w^{-j-1}
\notag\\
&=
\operatorname{Res}_{\partial}
\left(
\sum_{i,j\in I_{\alpha}}
A^{(L)}(\partial^jw^{-j-1}\circ f)\,
\partial^iz^{-i-1}
\right)
\notag\\
&=
\operatorname{Res}_{\partial}
\left(
A^{(L)}\bigl(\delta_{-\alpha}(\partial-w)\circ f\bigr)
\delta_{-\alpha}(\partial-z)
\right)
\notag\\
&=
\operatorname{Res}_{\partial}
\Bigl(
\bigl(L(\partial)\bigl(\delta_{-\alpha}(\partial-w)\circ f\bigr)\bigr)_+
L(\partial)\delta_{-\alpha}(\partial-z)
\notag\\
&\hspace{4.5em}
-
L(\partial)
\bigl(\bigl(\delta_{-\alpha}(\partial-w)\circ f\bigr)L(\partial)\bigr)_+
\delta_{-\alpha}(\partial-z)
\Bigr)
\end{align}
By Equations~\eqref{eqn:delta_symmetry_switch_variables}~\eqref{eqn:delta_distribution_identity_change_variables} we have:
\begin{equation}
    (L(\pr)\delta_{-\alpha}(\pr-w)\circ f)_+=L(w)(\delta_0(w-\pr))_+\circ f= L(w)i_w(w-\pr)^{-1}\circ f
\end{equation}
By Lemma \ref{lemma:delta_shift_generating_fct} and Equation~\eqref{eqn:delta_symmetry_switch_variables} we have that:
\begin{align}
    (\delta_{-\alpha}(\pr-w)\circ fL(\pr))_+=(\delta_0(\pr-w)\circ (L(w+\pr)^*f))_+=i_w(w-\pr)^{-1}\circ (L(w+\pr)^*f)
\end{align}
Combining both of these equations we obtain:
\begin{align}
    H^{L}(\pr)(z,w)=\mathrm{Res}_{\pr}\Bigl(L(w)i_w(w-\pr)^{-1}\circ fL(\pr)\delta_{\alpha}(z-\pr)-L(\pr)i_w(w-\pr)^{-1}\circ (L(w+\pr)^*f)\delta_{\alpha}(z-\pr)\Bigr)
\end{align}
Since we are multiplying everything in the residue by $\delta_{-\alpha}(z-\pr)$ this is equivalent to taking the $z$-symbol and so this implies that:

\begin{equation}
    H^L(\pr)(z,w)=\left(L(w)i_w(w-(z+\pr))^{-1}\circ L(z)-L(z+\pr)i_w(w-z-\pr)^{-1}\circ L(w+\pr)^*\right)
    \label{eqn:generating_series_Adler}
\end{equation}

\end{proof}
\begin{definition}[Adler-Type Operator]
    Let $(\cV,\{\cdot_{\lambda}\cdot\})$ be an $R$-linear differential algebra with an $R$-linear $\lambda$-bracket $\{\cdot_{\lambda}\cdot\}$. A pseudodifferential operator of $\alpha$-shifted degree $0$ $L(\pr)\in \cV_{\alpha}((\pr^{-1}))$ is of Adler-type if the following equation holds in $(\cV[\lambda])_{\alpha}((z^{-1},w^{-1}))$:
\begin{equation}
    \{L(z)_{\lambda}L(w)\}=H^{(L)}(\lambda)(z,w)
\end{equation}
\label{defn:adler_type}
\end{definition}
One can use the same arguments as in \cite[Section 2]{ADG_SKV} to obtain:
\begin{prop}[Compare with {\cite[Theorem 2.7]{ADG_SKV}}]
 Let $L\in \cA_{\alpha}((\pr^{-1}))$ be $\alpha$-shifted degree $0$. $L$ induces an $R$-linear Poisson vertex algebra structure on $P_{I_{\alpha}}^R$ defined through the generating series:
\begin{equation}
    \{L(z)_{\lambda}L(w)\}=H^{(L)}(\lambda)(z,w)
\end{equation}
\label{prop:adler_induce_PVA}
\end{prop}

\begin{notation}

    If $(P_{I_{\alpha}}^R,\cdot ,\{\cdot_{\lambda}\cdot\})$ is an $R$-linear Poisson vertex algebra induced from an Adler-type operator $L$, we will indicate this by denoting the $\lambda$-bracket as $\{\cdot_{\lambda}\cdot\}^L$.
\end{notation}
We will need the following technical results on filtrations, evaluations and the inverses of Adler-type operators. 
\begin{lemma}
\label{lemma:Adler_PVA_filtration}
The differential algebra $P_{I_{\alpha}}^R$ can be given a $\Z_{\geq 0}$-filtration $F^{\bullet}P_{I_{\alpha}}^R$ by declaring $\mathrm{deg}(u_{-\alpha+i})=i+1, \deg(\pr)=1$ and extending multiplicatively. If $L(z)=z^{\alpha}+v_{i_1}z^{\alpha-1}+\cdots $ is such that $\mathrm{deg}(v_{i_k})=k$ for all $k\in \Z_{\geq 0}$, then $F^{\bullet}P_{I_{\alpha}}^R$ is a filtration of the $R$-linear PVA $(P_{I_{\alpha}}^R,\{\cdot_{\lambda}\cdot\}^L)$. In particular, for every $i\in \Z_{\geq 0}$ $F^iP_{I_{\alpha}}^R$ is a free $R$-module of finite rank. 
\end{lemma}
\begin{proof}
The only non-trivial step is verifying that $F^{\bullet} P_I^R$ is compatible with $\lambda$-bracket induced by $L$. If $m,n\in \Z_{\geq 0}$, then the $\lambda$-bracket for
$ \{u_{-\alpha+m}{}_{\lambda}u_{-\alpha+n}\}
$  is by definition the coefficient of $z^{\alpha-(m+1)}w^{\alpha-(n+1)}$ of $H^{(L)}(\lambda)(z,w)$. On the other hand, from Lemma \ref{lemma:matrix_induced_adler_description} we know that $H^{(L)}(\pr)(z,w)$ can be written in terms of $L(-), i_w(w-(z+\pr))$. By declaring $\deg(z)=\deg(w)=1$, we see that each term in this series will have total ``degree" $2\alpha-1$, by the assumption on $L$.  In particular, we see that the coefficient:
\begin{equation*}
    (u_{-\alpha+m \ (k)}u_{-\alpha+n}) z^{\alpha-(m+1)}w^{\alpha-(n+1)}\pr^k
\end{equation*} has degree $2\alpha-1$. This implies that $(u_{-\alpha+m \ (k)}u_{-\alpha+n})$ has degree: 
\begin{equation*}2\alpha-1-(\alpha-m-1)-(\alpha-n-1)-k=m+n-k+1=(m+1)+(n+1)-k-1
\end{equation*}
as desired. 
\end{proof}
\begin{lemma}
    If $\alpha\in \mathbb{C}$, then $\mathrm{Ev}_{T=\alpha}(P_I^{\C[T]},\cdot, \{\cdot_{\lambda}\cdot\}^L)$ is a $\C$-PVA that is isomorphic to $(P^{\C}_I,\cdot,\{\cdot_{\lambda}\cdot\}^{\mathrm{Ev}_{T=\alpha}^*(L)})$. A similar statement holds for $\C[T,T^{-1}]$ when $\alpha\in \mathbb{C}^{\times}$. 

\label{lemma:evaluation_functor_PVA_adler_description}
\end{lemma}
\begin{proof}
    This follows immediately from Equation \eqref{eqn:generating_series_Adler} and Definition \ref{defn:adler_type}. 
\end{proof}

\begin{lemma}
\label{lemma:duality_Poisson_algebras}
Let $\cV$ be an $R$-linear PVA with free differential
generators $\{u_{i_{j}}:0\leq j\leq M\}$, where
$M\in \mathbb Z_{\geq 0}\cup\{\infty\}$. Let $\alpha\in R$ and let:
\begin{equation}
    L(\pr)
    =\pr^\alpha+u_{i_0}\pr^{\alpha-1}+u_{i_1}\pr^{\alpha-2}+\cdots\in \cV((\pr^{-1}))
\end{equation}
be a monic $\alpha$-shifted pseudodifferential operator. Then $L(\pr)$ is invertible with inverse: 
\begin{equation}
    L(\pr)^{-1}
    =
    \pr^{-\alpha}
    +h_{i_0}\pr^{-\alpha-1}
    +h_{i_1}\pr^{-\alpha-2}
    +\cdots .
\end{equation}
The coefficients $\{h_{i_j}:0\leq j\leq M\}$ also form a set of free differential generators of
$\cV$. More precisely, for each $0\leq k\leq M$, one has:
\begin{equation}
    h_{i_k}=-u_{i_k}+P_k(u_{i_0},\ldots,u_{i_{k-1}}),
    \label{eqn:trinagular_decomp_diff_gen}
\end{equation}
where $P_k$ is a differential polynomial in the preceding generators. If only a subset of the coefficients of $L(\pr)$ are free differential generators, say $\{u_{i_{j_k}}: 1\leq k\leq M'\}$, then $\{h_{i_{j_k}}:1\leq k\leq M'\}$ is also a set of free differential generators. 

If, in addition, $L(\pr)$ is of Adler-type for $\cV$, then $L(\pr)^{-1}$ is of Adler-type for the opposite PVA structure and the assignment:
\begin{equation}
    \Pi(u_{i_j})=h_{i_j}
    \qquad 0\leq j\leq M,
\end{equation}
extends to an anti-isomorphism:
\begin{equation}
    \Pi:
    (\cV,\cdot,\{\cdot_\lambda\cdot\})
    \longrightarrow
    (\cV,\cdot,\{\cdot_\lambda\cdot\}^{op}).
\end{equation}

\end{lemma}
 \begin{proof}

To prove the first three statements it suffices to show \eqref{eqn:trinagular_decomp_diff_gen}
holds. To see this, we expand $L(\pr)L(\pr)^{-1}$ coefficient wise. For the sake of brevity, we set $u_{i_{-1}}=1, h_{i_{-1}}=1$ so that:
\begin{equation}
    L(\pr)=\sum_{k=0}^{\infty}u_{i_{k-1}}\pr^{\alpha-k} \qquad L(\pr)^{-1}=\sum_{r=0}^{\infty}h_{i_{r-1}}\pr^{-\alpha-r}
\end{equation}
 Expanding we find that:
 \begin{align}
     L(\pr)L(\pr)^{-1}&= (\sum_{k=0}^{\infty}u_{i_{k-1}}\pr^{\alpha-k})(\sum_{r=0}^{\infty}h_{i_{r-1}}\pr^{-\alpha-r})\\
     &= \sum_{k,r=0}^{\infty}u_{i_{k-1}}\pr^{\alpha-k}\circ (h_{i_{r-1}}\pr^{-\alpha-r})\\
     &= \sum_{k,r=0}^{\infty}u_{i_{k-1}}\sum_{t=0}^{\infty} \binom{\alpha-k}{t}h_{i_{r-1}}^{(t)}\pr^{-(k+r+t)}\\
     &=\sum_{N=0}^{\infty}\bigl(\sum_{\substack{k,r, t\geq 0\\ k+t+r=N}} \binom{\alpha-k}{t} u_{i_{k-1}}h_{i_{r-1}}^{(t)}\bigr) \pr^{-N} 
 \end{align}
 This implies that for every $N\geq 0$ we have:
 \begin{equation}
    \sum_{\substack{k,r, t\geq 0\\ k+t+r=N}} \binom{\alpha-k}{t} u_{i_{k-1}}h_{i_{r-1}}^{(t)}=\delta_{N,0} 
 \end{equation}
 The claim then follows by induction, namely for the base case $k=1$, by looking at the degree $N=1$ coefficient we have:
 \begin{equation}
      \binom{\alpha}{0} u_{i_{-1}}h_{i_{0}}+\binom{\alpha}{1} u_{i_{-1}}h_{i_{-1}}^{(1)}+\binom{\alpha-1}{0} u_{i_{0}}h_{i_{-1}}^{(0)}=h_{i_0}+u_{i_0}=0\Rightarrow h_{i_0}=-u_{i_0}
 \end{equation}
 Suppose now it is true for $k_0$. Looking at the degree $N=k_0+1$ coefficient we find:
 \begin{equation}
     \sum_{\substack{k,r, t\geq 0\\ k+t+r=k_0+1}} \binom{\alpha-k}{t} u_{i_{k-1}}h_{i_{r-1}}^{(t)}=h_{i_{k_0}}+u_{i_{k_0}}+\sum_{\substack{t\geq 0, r,k<k_0+1\\ k+t+r=k_0+1}}\binom{\alpha-k}{t} u_{i_{k-1}}h_{i_{r-1}}^{(t)}=0
 \end{equation}
 By the inductive hypothesis we see that indeed $h_{i_{k_0}}$ can be written in the desired form. 
By \cite[Proposition 3.7]{NEW_DVD} $L^{-1}$ will be of opposite Adler-type, from which the anti-isomorphism follows. 
     
\end{proof}
\subsection{The Adler Presentation of $\cW(\fr{gl}_{T})$} 
\label{subsection:adler_presentation_gl_T}
We are now ready to extend the Adler--Gelfand--Dickey algebra to complex and generic rank. Namely, we first define the generic-rank object, which we call the indeterminate Adler--Gelfand--Dickey algebra and denote by $\cW(\fr{gl}_T)$.  By construction, this is the $R$-linear Poisson vertex algebra attached to the universal monic pseudodifferential operator of degree $T$. After specializing $T=\lambda$, this will match the Poisson structure coming from Drinfeld--Sokolov reduction of $\fr{gl}_{\lambda}$ given in \cite{Khesin_Malikov_DS_reduction}.

Recall that $P_{I_T}^{R}=R[u_i^{(n)}:i\in I_T, n\in \Z_{\geq 0}]$, $I_T=\{(-T,T),(-T+1,T),\cdots\}$, and $R$ is a localization of $\C[T]$. Set:
\begin{equation}
    L_T(\pr):=\pr^T+u_{(-T,T)}\pr^{T-1}+u_{(-T+1,T)}\pr^{T-2}+\cdots\in (P_{I_T}^R)_T((\pr^{-1}))
\end{equation}

\begin{definition}[Indeterminate Adler--Gelfand--Dickey Algebra]
    The indeterminate Adler--Gelfand--Dickey algebra is the $R$-linear Poisson vertex algebra $\cW(\fr{gl}_T):=(P^R_{I_T},\cdot,\{\cdot_{\lambda}\cdot\}^{L_T})$.
\end{definition}
By Lemma \ref{lemma:duality_Poisson_algebras} one can invert the $T$-shifted pseudodifferential operator $L_T(\pr)$ to obtain a second set of free differential generators $\{H_{m,T}:m\in \Z_{\geq 1}\}$ given through:
\begin{equation}
    L_T^{-1}(\pr)=\pr^{-T}-H_{1,T}\pr^{-T-1}+H_{2,T}\pr^{-T-2}+\cdots, \ \ \ \ \ \cW(\fr{gl}_T)=R[\pr][H_{i,T}:i\in \Z_{+}] 
\end{equation}

\begin{definition}[The Khesin--Malikov Poisson vertex algebra at $\alpha\in \mathbb{C}$]
For $\alpha\in \C$, let $I_{T=\alpha}:=\{(-\alpha,\alpha),(-\alpha+1,\alpha),\cdots\}$, and let $L_{T=\alpha}$ be the pseudodifferential operator on $P^{\C}_{I_{T=\alpha}}$ defined by:
\begin{equation}
    L_{T=\alpha}:=\pr^{\alpha}+u_{(-\alpha,\alpha)}\pr^{\alpha-1}+\cdots 
\end{equation}
The Khesin--Malikov Poisson vertex algebra at $\alpha\in \C$ is $\cW(\fr{gl}_{T=\alpha}):=(P_{I_{T=\alpha}}^{\C},\{\cdot_{\lambda}\cdot\}^{L_{T=\alpha}})$. 
\end{definition}
As expected, this arises from the evaluation functor at $\alpha\in \C$:
\begin{prop}
Let $\alpha\in \C$. There is an isomorphism of $\C$-PVAs:
\begin{equation}
    \mathrm{Ev}_{T=\alpha}(\cW(\fr{gl}_T))\cong \cW(\fr{gl}_{T=\alpha})
\end{equation}
induced by the differential algebra isomorphism:
\begin{equation}
    u_{(-T+i,T)}\mapsto u_{(-\alpha+i,\alpha)}
\end{equation}
\label{prop:evaluation_functor_UKM_PVA}
\end{prop}
\begin{proof}
    This follows from Lemma \ref{lemma:evaluation_functor_PVA_adler_description} and the fact that $\cW(\fr{gl}_{T=\alpha})$ is isomorphic to the PVA on $P_{I_T}^{\C}$ with $\lambda$-bracket $\{\cdot_{\lambda}\cdot \}^{\mathrm{Ev}_{T=\alpha}^*(L)}$. 
\end{proof}
An interesting feature of Khesin--Malikov PVAs is that for $T=n$ there is a natural surjection to the classical $\cW$-algebra associated to $\fr{gl}_n$. 
\begin{definition}[Adler Presentation of Classical $\cW$-Algebras in Type A]
Let $I_n=\{-n,\cdots, -1\}$ and set $L_n:=\pr^{n}+u_{-n}\pr^{n-1}+\cdots+u_{-1}$. The Adler presentation of $\cW(\fr{gl}_n)$ is the PVA:
\begin{equation}
    \cW^{\mathrm{Adl}}(\fr{gl}_n):=(P_{I_n}^{\C},\cdot, \{\cdot_{\lambda}\cdot\}^{L_n})
\end{equation}
\label{definition:Adler_Presenation_gl_n}
\end{definition}
\begin{prop}[Compare with {\cite[Proposition 2.13]{ADG_SKV}}]

    Let $\alpha\in \C$ and let $(\cV,\{\cdot_{\lambda}\cdot\})$ be a $\C$-Poisson vertex algebra with a monic  $\alpha$-shifted pseudodifferential operator of the form $L(\pr)=\pr^{\alpha}+v_1\pr^{\alpha-1}+\cdots$. If $L$ is of Adler-type, then there is a morphism of $\C$-Poisson vertex algebras $f_{L}:\cW(\fr{gl}_{T=\alpha})\rightarrow \cV$ given by $f_{L}(u_{(-\alpha+(i-1),\alpha)})=v_{i}$ for $i\geq 1$.

    \label{prop:universal_property_Adler_type}
\end{prop}
\begin{proof}
    One only needs to check that this also preserves the $\lambda$-bracket. To that end notice that $f_{L}(L_{T=\alpha}(\pr))=L(\pr)$. As $L$ is Adler-type this means that the $\lambda$-bracket for $\{v_{i_1}{}_{\lambda}v_{i_2}\}$ is built from $L(\pr)$ and so indeed it is preserved. 
\end{proof}

Recall that Feigin's Lie algebras of complex rank satisfy $\fr{gl}_{\lambda}=\fr{gl}_{-\lambda}$ for all $\lambda \in \C$. A similar relation holds for the associated $\cW$-algebras. 
\begin{prop}
    There is an anti-isomorphism of $\C$-Poisson vertex algebras: 
    \begin{equation}
        \Pi_{\alpha}:\cW(\fr{gl}_{T=\alpha})\rightarrow \cW(\fr{gl}_{T=-\alpha}),\quad L_{T=\alpha}(\pr)\mapsto L_{T=-\alpha}(\pr)
    \end{equation}
    \label{prop:anti_iso_complex_rank_Type_A}
\end{prop}
\begin{proof}

    This follows by considering the anti-isomorphism given in Lemma \ref{lemma:duality_Poisson_algebras} and applying Proposition \ref{prop:universal_property_Adler_type}. 
\end{proof}
The next proposition is fundamental for using $\cW(\fr{gl}_T)$ to interpolate the classical $\cW$-algebras. It is the $\cW$-algebra analogue of the asymptotic part of Theorem \ref{theorem:interpolate_segal_sugawara_type_A}, identifying the specialization of the indeterminate generators with the corresponding finite-rank Adler--Gelfand--Dickey generators in bounded degrees.
\begin{prop}

Define $\ssW{n}:P_{I_{T=n}}^\C\rightarrow P_{I_n}^{\C}$ on the free differential generators of $\cW(\fr{gl}_{T=n})$ by:
\begin{equation}
    \ssW{n}(u_{(-n+i,n)}):=\begin{cases}
        u_{-n+i}& \text{ if } 0\leq i\leq n-1 \\
        0 & \text{otherwise}
    \end{cases}
\end{equation}
 $\ssW{n}$ is a surjection of $\C$-Poisson vertex algebras. Furthermore, define a section of differential algebras by:
\begin{equation}
    \ssIW{n}:\cW^{\mathrm{Adl}}(\fr{gl}_n)\rightarrow \cW(\fr{gl}_{T=n}),  \ \ \ \ssIW{n}(u_{-n+i})=u_{(-n+i,n)}
\end{equation}
For fixed $i_0,j_0\in \Z_{\geq 0}$ and sufficiently large $n\in \Z_+$ the following holds:
\begin{equation}
    \ssIW{n}(\{u_{-n+i_0}{}_{\lambda}u_{-n+j_0}\})=\{u_{(-n+i_0,n)}{}_{\lambda}u_{(-n+j_0, n)}\}
\end{equation}
\label{prop:interpolate_gelfand_dickey_projection}
\end{prop}
\begin{proof}
   The first statement is just a special case of Proposition \ref{prop:universal_property_Adler_type}. For the second statement, apply Lemma \ref{lemma:Adler_PVA_filtration} so that $\cW^{\mathrm{Adl}}(\fr{gl}_{T=n})$ is a filtered PVA. In particular, this implies that the $\lambda$-bracket terms of $\{u_{(-T+i_0,T)}{}_{\lambda}u_{(-T+j_0,T)}\}$ will have degree at most $i_0+j_0+1$. Therefore, if $M_0\geq i_0+j_0+1$ the terms will be in the differential sub algebra generated by $u_{(-T,T)},\cdots, u_{(-T+M_0,T)}$. This implies that for $n\geq M_0$ the claim will  hold. 
\end{proof}

\subsection{The Adler Presentation of $\cW(\fr{po}_{T})$}
\label{subsection:adler_presentation_po_T} 
A useful heuristic is to regard $\cW(\fr{gl}_T)$ as the Poisson vertex algebra of functions on the space $\Psi\mathrm{DO}_T$ of pseudodifferential operators of degree $T$.  The subset of self-adjoint pseudodifferential operators should be thought of as a Poisson sub-manifold $\Psi\mathrm{DO}_T^{SA}$ and therefore carry its own induced Poisson structure.  In this subsection we realize the corresponding quotient Poisson vertex algebra and denote it by $\cW(\fr{po}_T)$. 

\begin{lemma}
Let $\alpha\in R$ and $L(\pr)=\pr^{\alpha}+\sum_{k=1}^{\infty}u_{(-\alpha+k-1)}\pr^{\alpha-k}$ be a monic $\alpha$-shifted pseudodifferential operator. If $L(\pr)$ is self-adjoint, then for every $n\in \Z_{\geq 0}$ there exists a differential polynomial $p_n(x_1,\cdots,x_{n})$ such that:
\begin{equation}
    u_{(-\alpha)}=p_0=0
\end{equation}
\begin{equation}
      u_{(-\alpha+2n)}=p_n(u_{(-\alpha+1)},\cdots, u_{(-\alpha+2n-1)})
\end{equation}
\label{lemma:differential_poly_SA_relation}
\end{lemma}
\begin{proof}
For ease of notation, let $u_{(-\alpha-1)}=1$ so that $L(\pr)=\sum_{k=0}^{\infty}u_{(-\alpha+k-1)}\pr^{\alpha-k}$. Recall that by Definition \ref{defn:delta_shifted_adjoint} the adjoint is defined as:
\begin{align}
    L(\pr)^*&=\sum_{k=0}^{\infty}(-1)^k\pr^{\alpha-k}\circ u_{(-\alpha+k-1)}  \\
    &= \sum_{k=0}^{\infty}\left(\sum_{\ell=0}^{\infty}(-1)^k\binom{\alpha-k}{\ell}u^{(\ell)}_{(-\alpha+k-1)}\pr^{\alpha-\ell-k}\right)\\
    &=\sum_{r=0}^{\infty}\left(\sum_{k=0}^r\binom{\alpha-k}{r-k}(-1)^ku_{(-\alpha+k-1)}^{(r-k)}\right)\pr^{\alpha-r}
\end{align}
Therefore, the equation $L(\pr)^*=L(\pr)$ is equivalent to the following equation holding for every $r\in \Z_+$:
\begin{equation}
    u_{(-\alpha+r-1)}=\sum_{k=0}^r\binom{\alpha-k}{r-k}(-1)^ku_{(-\alpha+k-1)}^{(r-k)}
    \label{eqn:Generators_universal_self_adjoint}
\end{equation}
In particular, if $r=2n$ is even then by subtracting $u_{(-\alpha+2n-1)}$ on both sides this is equivalent to the equation:
\begin{equation}
    0=\sum_{k=0}^{2n-1}\binom{\alpha-k}{2n-k}(-1)^ku_{(-\alpha+k-1)}^{(2n-k)}
\end{equation}
Similarly $r=2n+1$ is odd, solving for $u_{(-\alpha+(2n+1)-1)}=u_{(-\alpha+2n)}$ gives:
\begin{equation}
    u_{(-\alpha+2n)}=\frac{1}{2}\sum_{k=0}^{2n}\binom{\alpha-k}{2n+1-k}(-1)^ku_{(-\alpha+k-1)}^{(2n+1-k)}
\end{equation}
By induction we obtain the desired result. 
\end{proof}

One should think of the differential polynomials $p_1,\cdots, p_n,\cdots$ in Lemma \ref{lemma:differential_poly_SA_relation} as the defining equations of the sub-manifold $\Psi\mathrm{DO}_T^{SA}$.

\begin{definition}[Universal Self-Adjoint Khesin--Malikov Poisson vertex algebra]
Let $I_{SA}$ be the $R$-linear Poisson vertex algebra ideal in $\cW(\fr{gl}_T)$ generated by the elements 
\begin{equation}u_{(-T+2n,T)}-p_n(u_{(-T+1,T)},\cdots, u_{(-T+2n-1,T)}) \qquad n\in \Z_{\geq 0}
\end{equation} where the differential polynomials $p_n$ are those appearing in Lemma~\ref{lemma:differential_poly_SA_relation}. The universal self-adjoint Khesin--Malikov Poisson vertex algebra is:
\begin{equation}
    \cW(\fr{po}_T):=\frac{\cW(\fr{gl}_T)}{I_{SA}}
\end{equation}
\end{definition}

\begin{definition}[Universal Self-Adjoint Khesin--Malikov Poisson vertex algebra at $\alpha\in \C$]
For $\alpha\in \C$, let $I_{T=\alpha}:=\{(-\alpha,\alpha),(-\alpha+1,\alpha),\cdots\}$. Let $I_{SA,\alpha}$ be the $\C$-PVA ideal generated by the elements  
\begin{equation}
\mathrm{Ev}_{T=\alpha}(u_{(-T+2n,T)}-p_n(u_{(-T+1,T)},\cdots, u_{(-T+2n-1,T)})) \qquad n\in \Z_{\geq 0}
\end{equation} where the differential polynomial $p_n$ are those appearing in Lemma~\ref{lemma:differential_poly_SA_relation}. The universal self-adjoint Khesin--Malikov PVA at $\alpha\in \C$ is defined as:
\begin{equation}
    \cW(\fr{po}_{T=\alpha}):=\frac{\cW(\fr{gl}_{T=\alpha})}{I_{SA,\alpha}}
\end{equation}
\end{definition}

As expected, this arises from the evaluation functor at $\alpha\in \C$:
\begin{prop}
    Let $\alpha\in \C$. There is an isomorphism of $\C$-PVAs:
    \begin{equation}
        \mathrm{Ev}_{T=\alpha}(\cW(\fr{po}_T))\cong \cW(\fr{po}_{T=\alpha})
    \end{equation}
induced by the isomorphism from Proposition \ref{prop:evaluation_functor_UKM_PVA}.
\end{prop}
\begin{proof}
   To prove the claim, it suffices to show that $I_{SA}$ is flat as a $\C[T]$-module. By Lemma~\ref{lemma:Adler_PVA_filtration}, the Poisson vertex algebra $\cW(\fr{gl}_T)$ carries a filtration $F^{\bullet}\cW(\fr{gl}_T)$ of free finite rank $\C[T]$-modules. The defining relations:
\[u_{(-T+2n,T)} - p_n\bigl(u_{(-T+1,T)},\ldots,u_{(-T+2n-1,T)}\bigr)\]
belong to $F^{2n+1}\cW(\fr{gl}_T)$ for every $n\in \Z_{\geq 0}$. Hence, for each $i\in \Z_{\geq 0}$, the filtered piece: $F^i\cW(\fr{gl}_T)\cap I_{SA}$ is a finitely generated torsion-free $\C[T]$-module, and is therefore free. Since:
\[
    I_{SA}=\varinjlim_i \bigl(F^i\cW(\fr{gl}_T)\cap I_{SA}\bigr)
\]
it follows that $I_{SA}$ is a filtered colimit of free $\C[T]$-modules. Therefore $I_{SA}$ is flat over $\C[T]$.
\end{proof}
The next proposition justifies using the word "universal". 
\begin{prop}
Let $(\cV,\{\cdot_{\lambda}\cdot\})$ be a $\C$-Poisson vertex algebra with a monic pseudodifferential operator $L(\pr)=\pr^{\alpha}+v_1\pr^{\alpha-1}+\cdots$ of Adler-type. If it is self-adjoint, then there is a morphism of $\C$-Poisson vertex algebras $f_{L}:\cW(\fr{po}_{T=\alpha})\rightarrow \cV$  given by $f_{L}([u_{(-\alpha+i-1,\alpha)}])=v_i$, where $[u_{(-\alpha+i-1,\alpha)}]$ denotes the equivalence class of $u_{(-\alpha+i-1,\alpha)}$. 
\label{prop:universal_prop_sa_Khesin_malikov}
\end{prop}
\begin{proof}
By Proposition \ref{prop:universal_property_Adler_type} we see there exists a morphism of $\C$-Poisson vertex algebras $\tilde{f}_L:\cW(\fr{gl}_{T=\alpha})\rightarrow (\cV,\{\cdot_{\lambda}\cdot\})$ such that $\tilde{f}_L(u_{(-\alpha+i-1,\alpha)})=v_i$. On the other hand, we know that $L$ is self-adjoint which means that for $n\in \Z_{\geq 0}$ we have:
\begin{equation}
    v_{2n+1}=p_n(v_2,\cdots, v_{2n})
\end{equation}
where $p_n$ are the differential polynomials appearing in Lemma \ref{lemma:differential_poly_SA_relation}. Therefore, $\tilde{f}_L$ descends to a PVA morphism $f_L:\cW(\fr{po}_{T=\alpha})\rightarrow (\cV,\{\cdot_{\lambda}\cdot\})$. 
\end{proof}

$\cW(\fr{po}_{T})$ and its specializations can alternatively be described through the Adler formalism.
\begin{prop}
Denote the image of $L_T(\pr)$ under the quotient map $\cW(\fr{gl}_T)\rightarrow \cW(\fr{po}_T)$ by:
\begin{equation}
    L_T^{SA}(\pr):=\pr^T+\left(\sum_{r=1}^{\infty}[u_{(-T+2r-1,T)}]\pr^{T-2r}\right)+\left(\sum_{r=1}^{\infty}p_r([u_{(-T+1,T)}],\cdots [u_{(-T+2r-1,T)}])\pr^{T-2r-1}\right)\in (\cW(\fr{po}_T))_T((\pr^{-1}))
\end{equation}
where the $p_r$ are the differential polynomials appearing in Lemma \ref{lemma:differential_poly_SA_relation}. $L_T^{SA}(\pr)$ is of Adler-type in $\cW(\fr{po}_T)$ and the images $\{[u_{(-T+2m-1,T)}]:m\in \Z_+\}$ are free differential generators of $\cW(\fr{po}_T)$.
A similar statement holds for $\cW(\fr{po}_{T=\alpha})$ when specializing $T=\alpha\in \C$. In particular, the corresponding Adler operator $L_{T=\alpha}^{SA}(\pr)$ is obtained from $L_T^{SA}(\pr)$ by replacing all instances of $T$ with $\alpha$. 

\label{prop:differential_algebra_realization_SA_khesin_malikov}
\end{prop}
\begin{proof}
The image $L_T^{SA}(\partial)$ is of Adler-type as the quotient map $\mathcal W(\mathfrak{gl}_T)\to \mathcal W(\mathfrak{po}_T)$ is a morphism of Poisson vertex algebras and $L_T(\partial)$ is of Adler-type.   

We now show that the images of $u_{(-T+2m-1,T)}$ in $\cW(\fr{po}_T)$ are free differential generators by an interpolation argument.

Let $n\in \Z_{+}$ be arbitrary. The $\C$-PVA $\cW^{\mathrm{Mi}}(\fr{sp}_{2n})$, see Section \ref{section:finite_rank_generators_and_isomorphism_W_algebras}, will have free differential generators $\tilde{w}_2,\cdots,\tilde{w}_{2n}$ such that:
\begin{equation}
    (\pr+F_{11})\cdots(\pr+F_{nn})(\pr-F_{nn})\cdots(\pr-F_{11})=\pr^{2n}+\tilde{w}_2\pr^{2n-2}+\tilde{w}_3\pr^{2n-3}+\cdots + \tilde{w}_{2n}
\end{equation}
where $F_{ii}:=E_{i,i}-E_{2n-i+1, 2n-i+1}$ and $E_{i,j}$ are the standard matrices.
Denote the right-hand side by $L_{\Sp_{2n}}(\pr)$. It is easy to verify that:
\begin{equation}
    L_{\Sp_{2n}}(\pr)^*=L_{\Sp_{2n}}(\pr)
\end{equation}
By Proposition \ref{prop:universal_prop_sa_Khesin_malikov} we know there exists a map of $\C$-Poisson vertex algebras $f_{2n}:\cW(\fr{po}_{T=2n})\rightarrow \cW^{\mathrm{Mi}}(\fr{sp}_{2n})$ such that $f_{2n}([u_{(-2n+2i-1,2n)}])=\tilde{w}_{2i}$ for $1\leq i\leq n$. This shows that the set $\{[u_{(-2n+2i-1,2n)}]:1\leq i\leq n\}$ and their derivatives are algebraically independent. This implies that $\{[u_{(-T+2i-1,T)}]: i\in \Z_{+}\}$ and their derivatives are algebraically independent over $\C[T]$. To see this, observe that one may consider $\mathrm{Ev}_{T=\alpha}(\cW(\fr{po}_{T}))$ as a $\C[T]$-module where $T$ acts by $\alpha$. Furthermore, there is a natural $\C[T]$-PVA morphism 
\[e_{\alpha}:\cW(\fr{po}_T)\rightarrow \mathrm{Ev}_{T=\alpha}(\cW(\fr{po}_{T}))\cong \cW(\fr{po}_{T=\alpha}) \quad [u_{(-T+2i-1,T)}]\mapsto [u_{(-\alpha+2i-1,\alpha)}] \] 
given by sending $[u_{(-T+2i-1,T)}]$ to $[u_{(-\alpha+2i-1,\alpha)}]$. Suppose now that the elements $[u_{(-T+2i-1,T)}]$ and their derivatives are not algebraically independent. Then there exists a differential polynomial $q(x_1,\cdots, x_m)$ over $\C[T]$ and $i_1,\cdots, i_m$ such that:
\begin{equation}
    q([u_{(-T+2i_1-1,T)}],\cdots,[u_{(-T+2i_m-1,T)}])=0
\end{equation}
In this case, choose $n_0$ greater than $\max\{i_1,\cdots, i_m\}$ and outside the finite set of roots of the nonzero coefficients of $q$. We see that:
\begin{equation}
       (f_{2n_0}\circ e_{2n_0})(q([u_{(-T+2i_1-1,T)}],\cdots,[u_{(-T+2i_m-1,T)}]))=(\mathrm{Ev}_{T=2n_0}(q))(\tilde{w}_{2i_1},\cdots ,\tilde{w}_{2i_m})=0
\end{equation}
Since we have chosen $n_0$ to avoid the roots of the coefficients of $q$, the specialization $\mathrm{Ev}_{T=2n_0}(q)$ is a nonzero differential polynomial over $\C$, which is a contradiction.

Generation follows from the defining relations of $I_{SA}$, which express every coefficient $[u_{(-T+2n,T)}]$ as a differential polynomial in $[u_{(-T+1,T)}],\ldots,[u_{(-T+2n-1,T)}]$.
\end{proof}
\begin{notation}
    From now on we will abuse notation when talking about the differential generators of $\cW(\fr{po}_T)$ and its specializations by denoting $[u_{(-T+2m-1,T)}]$ simply as $u_{(-T+2m-1,T)}$. Similarly when $T=\alpha$.
\end{notation}
As in the general linear case, Proposition \ref{prop:differential_algebra_realization_SA_khesin_malikov} and Lemma \ref{lemma:duality_Poisson_algebras} imply that there is a second set of free differential generators $\{H_{(2r,T)}^{ev}:r\in \Z_{+}\}$ of $\cW(\fr{po}_T)$ given by the even $-T$-shifted degree coefficients of $L_T^{SA}(\pr)^{-1}$:
\begin{equation}
    L_T^{SA}(\pr)^{-1}=\pr^{-T}+(\sum_{r=1}^{\infty}H_{(2r,T)}^{ev}\pr^{-T-2r})+\cdots
\end{equation}
Recall that $\fr{po}_{\lambda}=\fr{po}_{-\lambda}$. For the same reasons as in Type $A$, the following holds.
\begin{prop}
    There is an anti-isomorphism of Poisson vertex algebras: 
    \begin{equation}
        \Pi_{\alpha}:\cW(\fr{po}_{T=\alpha})\rightarrow \cW(\fr{po}_{T=-\alpha}), \quad L_{T=\alpha}^{SA}(\pr)\mapsto L_{T=-\alpha}^{SA}(\pr)
    \end{equation}
    \label{proposition:anti_iso_complex_PVA_side}
\end{prop}
\begin{proof}
    This follows from Proposition \ref{prop:differential_algebra_realization_SA_khesin_malikov} which allows us to apply Lemma \ref{lemma:duality_Poisson_algebras} to obtain that $\cW(\fr{po}_{T=\alpha})$ is anti-isomorphic to $(\frac{P_{I_{T=\alpha}}^{\C}}{I_{SA,\alpha}},\cdot ,\{\cdot_{\lambda}\cdot\}^{(L_{T=\alpha}^{SA})^{-1}})$ which is isomorphic as a $\C$-PVA to $\cW(\fr{po}_{T=-\alpha})$ due to the universal property of $\cW(\fr{po}_{T=-\alpha})$. 
\end{proof}
We now recall the Adler presentation of the classical $\cW$-algebras in Types $B$ and $C$, which will be used later. 
\begin{definition}[Adler Presentation of Classical $\cW$-algebras in Type $B$]
    Let $I_{2n+1}^{even}:=\{-2n,\cdots, -2\}$ and define a  pseudodifferential operator on $P_{I_{2n+1}^{even}}^{\C}$ by: \begin{equation}L_{\fr{so}_{2n+1}}(\pr):=\pr^{2n+1}+u_{-2n}^B\pr^{2n-1}+p_{1}(u_{-2n}^B)\pr^{2n-2}+\cdots+ p_{n}(u_{-2n}^B,\cdots,u_{-2}^B)\pr^0\end{equation} where $p_i(x_1,\cdots,x_i)$ are the differential polynomials appearing in Lemma \ref{lemma:differential_poly_SA_relation}. We denote the induced PVA structure by $\cW^{\mathrm{Adl}}(\fr{so}_{2n+1})$ and call it the Adler Presentation of Classical $\cW$-algebras in Type $B$.
    \label{defn:Adler_Type_B_Finite}
\end{definition}
\begin{definition}[Adler Presentation of Classical $\cW$-algebras in Type $C$]
    Let $I_{2n}^{even}:=\{-2n+1,\cdots, -1\}$ and define a  pseudodifferential operator on $P_{I_{2n}^{even}}^{\C}$ by: 
    \begin{equation}L_{\fr{sp}_{2n}}(\pr):=\pr^{2n}+u_{-2n+1}^C\pr^{2n-2}+p_{1}(u_{-2n+1}^C)\pr^{2n-3}+\cdots+ u_{-1}^C\pr^0\end{equation}  where $p_i(x_1,\cdots,x_i)$ are the differential polynomials appearing in Lemma \ref{lemma:differential_poly_SA_relation}. We denote the induced PVA structure by $\cW^{\mathrm{Adl}}(\fr{sp}_{2n})$ and call it the Adler Presentation of Classical $\cW$-algebras in Type $C$.
    \label{defn:Adler_Type_C_Finite}
\end{definition}

The next proposition is the orthogonal and symplectic analogue of Proposition \ref{prop:interpolate_gelfand_dickey_projection}.  

\begin{prop}

Define a surjective $\C$-PVA homomorphism $\ssW{2n}:\cW(\fr{po}_{T=2n})\rightarrow \cW^{\mathrm{Adl}}(\fr{sp}_{2n})$ through the assignment:
\begin{equation}
    \ssW{2n}(u_{(-2n+2m-1,2n)}):=\begin{cases}
        u_{-2n+2m-1}^C & \text{if } 1\leq m\leq n \\
        0 & \text{ otherwise}
    \end{cases}
\end{equation}
and similarly define $\ssW{2n+1}:\cW(\fr{po}_{T=2n+1})\rightarrow \cW^{\mathrm{Adl}}(\fr{so}_{2n+1})$ by:
\begin{equation}
    \ssW{2n+1}(u_{(-2n-1+2m-1,2n+1)}):=\begin{cases}
        u_{-2n+2(m-1)}^B & \text{if } 1\leq m\leq n \\
        0 & \text{ otherwise}
    \end{cases}
\end{equation}

These maps admit natural sections as differential algebras. In the symplectic case, define:
\begin{equation}
    \ssIW{2n}:\cW^{\mathrm{Adl}}(\fr{sp}_{2n})\rightarrow \cW(\fr{po}_{T=2n}),  \quad \ssIW{2n}(u_{-2n+2m-1}^C)=u_{(-2n+2m-1,2n)}
\end{equation}
Then, for fixed $i_0,j_0\in \Z_{+}$, for all sufficiently large $n\in \Z_+$ the following holds:
\begin{equation}
    \ssIW{2n}(\{u^C_{-2n+2i_0-1}{}_{\lambda}u^C_{-2n+2j_0-1}\})=\{u_{(-2n+2i_0-1,2n)}{}_{\lambda}u_{(-2n+2j_0-1, 2n)}\}
\end{equation}
\label{prop:interpolate_gelfand_dickey_projection_symplectic_orthogonal}
A similar statement holds for $\fr{so}_{2n+1}$.
\end{prop}
\begin{proof}
  The proof is essentially identical to the proof of Proposition \ref{prop:interpolate_gelfand_dickey_projection}, and so is omitted. 
\end{proof}
\part{Feigin--Frenkel Duality in Complex Rank}

In Parts I and II we constructed two complex-rank PVAs which will now be compared. 
In Part I we defined the Feigin--Frenkel center at the critical level in complex rank:
\begin{equation}
    \fr z(\Lulhat{gl}{\alpha})
\end{equation}
and in Part II we defined the classical $\cW$-algebra at complex rank:
\begin{equation}
    \cW(\fr{gl}_{T=\alpha})
\end{equation}
The goal of Part III is to construct the interpolated Feigin--Frenkel isomorphism at the critical level: 
\[
    \fr f_{T=\alpha}:
    \fr z(\Lulhat{gl}{\alpha})
    \longrightarrow
    \cW(\fr{gl}_{T=\alpha})
\]

The construction of \(\fr f_{T=\alpha}\) is simple at the level of differential algebras. 
In Part I we interpolated Molev's higher Segal--Sugawara vectors to obtain a distinguished set of free differential generators for 
\(\fr z(\Lulhat{gl}{\alpha})\). 
In Part II we constructed a corresponding distinguished set of free differential generators for 
\(\cW(\fr{gl}_{T=\alpha})\). 
We define \(\fr f_{T=\alpha}\) by sending the former generators to the latter. 
This immediately gives an isomorphism of differential algebras. 
The nontrivial point is that this assignment is compatible with the $\lambda$-brackets.

To show that the $\lambda$-brackets are preserved we pass first to finite rank and then interpolate. 
For each $n\in \Z_+$, we need the finite-rank PVA isomorphism:
\begin{equation}
    \fr z(\widehat{\fr{gl}}_n)
    \longrightarrow
    \cW^{\mathrm{Adl}}(\fr{gl}_n)
\end{equation}
where the right-hand side denotes the Adler presentation of the classical \(\cW\)-algebra. 
On the other hand, the Feigin--Frenkel isomorphism at the critical level identifies
\(\fr z(\widehat{\fr{gl}}_n)\) with the Miura presentation, namely the sub-algebra cut out by screening operators:
\begin{equation}
    \widetilde{\fr f}_n:
    \fr z(\widehat{\fr{gl}}_n)
    \longrightarrow
    \cW^{\mathrm{Mi}}(\fr{gl}_n)
\end{equation}
The Adler and Miura presentations are themselves isomorphic via:
\begin{equation}
    \mu_{\fr{gl}_n}:
    \cW^{\mathrm{Adl}}(\fr{gl}_n)
    \longrightarrow
    \cW^{\mathrm{Mi}}(\fr{gl}_n)
\end{equation}
Thus the finite-rank isomorphism relevant for the interpolation argument is
\begin{equation}
    \fr f_n
    :=
    \mu_{\fr{gl}_n}^{-1}\circ \widetilde{\fr f}_n:
    \fr z(\widehat{\fr{gl}}_n)
    \longrightarrow
    \cW^{\mathrm{Adl}}(\fr{gl}_n)
\end{equation}

The main content of Section~\ref{section:finite_rank_generators_and_isomorphism_W_algebras} is to prove that the natural free differential generators of:
\begin{equation}
    \fr z(\widehat{\fr{gl}}_n), \qquad
    \cW^{\mathrm{Adl}}(\fr{gl}_n), \qquad
    \cW^{\mathrm{Mi}}(\fr{gl}_n)
\end{equation}
are carried to one another under these finite-rank isomorphisms. 
This finite-rank compatibility is the input which allows us, in Section~\ref{section:the_finale}, to conclude that the interpolated map
\(\fr f_{T=\alpha}\) is an isomorphism of PVAs.

We do the same thing in Types \(B\) and \(C\). 
The only difference is that the finite-rank Feigin--Frenkel isomorphism involves the Langlands dual Lie algebra with shifted rank. 
Thus, for \(\fr g_n\) of Type \(B\) or \(C\), the Feigin--Frenkel isomorphism is naturally:
\begin{equation}
    \widetilde{\fr f}_{\fr g_n}:
    \fr z(\widehat{\fr g}_n)
    \longrightarrow
    \cW^{\mathrm{Mi}}({}^L\fr g_{n'})
\end{equation}
where $n'=n-1$ in Type $B$ and $n'=n+1$ in Type C.
We compare it with the Adler presentation through:
\begin{equation}
    \mu_{{}^L\fr g_{n'}}:
    \cW^{\mathrm{Adl}}({}^L\fr g_{n'})
    \longrightarrow
    \cW^{\mathrm{Mi}}({}^L\fr g_{n'})
\end{equation}
So in these cases the finite-rank isomorphism is
\begin{equation}
    \fr f_{\fr g_n}
    :=
    \mu_{{}^L\fr g_{n'}}^{-1}\circ \widetilde{\fr f}_{\fr g_n}
\end{equation}
Thus Section~\ref{section:finite_rank_generators_and_isomorphism_W_algebras} has one job: to make the finite-rank comparison of generators explicit. 
Once this is done, Section~\ref{section:the_finale} interpolates this comparison back to complex rank and gives the desired PVA isomorphism.
\section{Finite-Rank Generators and Isomorphisms of \texorpdfstring{$\cW$}{W}}
\label{section:finite_rank_generators_and_isomorphism_W_algebras}

We have now assembled the complex-rank Adler side of the story.  To identify it with the center at the critical level, it remains to recall the finite-rank presentations of the classical $\cW$-algebras and Molev's explicit description of the Feigin--Frenkel isomorphism \cite{Molev}.  These finite-rank formulas are the templates that will be interpolated in the final section. For that purpose, the key results to keep in mind are Lemma \ref{lemma:Feigin_Frenkel_Adler_generator_gl_n} and Lemma \ref{lemma:FF_generators_TYPE_BC}.
\subsection{Preliminaries}
We first collect a small amount of notation from \cite[Chapter 12]{Molev} that will be used uniformly in Types $A$, $B$, and $C$.

\begin{definition}
Let $\fr{g}$ be a Lie algebra over $\C$, with invariant bilinear form $\langle-,-\rangle$. Let $(S(\fr{g}[t^{-1}]t^{-1}),\pr)$ denote the differential algebra with derivation defined on $a\in \fr{g}, n\in \Z_{+}$:
\begin{equation}
    \pr (at^{-n}):=n\cdot  at^{-n-1}
\end{equation}
The PVA $\cV(\fr{g})$ is $(S(\fr{g}[t^{-1}]t^{-1}),\pr)$ with $\lambda$-bracket defined on free differential generators for $a,b\in \fr{g}$ by:
\begin{equation}
    \{at^{-1}_{\lambda}bt^{-1}\}:=[a,b]t^{-1}+\lambda\cdot \langle a,b\rangle
\end{equation}
\end{definition}

\begin{definition}[Complete and Elementary Symmetric Functions]
Let $n\in \Z_{+}$, the complete and elementary symmetric functions in $n$-variables $x_1,\cdots, x_n$ are defined for $m\in \Z_+$ by:
\begin{equation}
    h_m(x_1,\cdots,x_n):=\sum_{i_1\leq \cdots \leq i_m}x_{i_1}\cdots x_{i_m},\qquad 
    e_m(x_1,\cdots,x_n):=\sum_{i_1>\cdots>i_m}x_{i_1}\cdots x_{i_m}
\end{equation}
We make the convention that $h_0=e_0=1$. 
\end{definition}
Recall the well known duality between the complete and elementary symmetric functions:
\begin{equation}
    \sum_{i=0}^m (-1)^i e_{m-i}\cdot h_{i}=\delta_{0,m}
    \label{eqn:complete_elementary_symmetric_duality}
\end{equation}
Often, we will use complete and elementary symmetric functions to obtain differential operators. For example, if $\fr{h}$ is the algebra of diagonal matrices in $\fr{gl}_n$ with elementary matrices $E_{11},\cdots, E_{nn}$ then:
\begin{equation}
    e_m(E_{nn}t^{-1}+\pr,\cdots, E_{11}t^{-1}+\pr) \in \cV(\fr{h})(\pr)
\end{equation}

\subsection{Type $A$}

We begin with Type $A$. 
As mentioned at the start of Part III, the finite-rank map relevant for interpolation is given by composing the usual Feigin--Frenkel isomorphism with an isomorphism comparing the Miura and Adler presentations. 
\begin{equation}
    \fr f_{\fr{gl}_n}
    :=
    \mu_{\fr{gl}_n}^{-1}\circ \widetilde{\fr f}_{\fr{gl}_n}:
    \fr z(\widehat{\fr{gl}}_n)
    \longrightarrow
    \cW^{\mathrm{Adl}}(\fr{gl}_n).
\end{equation}
The Miura presentation is the intermediate object where the comparison takes place. Therefore, we first list a set of generators of the Miura presentation and then describe the isomorphisms $\mu_{\fr{gl}_n}$ and $\tilde{\fr{f}}_{\fr{gl}_n}$ on the level of generators.

\genheading{Generators of the Miura Presentation} The Miura presentation $\cW^{\mathrm{Mi}}(\fr{gl}_n)$ is the Poisson vertex sub-algebra of $\cV(\fr{h})$, where $\fr{h}$ is the Cartan sub-algebra of $\fr{gl}_n$, generated by the coefficients of the differential operator:
    \begin{equation}
        L_{\fr{gl_n}}(\pr):=(\pr+E_{11}t^{-1})\cdots (\pr+E_{nn}t^{-1})=\pr^n+\tilde{w}_1\pr^{n-1}+\cdots+\tilde{w}_n
    \end{equation}
    This is the intersection of kernels of certain screening operators \cite[Theorem 12.4.2]{Molev}.
    Alternatively, the generators $\tilde{w}_1,\cdots, \tilde{w}_n$ can be described through elementary symmetric polynomials. 
    \begin{lemma}[{\cite[Chapter 12]{Molev}}]
        Define for $1\leq i\leq n$:
    \begin{equation}
        \tilde{w}_i=e_i(E_{11}t^{-1}+\pr,\cdots, E_{nn}t^{-1}+\pr)\circ 1
    \end{equation}
    \begin{equation}
        \mathcal{H}_i:= h_i(E_{11}t^{-1}+\pr,\cdots, E_{nn}t^{-1}+\pr)\circ 1
    \end{equation}
    where $\circ $ denotes the product of pseudodifferential operators (see Definition \ref{definition:product_pseudodifferential_operators}).
    In this case, the sets $\{\tilde{w}_1,\cdots, \tilde{w}_n\}$ and $\{\mathcal{H}_1,\cdots, \mathcal{H}_n\}$ are free differential generators.
        \end{lemma}
    Molev's result provides two sets of free differential generators for the Miura presentation. We now show how these two sets are related: the complete symmetric generators occur as the initial coefficients of the inverse operator $L_n(\pr)^{-1}$.
    \begin{lemma}
$L_n(\pr)$ is an invertible pseudodifferential operator with inverse:
        \begin{equation}
            L_n^{-1}(\pr)=\pr^{-n}+
            (-1)\cH_1\pr^{-n-1}+\cdots+(-1)^n\cH_n\pr^{-2n}+\cdots 
        \end{equation}
    \end{lemma}
    \begin{proof}
    As the terms of $L_n^{-1}(\pr)$ are constructed recursively, it suffices to show that the pseudodifferential operator:
    \begin{equation}
        (\pr^n+\tilde{w}_{1}\pr^{n-1}+\cdots+\tilde{w}_{n})(\pr^{-n}+
            (-1)\cH_1\pr^{-n-1}+\cdots+(-1)^n\cH_n\pr^{-2n})
        \label{eqn:lemma_duality_generators_step_1}
    \end{equation}
  has $\pr^0$ coefficient $1$ and the coefficients of $\pr^{-1},\cdots, \pr^{-n}$ are zero. To see this consider the following equation over the indeterminate $z$: 
        \begin{equation}
            \left(\sum_{m=0}^ne_m(E_{11}t^{-1}+\pr,\cdots, E_{nn}t^{-1}+\pr)z^{n-m}\right)\left(\sum_{k=0}^n((-1)^k\cH_k)z^{-n-k}\right)
        \end{equation}
        Simplifying we obtain:
        \begin{align}
            \sum_{t=0}^{2n}\left(\sum_{k=\max(0,t-n)}^{\min(n,t)}(-1)^ke_{t-k}(E_{11}t^{-1}+\pr,\cdots, E_{nn}t^{-1}+\pr)\cH_k\right)z^{-t} 
        \end{align}
        The first $n+1$-terms in this summation simplify to: 
        \begin{align}
           \hspace{2cm} \sum_{t=0}^{n}\left(\sum_{k=0}^t(-1)^ke_{t-k}(E_{11}t^{-1}+\pr,\cdots, E_{nn}t^{-1}+\pr)h_k(E_{11}t^{-1}+\pr,\cdots, E_{nn}t^{-1}+\pr)\circ 1\right)z^{-t} &=\\
            \sum_{t=0}^n\delta_{t,0}z^{-t}&=1
        \end{align}
        where in the last step we used Equation \eqref{eqn:complete_elementary_symmetric_duality}.
        On the other hand, it was pointed out in the proof of \cite[Proposition 12.4.4]{Molev} that:
        \begin{equation}
            \sum_{m=0}^ne_m(E_{11}t^{-1}+\pr,\cdots, E_{nn}t^{-1}+\pr)z^{n-m}=\sum_{m=0}^n\tilde{w}_m(z+\pr)^{n-m}
        \end{equation}
        Therefore, by Lemma~\ref{lemma:symbol_of_product}, the coefficients of \(\pr^0,\pr^{-1},\ldots,\pr^{-n}\) in \eqref{eqn:lemma_duality_generators_step_1} are respectively $1,0,\ldots,0$. This proves the claimed initial terms of \(L_n(\pr)^{-1}\), and the remaining terms are determined recursively.
    \end{proof}

 \genheading{Generators of the Adler Presentation and their images in the Miura Presentation} 
By Definition \ref{definition:Adler_Presenation_gl_n} the Adler presentation $\cW^{\mathrm{Adl}}(\fr{gl}_n)$ has free differential generators $u_{-n},\cdots, u_{-1}$ given as the coefficients of the differential operator:
\begin{equation}
    L_n=\pr^n+u_{-n}\pr^{n-1}+\cdots+u_{-1}
\end{equation}
By Lemma \ref{lemma:duality_Poisson_algebras} $\cW^{\mathrm{Adl}}(\fr{gl}_n)$ has a second set of free differential generators $H_1,\cdots, H_n$ given by some of the coefficients of the inverse pseudodifferential operator:
   \begin{equation}
       L_{_n}^{-1}(\pr)=\pr^{-n}+(-H_1)\pr^{-n-1}+\cdots +((-1)^nH_n)\pr^{-2n}+\cdots
       \label{eqn:alternative_generators_type_A}
    \end{equation}  
   By work of De Sole, Kac, Valeri \cite{ADG_SKV}, we have an isomorphism: $\mu_{\fr{gl}_n}: \cW^{\mathrm{Adl}}(\fr{gl}_n)\rightarrow \cW^{\mathrm{Mi}}(\fr{gl}_n)$/
\begin{theorem}[{\cite[Theorem 2.23]{ADG_SKV}}]

    There is an isomorphism of $\C$-linear Poisson vertex algebras  $\mu_{\fr{gl}_n}:\cW^{\mathrm{Adl}}(\fr{gl}_n)\rightarrow \cW^{\mathrm{Mi}}(\fr{gl}_n)$ given by:
    \begin{equation}
        \mu_{\fr{gl}_n}(L(\pr))=(\pr+E_{11}t^{-1})\cdots (\pr+E_{nn}t^{-1})=\pr^n+\mu_{\fr{gl}_n}(u_{-n})\pr^{n-1}+\cdots +\mu_{\fr{gl}_n}(u_{-1})
    \end{equation}
In particular, for $0\leq i\leq n-1$, 
    $\mu_{\fr{gl}_n}(u_{-n+i})=\tilde{w}_{i+1}$, and $\mu_{\fr{gl}_n}(H_{i+1})= \cH_{i+1}$
\end{theorem}
\genheading{Generators of the Feigin--Frenkel Center and their images in the Miura Presentation} Recall from Section \ref{subsection:molevs_invariants_in_finite_rank} Molev constructed two sets of free differential generators of $\fr{z}(\widehat{\fr{gl}}_n)$, given by $\{\phi_{1,n},\cdots, \phi_{n,n}\}$, $\{\psi_{1,n},\cdots, \psi_{n,n}\}$ each respectively given by the formulas for $1\leq m\leq n$:
   \begin{equation}
       \phi_{m,n}:= \sum_{\mu\vdash m, \ } c_{\mu,m}\cdot Q_{m,\ell}(n)\cdot D^{A}(\mu)
   \end{equation}
   \begin{equation}
      \psi_{m,n}:= \sum_{\mu\vdash m, \ } c_{\mu,m}\cdot (-1)^{\ell+m}Q_{m,\ell}(-n)\cdot P^{A}(\mu)
   \end{equation}
   where $P^A(\mu), D^A(\mu)$ are as described in Equations \eqref{eqn:P_immanant_type_A},\eqref{eqn:D_immanant_type_A} and the polynomial $Q_{m,\ell}(x)$ is as in Equation \eqref{eqn:binomical_coefficient_over_ring}. We now use Molev's finite-rank description of the Feigin--Frenkel isomorphism. 
In our notation this is an isomorphism of Poisson vertex algebras:
\begin{equation}
    \widetilde{\fr f}_{\fr{gl}_n}:
    \fr z(\widehat{\fr{gl}}_n)
    \longrightarrow
    \cW^{\mathrm{Mi}}(\fr{gl}_n).
\end{equation}
We only need its effect on the generators \(\phi_{i,n}\) and \(\psi_{i,n}\).
\begin{prop}[{\cite[Proposition 13.1.3]{Molev}}]
The Feigin--Frenkel isomorphism maps the following free differential generators to one another:
\begin{equation}
    \tilde{\fr{f}}_{\fr{gl}_n}(\phi_{i,n})=\tilde{w}_{i},\qquad 
    \tilde{\fr{f}}_{\fr{gl}_n}(\psi_{i,n})=\cH_{i}
\end{equation}
\end{prop}
In summary, we have found: 
\begin{lemma}
For all $n\in \mathbb{Z}_+$ there is an isomorphism of Poisson vertex algebras $\fr{f}_{\fr{gl}_n}:\fr{z}(\widehat{\fr{gl}}_n)\rightarrow \cW^{\mathrm{Adl}}(\fr{gl}_n)$ such that:
\begin{equation}
     \ \ \ \fr{f}_{\fr{gl}_n}(\phi_{i,n})=u_{-n+i-1}, \ \ \ \ \fr{f}_{\fr{gl}_n}(\psi_{i,n})=H_{i} \ \ \ \text{ for }1\leq i\leq n 
\end{equation}
\label{lemma:Feigin_Frenkel_Adler_generator_gl_n}
\end{lemma}

\subsection{Types $B$ and $C$}
\label{subsection:typesBC_Adler_presentation}
We now turn to the orthogonal and symplectic families.  The overall pattern is parallel to Type $A$, except now the Langlands-dual occurs and there is a rank shift. For type $B$ the finite-rank map relevant for interpolation is:
\begin{equation}
    \fr{f}_{\fr{so}_{2n+1}}:=\mu_{\fr{sp}_{2n}}^{-1}\circ \tilde{\fr{f}}_{\fr{so}_{2n+1}}: \fr{z}(\widehat{\fr{so}}_{2n+1})\rightarrow \cW^{\mathrm{Adl}}(\fr{sp}_{2n})
\end{equation}
and for Type $C$ the relevant finite-rank map is: 
\begin{equation}
    \fr{f}_{\fr{sp}_{2n}}:=\mu_{\fr{so}_{2n+1}}^{-1}\circ \tilde{\fr{f}}_{\fr{sp}_{2n}}: \fr{z}(\widehat{\fr{sp}}_{2n})\rightarrow \cW^{\mathrm{Adl}}(\fr{so}_{2n+1})
\end{equation}
As in the Type $A$ case, the generators of the Miura presentation and understanding where the various isomorphisms send these generators is key. 

\genheading{Generators of the Miura Presentation} 
Let $\fr{g}=\fr{so}_{2n+1}$ or $\fr{sp}_{2n}$. For $1\leq i\leq n$ let $F_{ii}$ be as in Section \ref{subsection:molevs_invariants_in_finite_rank}. Denote $\fr{h}_{\fr{g}}$ to be the Cartan sub-algebra spanned by the $F_{ii}$ in $\fr{g}$, where $\fr{g}=\fr{so}_{2n+1}$ or $\fr{sp}_{2n}$. 

The Miura presentation in Type $B$, denoted by $\cW^{\mathrm{Mi}}(\fr{so}_{2n+1})$, is the Poisson vertex sub-algebra of $\cV(\fr{h}_{\fr{so}_{2n+1}})$ with free differential generators given by the coefficients of the odd-powers of $\pr$ \cite[Theorem 12.2.3]{Molev}:
    \begin{equation}
    \begin{split}
    L_{n}^B(\pr):&=(\pr+F_{11})\cdots(\pr+F_{nn})\pr(\pr-F_{nn})\cdots(\pr-F_{11}) \\&=\pr^{2n+1}+\tilde{w}_2^B\pr^{2n-1}+\tilde{w}_3^B\pr^{2n-2}+\cdots+ \tilde{w}_{2n+1}^B
    \end{split}
    \end{equation}
Alternatively, the generators $\tilde{w}_2^B,\cdots, \tilde{w}_{2n}^B$ can be described through elementary symmetric polynomials. 
\begin{lemma}[{\cite[Chapter 12]{Molev}}]
Define for $1\leq j\leq 2n$:
\begin{equation}
    \tilde{w}_{j}^B:=e_j(F_{11}t^{-1}+\pr,\cdots,F_{nn}t^{-1}+\pr,\pr,-F_{nn}t^{-1}+\pr,\cdots, -F_{11}t^{-1}+\pr)\circ 1
\end{equation}
\begin{equation}
    \mathcal{H}_{j}^B:=h_j(F_{11}t^{-1}+\pr,\cdots,F_{nn}t^{-1}+\pr,\pr,-F_{nn}t^{-1}+\pr,\cdots, -F_{11}t^{-1}+\pr)\circ 1
\end{equation}
where $\circ $ denotes the product of pseudodifferential operators (see Definition \ref{definition:product_pseudodifferential_operators}).
    In this case, the sets $\{\tilde{w}_2^B,\cdots, \tilde{w}_{2n}^B\}$ and $\{\mathcal{H}_2^B,\cdots, \mathcal{H}_{2n}^B\}$ are free differential generators of $\cW^{\mathrm{Mi}}(\fr{so}_{2n+1})$. 
\end{lemma}
As in the Type $A$ case, the second set of generators comes from certain coefficients of ($L_{n}^B(\pr))^{-1}$.
\begin{lemma}
    $L_{n}^B(\pr)$ is an invertible pseudodifferential operator with inverse
    \begin{equation}
        L_n^B(\pr)^{-1}=\pr^{-2n-1}+\cH_2^B\pr^{-2n-3}+\cdots + (-1)^{2n+1}\cH_{2n+1}^B\pr^{-4n-2}+\cdots 
    \end{equation}
\end{lemma}
\begin{proof}
    The proof is almost identical to the Type $A$ case, and so it is omitted. 
\end{proof}

The Miura presentation in Type $C$, denoted by $\cW^{\mathrm{Mi}}(\fr{sp}_{2n})$, is the Poisson vertex sub-algebra of $\cV(\fr{h}_{\fr{sp}_{2n}})$ with free differential generators given by the even-degree coefficients of $\pr$ \cite[Theorem 12.2.6]{Molev}:
    \begin{equation}
    \begin{split}
        L_n^C(\pr):=(\pr+F_{11})\cdots(\pr+F_{nn})(\pr-F_{nn})\cdots (\pr-F_{11})=\pr^{2n}+\tilde{w}_2^C\pr^{2n-2}+\tilde{w}_3^C\pr^{2n-3}+\cdots+\tilde{w}_{2n}^C
    \end{split}
    \end{equation}
    Alternatively, the free differential generators $\tilde{w}_2^C,\cdots, \tilde{w}_{2n}^C$ can be described through elementary symmetric polynomials. 
\begin{lemma}[{\cite[Chapter 12]{Molev}}]
Define for $1\leq j\leq 2n$:
\begin{equation}
    \tilde{w}_{j}^C:=e_j(F_{11}t^{-1}+\pr,\cdots,F_{nn}t^{-1}+\pr,-F_{nn}t^{-1}+\pr,\cdots, -F_{11}t^{-1}+\pr)\circ 1
\end{equation}
\begin{equation}
    \mathcal{H}_{j}^C:=h_j(F_{11}t^{-1}+\pr,\cdots,F_{nn}t^{-1}+\pr,-F_{nn}t^{-1}+\pr,\cdots, -F_{11}t^{-1}+\pr)\circ 1
\end{equation}
where $\circ $ denotes the product of pseudodifferential operators (see Definition \ref{definition:product_pseudodifferential_operators}).
    In this case, the sets $\{\tilde{w}_2^C,\cdots, \tilde{w}_{2n}^C\}$ and $\{\mathcal{H}_2^C,\cdots, \mathcal{H}_{2n}^C\}$ are free differential generators of $\cW^{\mathrm{Mi}}(\fr{sp}_{2n})$. 
\end{lemma}
As in the Type $A$ case, the second set of generators comes from certain coefficients of ($L_{n}^C(\pr))^{-1}$.
\begin{lemma}
    $L_{n}^C(\pr)$ is an invertible pseudodifferential operator with inverse
    \begin{equation}
        L_{n}^C(\pr)^{-1}=\pr^{-2n}+\cH_2^C\pr^{-2n-2}+\cdots + (-1)^{2n+1}\cH_{2n+1}^C\pr^{-4n-1}+\cdots 
    \end{equation}
\end{lemma}
\begin{proof}
    The proof is almost identical to the Type $A$ case, and so it is omitted. 
\end{proof}
\genheading{Generators of the Adler Presentation and their images in the Miura Presentation} 

By Definition~\ref{defn:Adler_Type_B_Finite} the Adler presentation in Type $B$, denoted by$\cW^{\mathrm{Adl}}(\fr{so}_{2n+1})$, has free differential generators $u_{-2n}^B,\cdots, u_{-2}^B$ which are packaged together through a differential operator:
\begin{equation}
L_{\fr{so}_{2n+1}}(\pr):=\pr^{2n+1}+u_{-2n}^B\pr^{2n-1}+p_{1}(u_{-2n})\pr^{2n-2}+\cdots+ p_{n}(u_{-2n},\cdots,u_{-2})\pr^0
\end{equation}
where $p_i(x_1,\cdots,x_i)$ are the differential polynomials from Lemma \ref{lemma:differential_poly_SA_relation}. By Lemma \ref{lemma:duality_Poisson_algebras} $\cW^{\mathrm{Adl}}(\fr{so}_{2n+1})$ has a second set of free differential generators $H_{2}^B,\cdots, H_{2n}^B$ given by some of the coefficients of the inverse pseudodifferential operator:
\begin{equation}
    L_{\fr{so}_{2n+1}}^{-1}(\pr)=\pr^{-2n-1}+(H_2^B)\pr^{-2n-3}+\cdots+(H_{2n}^B)\pr^{-4n-1}+\cdots 
\end{equation}

By Definition \ref{defn:Adler_Type_C_Finite} the Adler presentation, in Type $C$, denoted by $\cW^{\mathrm{Adl}}(\fr{sp}_{2n})$, has free differential generators $u_{-2n+1}^C,\cdots, u_{-1}^C$ which are packaged together through a differential operator:
\begin{equation}
    L_{\fr{sp}_{2n}}(\pr)=\pr^{2n}+u_{-2n+1}^C\pr^{2n-2}+\cdots+u_{-1}^C\pr^0
\end{equation}
By Lemma \ref{lemma:duality_Poisson_algebras} $\cW^{\mathrm{Adl}}(\fr{sp}_{2n})$ has a second set of free differential generators $H_2^C,\cdots, H_{2n}^C$ given by some of the coefficients:
\begin{equation}
    L^{-1}_{\fr{sp}_{2n}}(\pr)=\pr^{-2n}+H_{2}^C\pr^{-2n-2}+\cdots+ H_{2n}^C\pr^{-4n}+\cdots
\end{equation}

Both sets of free differential generators in Type $B$ and Type $C$ are sent to the corresponding Miura presentation through the following isomorphism. 
\begin{prop}

 Let $\fr{g}_N$ denote $\fr{so}_{2n+1}$ in the Type $B$ case, $\fr{sp}_{2n}$ in the Type $C$ case. There is an isomorphism of PVAs $\mu_{\fr{g}_N}:\cW^{\mathrm{Adl}}(\fr{g}_N)\rightarrow \cW^{\mathrm{Mi}}(\fr{g}_N)$ given in Types $B$ and $C$ respectively by:
\begin{equation}
    \mu_{\fr{g}_N}(L_{\fr{so}_{2n+1}}(\pr))=(\pr+F_{11})\cdots(\pr+F_{nn})\pr(\pr-F_{nn})\cdots(\pr-F_{11}) 
\end{equation}
\begin{equation}
    \mu_{\fr{g}_N}(L_{\fr{sp}_{2n}})=(\pr+F_{11})\cdots(\pr+F_{nn})(\pr-F_{nn})\cdots (\pr-F_{11})
\end{equation}
In particular, for $1\leq i\leq n$:
\begin{equation}
    \mu_{\fr{g}_N}(u_{-N+2i-1}^{B/C})=\tilde{w}_{2i}^{B/C}, \ \  \ \mu_{\fr{g}_N}(H^{B/C}_{2i})=\cH_{2i}^{B/C}
\end{equation}
\end{prop}
\begin{proof}
   The differential operators on the left-hand side are self-adjoint in the sense of Definition~\ref{defn:delta_shifted_adjoint}. Therefore, by Proposition \ref{prop:universal_prop_sa_Khesin_malikov}, there is a surjection $\cW(\fr{po}_{T=N})\rightarrow \cW^{\mathrm{Mi}}(\fr{g}_N)$ of PVAs. This morphism factors through to a morphism $\mu_{\fr{g}_N}:\cW^{\mathrm{Adl}}(\fr{g}_N)\rightarrow \cW^{\mathrm{Mi}}(\fr{g}_N)$. Since $\mu_{\fr{g}_N}$ maps free differential generators to the free differential generators of $\cW^{\mathrm{Mi}}(\fr{g}_N)$, it is an isomorphism.
\end{proof}
\genheading{Generators of the Feigin--Frenkel Center and their images in the Miura Presentation} Recall from Section \ref{subsection:molevs_invariants_in_finite_rank} that Molev constructed a set of free differential generators $\{\phi^B_{2,2n+1},\cdots, \phi_{2n,2n+1}^B\}$ of $\fr{z}(\widehat{\fr{so}}_{2n+1})$ and a set of free differential generators $\{\phi_{2,2n}^C,\cdots, \phi_{2n,2n}^C\}$ of $\fr{z}(\widehat{\fr{sp}}_{2n})$ given by \cite[Theorem 2.3]{Molev_Type_B_C_Invaraints}: 
   \begin{equation}\phi_{m,2n+1}^{B}:=\sum_{\substack{\mu\vdash m \\\ even}}c_{\mu,m}\cdot Q_{m,\ell}(-2n)\cdot P^B(\mu)\end{equation}
   \begin{equation}\phi_{m,2n}^{C}:=\sum_{\substack{\mu \vdash m, \\  \ even}} c_{\mu,m}\cdot Q_{m,\ell}(2n+1)\cdot D^C(\mu)\end{equation} 
 where $P^B(\mu), D^C(\mu)$ are as in Equations \eqref{eqn:P_immanant_type_B} \eqref{eqn:D_immanant_type_C} and the polynomial $Q_{m,\ell}(x)$ is as in Equation \eqref{eqn:binomical_coefficient_over_ring}. We now use Molev's finite-rank description of the Feigin--Frenkel isomorphism in Types $B$ and $C$. In our notation these are the isomorphisms:
\begin{equation}
    \widetilde{\fr f}_{\fr{so}_{2n+1}}:
    \fr z(\widehat{\fr{so}}_{2n+1})
    \longrightarrow
    \cW^{\mathrm{Mi}}(\fr{sp}_{2n}) \qquad \widetilde{\fr f}_{\fr{sp}_{2n}}:
    \fr z(\widehat{\fr{sp}}_{2n})
    \longrightarrow
    \cW^{\mathrm{Mi}}(\fr{so}_{2n+1})
\end{equation}
We only need their effects on the free differential generators:

\begin{prop}[{\cite[Propositions 13.1.7, 13.1.14]{Molev}}]
    The Feigin--Frenkel isomorphism will map the following free differential generators to one another for $1\leq i\leq n$:
    \begin{equation}
        \tilde{\fr{f}}_{\fr{so}_{2n+1}}^B(\phi_{2i,2n+1}^B)=\cH_{2i}^C ,\qquad 
        \tilde{\fr{f}}_{\fr{sp}_{2n}}^C(\phi_{2i,2n}^C)=\tilde{w}_{2i}^B
    \end{equation}
\label{prop:FF_type_B}
\end{prop}

In summary, we have found:
\begin{lemma}
   Let $\fr{g}=\fr{so}_{2n+1}$ or $\fr{sp}_{2n}$ for $n\in \Z_+$. There are isomorphisms of Poisson vertex algebras: 
   \begin{equation}
   \fr{f}_{\fr{so}_{2n+1}}^{B}:\fr{z}(\widehat{\fr{so}}_{2n+1})\rightarrow \cW^{\mathrm{Adl}}(\fr{sp}_{2n}),\quad \fr{f}_{\fr{sp}_{2n}}^{C}:\fr{z}(\widehat{\fr{sp}}_{2n})\rightarrow \cW^{\mathrm{Adl}}(\fr{so}_{2n+1})
   \end{equation} that acts on generators $1\le i\leq n$:
   \begin{equation}
       \fr{f}^B_{\fr{so}_{2n+1}}(\phi^B_{2i,2n+1})= H_{2i}^{C},
\qquad 
       \fr{f}^C_{\fr{sp}_{2n}}(\phi^C_{2i,2n})= u_{-2n+2(i-1)}^{B}
   \end{equation}
\label{lemma:FF_generators_TYPE_BC}
\end{lemma}
\section{Interpolating the Feigin--Frenkel Isomorphism at the Critical Level}
\label{section:the_finale}

As outlined at the start of Part III, we now use the finite-rank description from the previous section to extend the Feigin--Frenkel isomorphism at the critical level to complex rank by interpolating. 
\subsection{Type A}
We begin with Type $A$, where the generator matching is simplest. On the affine side, by Theorem \ref{theorem:interpolate_segal_sugawara_type_A} we have two sets of free differential generators for the interpolated Feigin--Frenkel center at the critical level over $R=\C[T,T^{-1}]$:
\begin{equation}
    \fr{z}(\Lulhat{gl}{T})=R[\mathsf{T}][\phi_{i,T}:i\in\Z_+]=R[\mathsf{T}][\psi_{i,T}:i\in \Z_+]
\end{equation}
On the $\cW$-algebra side, by Lemma \ref{lemma:duality_Poisson_algebras} we have two sets of free differential generators of the indeterminate Adler--Gelfand--Dickey algebra $\cW(\fr{gl}_T)$:
\begin{equation}
    \cW(\fr{gl}_T)=R[\pr][u_{(-T+(i-1),T)}:i\in \Z_+]=R[\pr][H_{i,T}:i\in \Z_+]
\end{equation}
By Lemma \ref{lemma:Feigin_Frenkel_Adler_generator_gl_n}, the finite-rank Feigin--Frenkel isomorphism sends the tuple of free differential generators $(\phi_{1,n},\cdots, \phi_{n,n})$ exactly to the tuple $(u_{-n},\cdots, u_{-1})$. Using this, we can naturally extend the Feigin--Frenkel isomorphism at the critical level to complex rank in Type $A$. 
\begin{theorem}[Interpolated Feigin--Frenkel Isomorphism at the Critical Level: Type $A$]
\label{theorem:FF_inter_type_A}

There is an isomorphism of $\C[T,T^{-1}]$-Poisson vertex algebras given by:
    \begin{equation}\fr{f}_T:\fr{z}(\Lulhat{gl}{T})\rightarrow  \cW(\fr{gl}_{T}), \qquad \fr{f}_T(\phi_{i,T})=u_{(-T+i-1,T)}  \quad i\geq 1
    \end{equation}
Specializing to $\alpha\in \C^{\times}$ induces an isomorphism of $\C$-Poisson vertex algebras:
\begin{equation}
    \fr{f}_{T=\alpha}:\fr{z}(\Lulhat{gl}{\alpha})\rightarrow \cW(\fr{gl}_{T=\alpha}), \qquad \fr{f}_{T=\alpha}(\phi_{i,T=\alpha})=u_{(-\alpha+i-1,\alpha)}  \quad i\geq 1
\end{equation}
These isomorphisms interpolate the usual Feigin--Frenkel isomorphisms in the sense that Diagram \ref{diagram:interpolating_FF_type_A} commutes. 
\begin{figure}[H]
\[\begin{tikzcd}
	{\fr{z}(\Lulhat{gl}{n})} && {\cW(\fr{gl}_{T=n})} \\
	\\
	{\fr{z}(\widehat{\fr{gl}}_{n})} && {\cW^{\mathrm{Adl}}(\fr{gl}_n)}
	\arrow["{\fr{f}_{T=n}}", from=1-1, to=1-3]
	\arrow["{\ssV{n}}"', from=1-1, to=3-1]
	\arrow["{\ssW{n}}", from=1-3, to=3-3]
	\arrow["{\fr{f}_{\fr{gl}_n}}"', from=3-1, to=3-3]
\end{tikzcd}\]
\caption{}
\label{diagram:interpolating_FF_type_A}
\end{figure}
\end{theorem}
\begin{proof}
    By the discussion preceding the theorem,  $\fr{f}_{T}$ is an isomorphism of differential algebras over $\C[T,T^{-1}]$. Similarly, after specializing $T=\alpha\in \C^{\times}$, $\fr{f}_{T=\alpha}$ is an isomorphism of differential algebras over $\C$. Moreover, by Theorem \ref{theorem:interpolate_segal_sugawara_type_A} and Lemma \ref{lemma:evaluation_functor_PVA_adler_description}, evaluation is compatible with the Poisson vertex algebra structures:
\[
\mathrm{Ev}_{T=\alpha}(\fr{z}(\Lulhat{gl}{T}))\cong \fr{z}(\Lulhat{gl}{\alpha}),
\qquad
\mathrm{Ev}_{T=\alpha}(\cW(\fr{gl}_T))\cong \cW(\fr{gl}_{T=\alpha}).
\]
Thus, once $\fr{f}_T$ is shown to preserve the $\lambda$-bracket over $\C[T,T^{-1}]$, its specialization $\fr{f}_{T=\alpha}$ is automatically an isomorphism of $\C$-PVAs under these identifications. It remains to prove bracket preservation for $\fr{f}_T$. Fix $i,j\in \Z_+$. By the asymptotic component of Theorem \ref{theorem:interpolate_segal_sugawara_type_A}, there exists $N_{(i,j)}\in \Z_{\geq 0}$ such that for  $N'\geq N_{(i,j)}$ there is a  differential algebra homomorphism $\ssIV{N'}:\fr{z}(\widehat{\fr{gl}}_{N'})\rightarrow \fr{z}(\Lulhat{gl}{N'})$ for which the following holds:
    \begin{equation}
        \mathrm{Ev}_{T=N'}(\{\phi_{i,T}{}_{\lambda}\phi_{j,T}\})=\{\phi_{i,T=N'}{}_{\lambda}\phi_{j,T=N'}\}= \ssIV{N'}(\{\phi_{i,N'}{}_{\lambda}\phi_{j,N'}\})
        \label{equation:step_1_type_A_Finale_equation}
    \end{equation} 
    By the asymptotic component of Proposition \ref{prop:interpolate_gelfand_dickey_projection}, there exists $M_{(i,j)}\in \Z_{\geq 0}$ such that for $N''\geq M_{(i,j)}$ there is an injective differential algebra homomorphism $\ssIW{N''}:\cW^{\mathrm{Adl}}(\fr{gl}_{N''})\rightarrow \cW(\fr{gl}_{T=N''})$ for which the following holds:
    \begin{equation}
        \mathrm{Ev}_{T=N''}(\{u_{(-T+i-1,T)}{}_{\lambda}u_{(-T+j-1,T)}\})=
 \{u_{(-N''+i-1,N'')}{}_{\lambda}u_{(-N''+j-1,N'')}\}=\ssIW{N''}(\{u_{-N''+i-1}{}_{\lambda}u_{-N''+j-1}\})
    \end{equation}
Set $N_0= \max\{N_{(i,j)}, M_{(i,j)}\}$, and notice that for $N\geq N_0$:
\begin{equation}
    \fr{f}_{T=N}\circ \ssIV{N}=\ssIW{N}\circ\fr{f}_{\fr{gl}_N}
    \label{equation:step_2_type_A_Finale_equation}
\end{equation}
as it holds for the free differential generators of $\fr{z}(\widehat{\fr{gl}}_N)$. Consider: 
\begin{equation}
\begin{aligned}
\Delta_{i,j}(T)&:= \fr{f}_T\left(\{\phi_{i,T}{}_{\lambda}\phi_{j,T}\}\right)-(\{\fr{f}_T(\phi_{i,T}){}_{\lambda}\fr{f}_T(\phi_{j,T})\})\\
& \ =\fr{f}_T\left(\{\phi_{i,T}{}_{\lambda}\phi_{j,T}\}\right)-(\{u_{(-T+i-1,T)}{}_{\lambda}u_{(-T+j-1,T)}\})
 \end{aligned}
 \end{equation}
By Theorem \ref{theorem:invariant_ind_PVA_structure}, the \(\lambda\)-coefficients of
\(\Delta_{i,j}(T)\) are Laurent polynomials in \(T\). Thus, to show that $\Delta_{i,j}(T)$ is zero, it suffices to show that they vanish at sufficiently large integers. 
Equations~\eqref{equation:step_1_type_A_Finale_equation}, \eqref{equation:step_2_type_A_Finale_equation} imply: 
\begin{equation}
\begin{aligned}
    \mathrm{Ev}_{T=N}(\Delta_{i,j}(T))
    &=
    (\fr{f}_{T=N}\circ \ssIV{N})(\{\phi_{i,N}{}_{\lambda}\phi_{j,N}\})
    -
   (\{u_{(-N+i-1,N)}{}_{\lambda}u_{(-N+j-1,N)}\})  \\
    &=
    (\ssIW{N}\circ \fr{f}_{\fr{gl}_N})(\{\phi_{i,N}{}_{\lambda}\phi_{j,N}\})
    -
    \ssIW{N} (\{u_{-N+i-1}{}_{\lambda}u_{-N+j-1}\}) \\
    &=0,
\end{aligned}
\end{equation}
    Hence, the coefficients are zero and so $\fr{f}_T\left(\{\phi_{i,T}{}_{\lambda}\phi_{j,T}\}\right)=\{\fr{f}_T(\phi_{i,T}){}_{\lambda}\fr{f}_T(\phi_{j,T})\}$. As $i,j$ were arbitrary, this implies that $\fr{f}_T$ is an isomorphism of $\C[T,T^{-1}]$-PVAs.

    The last claim follows by computing $(\fr{f}_{\fr{gl}_n}\circ \ssV{n})(\phi_{i,T=n})$ for $1\leq i\leq n$ through Theorem \ref{theorem:interpolate_segal_sugawara_type_A}, Lemma \ref{lemma:Feigin_Frenkel_Adler_generator_gl_n}, and computing $(\ssW{n}\circ \fr{f}_{T=n})(\phi_{i,T=n})$ through Proposition \ref{prop:interpolate_gelfand_dickey_projection}. 
\end{proof}
Notice that we have two sets of free differential generators for the interpolated Feigin--Frenkel center $\fr{z}(\Lulhat{gl}{T})$ and two sets of free differential generators for the indeterminate Adler--Gelfand--Dickey algebra $\cW(\fr{gl}_T)$. In each case these two sets of generators arise from a duality induced by the interpolated parity functor, see Corollary \ref{cor:anti_isomorphism_centers_generators} and Proposition \ref{prop:anti_iso_complex_rank_Type_A}. The following result shows these dualities are compatible with the interpolated Feigin--Frenkel isomorphism. Interestingly, this is something that can only be seen in the interpolated setting.
\begin{theorem}
\label{theorem:FF_inter_type_A_parity}
    The interpolated Feigin--Frenkel isomorphism $\fr{f}_T:\fr{z}(\Lulhat{gl}{T})\rightarrow \cW(\fr{gl}_{T})$ commutes with the composition of the interpolated parity functor and Cartan morphism. More precisely, Diagram \ref{diagram:parity_functor_type_A} commutes. A similar result holds when $T$ is specialized to $\alpha\in \C^{\times}$.
    \begin{figure}[H]
    \[\begin{tikzcd}
	{\mathfrak{z}(\Lulhat{gl}{T})} & {\cW(\fr{gl}_{T})} \\
	{\mathfrak{z}(\Lulhat{gl}{-T})} & {\cW(\fr{gl}_{-T})}
	\arrow["{\mathfrak{f}_T}", from=1-1, to=1-2]
	\arrow["\Pi_{T}", from=1-2, to=2-2]
	\arrow[" \Pi_{T}\circ \nu_T", from=1-1, to=2-1]
	\arrow["{\mathfrak{f}_{-T}}", from=2-1, to=2-2]
\end{tikzcd}\]
\caption{}
\label{diagram:parity_functor_type_A}
\end{figure}

\end{theorem}
\begin{proof}
It suffices to calculate $\fr{f}_{-T}\circ (\Pi_T\circ \nu_T)$ on the generators $\phi_{i,T}$. By Corollary \ref{cor:anti_isomorphism_centers_generators} and by interpolating Lemma \ref{lemma:Feigin_Frenkel_Adler_generator_gl_n}:
\begin{equation}
    (\fr{f}_{-T}\circ ( \Pi_T\circ \nu_T))(\phi_{i,T})=(-1)^i\fr{f}_{-T}(\psi_{i,-T})=(-1)^iH_{i,-T}
\end{equation}
On the other hand we have by Theorem \ref{theorem:FF_inter_type_A} and Lemma \ref{lemma:duality_Poisson_algebras}:
\begin{equation}
    (\Pi_{T}\circ \fr{f}_T)(\phi_{i,T})=\Pi_T(u_{(-T+i-1,T)})=(-1)^iH_{i,-T}
\end{equation}
The specialization statement follows by evaluating $T=\alpha$.
\end{proof}
\subsection{Types $B$ and $C$}

For Types $B$ and $C$, the generator matching is slightly more involved. On the affine side, by Theorem \ref{thm:interpolated_segal_sugawara_Type_B_C} we have a set of free differential generators for the interpolated Feigin--Frenkel centers at the critical level over $R=\C[T]$: 
\begin{equation}
    \fr{z}(\Lulhat{so}{T})=R[\mathsf{T}][\phi_{2i,T}^B:i\in \Z_+], \qquad
    \fr{z}(\Lulhat{sp}{T})=R[\mathsf{T}][\phi_{2i,T}^C:i\in \Z_+]
\end{equation}
On the $\cW$-algebra side, by Proposition \ref{prop:differential_algebra_realization_SA_khesin_malikov} and Lemma \ref{lemma:duality_Poisson_algebras} we have two sets of free differential generators of the universal self-adjoint Khesin--Malikov PVA over $R$:
\begin{equation}
    \cW(\fr{po}_T)=R[\pr][u_{(-T+2i-1,T)}:i\in \Z_+]
= R[\pr][H_{(2i,T)}^{ev}:i\in \Z_+]\end{equation}
By Lemma \ref{lemma:FF_generators_TYPE_BC}, the finite-rank Feigin--Frenkel isomorphism maps between these free differential generators:
\begin{equation}
    \fr{f}_{\fr{so}_{2n+1}}^B:(\phi_{2,2n+1}^B,\cdots, \phi_{2n,2n+1}^B)\mapsto (H_{2}^C,\cdots, H_{2n}^C), \qquad
    \fr{f}_{\fr{sp}_{2n}}^C:(\phi_{2,2n}^C,\cdots, \phi_{2n,2n}^C)\mapsto (u_{-2n}^B,\cdots, u_{-2}^B)
\end{equation}
Using this, we can naturally extend the Feigin--Frenkel isomorphism at the critical level to complex rank in Types $B$ and $C$. 
\begin{theorem}[Interpolated Feigin--Frenkel Isomorphism at the Critical Level: Types $B$ and $C$]
\label{theorem:FF_inter_type_BC}

    There are isomorphisms of $\C[T]$-Poisson vertex algebras given on generators by:
    \begin{equation}
        \fr{f}_T^B:\fr{z}(\Lulhat{so}{T})\rightarrow \cW(\fr{po}_{T-1}), \qquad \fr{f}_T^B(\phi^B_{2i,T})=H_{(2i,T-1)}^{ev} \quad i\geq1
    \end{equation}
    \begin{equation}
        \fr{f}_T^C:\fr{z}(\Lulhat{sp}{T})\rightarrow \cW(\fr{po}_{T+1}), \qquad \fr{f}_T^C(\phi^C_{2i,T})=u_{(-T+2(i-1),T+1)} \quad i\geq 1
    \end{equation}
    Specializing to $\alpha\in \C$ induces isomorphisms of $\C$-Poisson vertex algebras:
    \begin{equation}
        \fr{f}_{T=\alpha}^B:\fr{z}(\Lulhat{so}{\alpha})\rightarrow \cW(\fr{po}_{T=\alpha-1}), \qquad \fr{f}_{T=\alpha}^B(\phi^B_{2i,T=\alpha})= H_{(2i,\alpha-1)}^{ev} \quad i\geq1
    \end{equation}
    \begin{equation}
        \fr{f}_{T=\alpha}^C:\fr{z}(\Lulhat{sp}{\alpha})\rightarrow \cW(\fr{po}_{T=\alpha+1}), \qquad \fr{f}_{T=\alpha}^C(\phi^C_{2i,T=\alpha})= u_{(-\alpha+2(i-1),\alpha+1)} \quad i\geq1 
    \end{equation}
These isomorphisms interpolate the usual Feigin--Frenkel isomorphisms in
the sense that the following diagrams commute:
\begin{figure}[H]
\centering
\noindent\hspace*{-1.4em}
\begin{minipage}{0.53\textwidth}
\centering
\[
\begin{tikzcd}[
	column sep=1.85em,
	row sep=large,
	cells={nodes={inner xsep=1pt}}
]
	{\fr{z}(\Lulhat{so}{2n+1})}
		&&
	{\cW(\fr{po}_{T=2n})}
	\\
	\\
	{\fr{z}(\widehat{\fr{so}}_{2n+1})}
		&&
	{\cW^{\mathrm{Adl}}(\fr{sp}_{2n})}
	\arrow["{\scriptstyle \fr{f}^B_{T=2n+1}}", from=1-1, to=1-3]
	\arrow["{\scriptstyle \ssV{2n+1}}"', from=1-1, to=3-1]
	\arrow["{\scriptstyle \ssW{2n}}", from=1-3, to=3-3]
	\arrow["{\scriptstyle \fr{f}^B_{\fr{so}_{2n+1}}}"', from=3-1, to=3-3]
\end{tikzcd}
\]
\end{minipage}
\hspace*{-2.2em}
\begin{minipage}{0.53\textwidth}
\centering
\[
\begin{tikzcd}[
	column sep=1.85em,
	row sep=large,
	cells={nodes={inner xsep=1pt}}
]
	{\fr{z}(\Lulhat{sp}{2n})}
		&&
	{\cW(\fr{po}_{T=2n+1})}
	\\
	\\
	{\fr{z}(\widehat{\fr{sp}}_{2n})}
		&&
	{\cW^{\mathrm{Adl}}(\fr{so}_{2n+1})}
	\arrow["{\scriptstyle \fr{f}^C_{T=2n}}", from=1-1, to=1-3]
	\arrow["{\scriptstyle \ssV{2n}}"', from=1-1, to=3-1]
	\arrow["{\scriptstyle \ssW{2n+1}}", from=1-3, to=3-3]
	\arrow["{\scriptstyle \fr{f}^C_{\fr{sp}_{2n}}}"', from=3-1, to=3-3]
\end{tikzcd}
\]
\end{minipage}
\caption{Finite-rank specialization of the interpolated Feigin--Frenkel isomorphisms in Types \(B\) and \(C\).}
\label{diagram:interpolated_FF_type_BC}
\end{figure}
\end{theorem}
\begin{proof}
      By the discussion preceding the theorem,  $\fr{f}_{T}^B,\fr{f}_T^C$ are isomorphisms of differential algebras over $\C[T]$. Similarly, after specializing $T=\alpha\in \C$, $\fr{f}_{T=\alpha}^B,\fr{f}_{T=\alpha}^C$ are isomorphisms of differential algebras over $\C$. Moreover, by Theorem \ref{thm:interpolated_segal_sugawara_Type_B_C} and Proposition \ref{prop:evaluation_functor_UKM_PVA}, evaluation is compatible with the $\C$-PVA structures:
\[
\mathrm{Ev}_{T=\alpha}(\fr{z}(\Lulhat{so}{T}))\cong \fr{z}(\Lulhat{so}{\alpha}),
\qquad
\mathrm{Ev}_{T=\alpha}(\fr{z}(\Lulhat{sp}{T}))\cong \fr{z}(\Lulhat{sp}{\alpha}) 
\]
\[
\mathrm{Ev}_{T=\alpha}(\cW(\fr{po}_T))\cong \cW(\fr{po}_{T=\alpha})
\]
Thus, once $\fr{f}_T^B$ and $\fr{f}^C_T$ are shown to preserve the \(\lambda\)-bracket over $\C[T]$, their specializations $\fr{f}_{T=\alpha}^B$ and $\fr{f}_{T=\alpha}^C$ are automatically isomorphisms of $\C$-PVAs under these identifications.

It remains to prove bracket preservation for $\fr{f}_T^B$ and $\fr{f}_T^C$ on their respective free differential generators. After accounting for the rank shift, the Type $C$ case is identical, mutatis mutandis, to the Type $A$ argument, and so is omitted. We prove the Type $B$ case. Fix $i,j\in \Z_+$. By the asymptotic component of Theorem~\ref{thm:interpolated_segal_sugawara_Type_B_C}, there exists $N_{(i,j)}\in \Z_{\geq 0}$ such that for odd $N'\geq N_{(i,j)}$ there is a  differential algebra homomorphism $\ssIV{N'}:\fr{z}(\widehat{\fr{so}}_{N'})\rightarrow \fr{z}(\Lulhat{so}{N'})$ for which the following holds:
    \begin{equation}
        \mathrm{Ev}_{T=N'}(\{\phi^B_{2i,T}{}_{\lambda}\phi^B_{2j,T}\})=\{\phi^B_{2i,T=N'}{}_{\lambda}\phi^B_{2j,T=N'}\}= \ssIV{N'}(\{\phi^B_{2i, N'}{}_{\lambda}\phi^B_{2j, N'}\})
    \end{equation} 
    By the asymptotic component of Proposition \ref{prop:interpolate_gelfand_dickey_projection_symplectic_orthogonal}, there exists $M_{(i,j)}\in \Z_{\geq 0}$ such that for even $N''\geq M_{(i,j)}$ there is a differential algebra homomorphism $\ssIW{N''}:\cW^{\mathrm{Adl}}(\fr{sp}_{N''})\rightarrow \cW(\fr{po}_{T=N''})$ for which the following holds:
   \begin{equation}
\begin{aligned}
\mathrm{Ev}_{T=N''}
\bigl(\{u_{(-T+2i-1,T)}{}_{\lambda}u_{(-T+2j-1,T)}\}\bigr)
&=
\{u_{(-N''+2i-1,N'')}{}_{\lambda}
u_{(-N''+2j-1,N'')}\} \\
&=
\ssIW{N''}
\bigl(\{u^C_{-N''+2i-1}{}_{\lambda}
u^C_{-N''+2j-1}\}\bigr).
\end{aligned}
\end{equation}
After increasing the bound if necessary, the same comparison holds for the $H^{ev}$-generators. Namely, there exists $L_{(i,j)}\in \Z_{\geq0}$ such that, for every even $M\geq L_{(i,j)}$:
\begin{equation}
       \mathrm{Ev}_{T=M}( \{H^{ev}_{(2i,T)}{}_{\lambda}H^{ev}_{(2j,T)}\})=
 \{H^{ev}_{(2i,M)}{}_{\lambda}H^{ev}_{(2j,M)}\}=\ssIW{M}(\{H^C_{2i}{}_{\lambda}H^C_{2j}\})
    \end{equation}
By Lemma \ref{lemma:duality_Poisson_algebras},  $H_{2i}^C$ and $H_{2j}^C$ can be expressed as differential polynomials in $u^C_{-M+1},u^C_{-M+3},\cdots, u^C_{-M+2i-1}$ and $u^C_{-M+1},u^C_{-M+3},\cdots, u^C_{-M+2j-1}$, respectively. Thus, one may take $L_{(i,j)}:=\max\{M_{(r,s)}: 1\leq r\leq i, 1\leq s\leq j\}$. 
Letting $N_0\geq \max\{N_{(i,j)}, L_{(i,j)}+1\}$, for odd $N\geq N_0$ we have: 
\begin{equation}
    \fr{f}^B_{T=N}\circ \ssIV{N}=\ssIW{N-1}\circ\fr{f}^B_{\fr{so}_N}
\end{equation}
The same argument as in Type $A$ now shows that $\fr{f}^B_T$ preserves the $\lambda$-bracket on the free differential generators.

For Type $B$, the last claim follows by computing $(\fr{f}^{B}_{\fr{so}_{2n+1}}\circ \ssV{2n+1})(\phi^{B}_{2i,T=2n+1})$ for $1\leq i\leq n$ through Theorem \ref{thm:interpolated_segal_sugawara_Type_B_C}, Lemma \ref{lemma:FF_generators_TYPE_BC}, and computing $(\ssW{2n}\circ \fr{f}^B_{T=2n+1})(\phi^B_{2i,T=2n+1})$ through Proposition \ref{prop:interpolate_gelfand_dickey_projection_symplectic_orthogonal}. The Type $C$ claim follows similarly. 
\end{proof}
 
\begin{theorem}
\label{theorem:FF_inter_type_BC_parity}
    The interpolated Feigin--Frenkel isomorphisms $\fr{f}_T^B, \fr{f}_T^C$ are related through the interpolated parity functor. More precisely, Diagram \ref{diagram:interpolated_FF_parity_functor_typeBC} commutes and continues to hold when $T$ is evaluated to $\alpha\in \C$. 
   \begin{figure}[H]
    \[\begin{tikzcd}
	{\mathfrak{z}(\Lulhat{sp}{T})} & {\cW(\fr{po}_{T+1})} \\
	{\mathfrak{z}(\Lulhat{so}{-T})} & {\cW(\fr{po}_{-T-1})}
	\arrow["{\mathfrak{f}_T^C}", from=1-1, to=1-2]
	\arrow["\Pi_{T+1}", from=1-2, to=2-2]
	\arrow["\Pi_{T}", from=1-1, to=2-1]
	\arrow["{\mathfrak{f}_{-T}^B}", from=2-1, to=2-2]
\end{tikzcd}\]
\caption{}
\label{diagram:interpolated_FF_parity_functor_typeBC}
\end{figure}
\end{theorem}
\begin{proof}
   The proof is identical to the Type $A$ argument, using Corollary~\ref{corollary:anti_isomorphism_generators}, Lemma~\ref{lemma:duality_Poisson_algebras}, Theorem~\ref{theorem:FF_inter_type_BC} and Lemma~\ref{lemma:FF_generators_TYPE_BC}. 
\end{proof}
\appendix
\include{FFDIC_Appendix_v3_submission}
\bibliographystyle{ieeetr}
\bibliography{FFDIC_v3_submission_ref}

\end{document}

%% file: FFDIC_Appendix_v3_submission.tex
\section{}
In this appendix we provide the proofs and formalism needed for the results in Section~\ref{section:vertex_algebras_universal_affine_interpolated}. More precisely, we collect the technical ingredients required for the categorical construction of the vacuum module. We first recall the necessary background on PROPs and operads. We then use this formalism to reduce the categorical Poincar\'e--Birkhoff--Witt theorem in the Ind-completion of a symmetric pseudo-tensor category to corresponding statements in suitable PROP categories. With the PBW theorem in place, we explain the equivalence between the representation theories of $\fr{g}$ and $U(\fr{g})$, and construct generalized Verma modules associated with Lie algebra triples.

Throughout this appendix, $\mathbb{F}$ denotes a field of characteristic $0$, and $R$ denotes a commutative $\mathbb{F}$-algebra of characteristic $0$.

\subsection{PROPs, Operads and Lie Formulas}

For the basics of operads and PROPs, we refer respectively to  \cite{Operads_Algebra_topology_physics}, \cite{PROP}. 

\begin{definition}
An operad in the category of $R$-modules is a collection of right $R[\Sigma_n]$-modules:
\[
\mathcal{P}=\{\mathcal{P}(n)\}_{n\ge 0},
\]
together with $R$-linear maps (the $\circ_i$-compositions):
\[
\circ_i:\mathcal{P}(m)\otimes \mathcal{P}(n)\longrightarrow \mathcal{P}(m+n-1),
\]
for $1\leq i\leq m$ and $n\geq 0$, satisfying the associativity, equivariance, and unit axioms. See \cite[Definition~1.16, \S 1.7]{Operads_Algebra_topology_physics}.
\end{definition}

\begin{definition}[$R$-linear PROP]
A $\mathrm{Mod}(R)$-enriched PROP is a symmetric strict monoidal category $P=(P,\odot,S,1)$ enriched over $\mathrm{Mod}(R)$ such that:
\begin{enumerate}
    \item the objects are identified with $\Z_{\geq 0}$;
    \item the monoidal product satisfies $m\odot n=m+n$ for all $m,n\in \Z_{\geq 0}$.
\end{enumerate}
For a PROP $P$ and $m,n\in \Z_{\geq 0}$, we write:
\begin{equation}
    \mathbf{P}(m,n):=\mathrm{Mor}_{\mathbf{P}}(m,n)
\end{equation}
Each $\mathbf{P}(m,n)$ is an $R$-module, and it carries the structure of an $(S_m,S_n)$-bimodule. By an \emph{$R$-linear PROP}, we mean a $\mathrm{Mod}(R)$-enriched PROP.
\end{definition}

PROPs can be constructed from operads.
\begin{example}[{\cite[Example~60]{PROP}}]
Let $\mathcal{P}$ be an operad. Then there exists a unique PROP $\mathbf{P}$ generated by $\mathcal{P}$ such that $\mathbf{P}(1,n)=\mathcal{P}(n)$. Its components are given by:
\begin{equation}
\mathbf{P}(m,n):=\bigoplus_{r_1+\cdots+r_m=n}\bigl[\mathcal{P}(1,r_1)\otimes\cdots\otimes \mathcal{P}(1,r_m)\bigr]\times_{S_{r_1}\times\cdots\times S_{r_m}} S_n
\label{eqn:PROP_operad_hom_description}
\end{equation}
One should think of this as concatenating together $m$ operations whose total number of inputs is $n$, and then inducing so as to allow all possible permutations of the inputs.
\end{example}
\begin{remark}
    We follow the convention that, for a PROP $\mathbf{P}$, the component $\mathbf{P}(m,n)$
consists of operations with $n$ inputs and $m$ outputs. 
\end{remark}
\begin{notation}
For notational convenience, we will occasionally write $\mathbf{1}_{\mathrm{P}}$ for the monoidal unit of an $R$-linear PROP $\mathrm{P}$, and $\bullet_{\mathrm{P}}$ for the generating object, i.e.\ the object identified with $1\in \Z_{\geq 0}$. If a PROP is defined from an operad $\mathcal{P}$, we denote the generating object simply as $\bullet_{\mathcal{P}}$.
\end{notation}

In particular, if we start with an operad $\mathcal{P}$ and then take the additive completion and Ind-completion of the associated PROP $\mathbf{P}$, we obtain:
\begin{align}
    \mathrm{End}_{\mathrm{Ind}(\mathbf{P}^{\mathrm{Add}})}(\mathcal{T}(\bullet_{\mathcal{P}}))
    &=\bigoplus_{n,m=0}^{\infty}\mathbf{P}(m,n) \label{eqn:tensor_algebra_decomp}\\
    &=\bigoplus_{n,m=0}^{\infty}\bigoplus_{r_1+\cdots+r_m=n}\bigl[\mathcal{P}(1,r_1)\otimes\cdots\otimes \mathcal{P}(1,r_m)\bigr]\times_{S_{r_1}\times\cdots\times S_{r_m}} S_n
\end{align}

The operad $\mathrm{Lie}$ is a particularly important example. Recall that $\mathrm{Lie}$ is generated by a single binary operation required to be skew-symmetric in the inputs and to satisfy the Jacobi identity. Similarly, one has the operad of unital associative algebras, denoted $u\mathrm{Ass}$.

The following observation is essential for applying the formalism of PROPs in our setting.
\begin{definition}
Let $\mathcal{P}$ be an operad enriched in $\mathrm{Mod}(R)$, and let $\mathrm{P}_{\mathcal{P}}$ be the associated PROP. The category of $\mathcal{P}$-algebras in a symmetric monoidal category $\cC$ is, by definition:
\begin{equation}
    \mathrm{Alg}_{\mathcal{P}}(\cC):=\mathrm{Fun}^{\otimes}_{R}(\mathrm{P}_{\mathcal{P}},\cC)
\end{equation}
where $\mathrm{Fun}^{\otimes}_R(\mathrm{P}_{\mathcal{P}},\cC)$ denotes the category of $R$-linear braided strong monoidal functors from $\mathrm{P}_{\mathcal{P}}$ to $\cC$.
\end{definition}

\begin{remark}
This reduces universal statements about algebraic objects internal to a symmetric monoidal category to statements in the corresponding PROP categories. For example, statements about operadic Lie algebras in a symmetric pseudo-tensor category reduce to statements about the PROP associated with $\mathrm{Lie}$. Similarly, unital associative algebras in a symmetric pseudo-tensor category arise from functors out of the PROP associated with $u\mathrm{Ass}$.
\end{remark}

\begin{notation}
Of particular interest are the categories of Lie algebras and unital associative algebras in a symmetric monoidal category $\cC$ associated with the operads $\mathrm{Lie}$ and $u\mathrm{Ass}$, respectively. To streamline the exposition, we use the shorthand:
\begin{equation}
    \mathrm{Lie}_{\cC}, \qquad \mathrm{uAss}_{\cC}
\end{equation}
\end{notation}

\begin{definition}
If $\cC$ is a symmetric $R$-linear monoidal category, then there is a natural functor:
\begin{equation}
    C:\mathrm{uAss}_{\cC}\rightarrow \mathrm{Lie}_{\cC}
\end{equation}
given by passing to the commutator bracket.

In the next subsection we construct a left adjoint to this functor under suitable hypotheses.
\end{definition}

\begin{example}
The PROP associated with $\mathrm{Lie}$ admits the following combinatorial description. Let $\mathrm{FreeLie}(x_1,\dots,x_n)$ denote the free Lie algebra on $n$ letters. Each Lie monomial has a multidegree $(m_1,\dots,m_n)$, where $m_i$ records how many times $x_i$ occurs. A Lie polynomial is called \emph{multilinear} if it has multidegree $(1,\dots,1)$.

One may identify $\mathrm{Lie}(n)$ with the free $R$-module spanned by multilinear Lie monomials in $n$ letters; see \cite[p.~58]{Operads_Algebra_topology_physics}. In particular, $\mathrm{Lie}(m,n)$ may be viewed as the space spanned by tensor products of $m$ multilinear Lie polynomials in $n$ letters. Note, in particular, that the $n$ letters may be permuted arbitrarily.

\begin{center}
    $\mathrm{Lie}(m,n)$ is the free $R$-module with basis given by $m$ multilinear Lie polynomials in $n$ variables.
\end{center}

By combining this observation with Equation~\eqref{eqn:tensor_algebra_decomp}, we find that the endomorphism algebra of the tensor algebra is given by:
\begin{center}
    $\mathrm{End}_{\mathrm{Ind}(\mathcal{P}_{\mathrm{Lie}}^{\mathrm{Add}})}(\mathcal{T}(\bullet_{\mathrm{Lie}}))$ is the free $R$-module with basis obtained by concatenating multilinear Lie polynomials.
\end{center}

Combinatorially, one may think of a Lie polynomial in $n$ variables as a binary rooted tree with $n$ leaves. Composition may then be viewed as stacking concatenations of such trees.
\end{example}

\begin{example}
Another important example is the PROP associated with the operad of unital associative algebras. It may be described as follows. Let $\mathrm{FreeAlg}(x_1,\dots,x_n)$ be the free unital associative algebra on $n$ letters over $R$. Each associative monomial has a multidegree $(m_1,\dots,m_n)$, where $m_i$ records how many times $x_i$ occurs. An associative monomial is called multilinear if it has multidegree $(1,\dots,1)$.

One may identify $\mathrm{Ass}(n)$ with the free $R$-module spanned by multilinear associative monomials in $n$ letters. In particular, $\mathrm{Ass}(m,n)$ may be viewed as the space spanned by tensor products of $m$ multilinear associative polynomials in $n$ letters. Again, the $n$ letters may be permuted arbitrarily. Consequently, the endomorphism algebra of the tensor algebra is given by:
\begin{center}
    $\mathrm{End}_{u\mathrm{Ass}}(\mathcal{T}(\bullet))$ is the free $R$-module with basis obtained by concatenating multilinear associative polynomials.
\end{center}
\end{example}
We will use the following translation lemma.
\begin{lemma}
Consider the operad $\mathrm{Lie}$ and the associated $R$-linear PROP $P_{\mathrm{Lie}}$. A pure Lie formula with $n$ inputs and $m$ outputs is a formula built only from sums, tensor products, the braiding, and Lie brackets on $n$ distinct variables, producing $m$ outputs. Let $\mathcal{F}_{n,m}$ be the $R$-module spanned by such formulas, modulo the multilinearity, antisymmetry, and Jacobi relations. If $k=m_1+\cdots+m_r$, then one may compose formulas $B_i\in \mathcal{F}_{n,m_i}$ for $1\leq i\leq r$ with $A\in \mathcal{F}_{k,t}$ to obtain:
\[
A(B_1,\dots,B_r)\in \mathcal{F}_{n,t}
\]
This defines a $\mathrm{Mod}(R)$-enriched category $\mathcal{F}$ with objects $\Z_{\geq 0}$ and morphism spaces $\mathcal{F}_{n,m}$. Furthermore, there is an equivalence of $\mathrm{Mod}(R)$-enriched categories:
\begin{equation}
    \mathrm{Inp}:\mathrm{P}_{\mathrm{Lie}}\rightarrow \mathcal{F}
\end{equation}

Recall that $\mathcal{F}^{\mathrm{Add}_R}$ denotes the $R$-linear additive completion. We refer to morphisms in this category as Lie formulas. Consequently, we have:
\begin{equation}
    \mathrm{P}_{\mathrm{Lie}}^{\mathrm{Add}_R}\cong \mathcal{F}^{\mathrm{Add}_R}
\end{equation}
\label{lemma:translation_lemma}
\end{lemma}

\begin{proof}
This follows from the fact that the Lie operad is the free operad modulo the Lie relations.
\end{proof}

\begin{corollary}
If a Lie formula holds on $\Z_{\geq 0}$ letters, then the corresponding equality holds in the Ind-completion of $\mathrm{P}_{\mathrm{Lie}}^{\mathrm{Add}_R}$.
\label{corollary:translation_cor}
\end{corollary}
\begin{proof}
This follows from the fact that Lie formulas may be interpreted as Lie polynomials on $\Z_{\geq 0}$ letters.
\end{proof}

Analogous statements hold for the operad $u\mathrm{Ass}$.

\begin{remark}
It is important not to confuse the monoidal structure in the PROP associated with $u\mathrm{Ass}$ with the multiplicative structure coming from the associative algebra itself.
\end{remark}

\begin{example}
An example of a pure Lie bracket formula is:
\begin{equation}
    [x_1,[x_2,x_3]]+[[x_2,x_1],x_3]
\end{equation}
An example of a pure unital associative algebra formula is:
\begin{equation}
    x_1*x_2+x_2*x_1
\end{equation}
where $*$ denotes the product in the free associative algebra.
\end{example}

Recall that if $\cC$ is an additive $R$-linear category, then $\cC^{\mathrm{Kar}_R}$ denotes its Karoubian completion.

\begin{notation}
For an $R$-linear PROP $P$, we use the abbreviation:
\begin{equation}
    U(P):=\mathrm{Ind}\bigl((P^{\mathrm{Add}_R})^{\mathrm{Kar}_R}\bigr)
\end{equation}
\end{notation}

The following is straightforward, so we omit the proof.
\begin{prop}
Let $\mathcal{P}$ be an operad, let $\cC$ be a symmetric pseudo-tensor category over $R$, and let $A$ be a $\mathcal{P}$-algebra in $\cC$. Denote the corresponding functor by:
\begin{equation}
    F:P_{\mathcal{P}}\rightarrow \cC, \qquad F(\bullet_{\mathcal{P}}):=A
\end{equation}
Then there is an induced functor:
\begin{equation}
    \tilde{F}:U(P_{\mathcal{P}})\rightarrow \mathrm{Ind}(\cC)
\end{equation}
\end{prop}
\subsection{Universal Enveloping Algebras}
\label{subsection:universal_enveloping_algebra}

In \cite[\S 1.3.7]{Quantum_Fields_Strings_Lie}, the authors explain how to define the universal enveloping algebra in the category of super vector spaces. Their construction starts from the symmetric algebra $S(\fr{g})$ of a Lie (super)algebra $\fr{g}$ and defines the multiplication inductively with respect to the natural filtration $F^{\bullet}S(\fr{g})$ by polynomial degree. The same approach works in the Ind-completion of a symmetric pseudo-tensor category in characteristic $0$. More precisely, the PROP formalism allows us to define universal enveloping algebras in this setting.

\begin{prop}
Let $\fr{g}$ denote the free Lie algebra on $\Z_{\geq 0}$ letters. The universal enveloping algebra structure
\[
U(\fr{g}):=(S(\fr{g}),m_{\star},1)
\]
defined in \cite[\S 1.3.7]{Quantum_Fields_Strings_Lie} induces a unital associative algebra structure on $S(\bullet_{\mathrm{Lie}})$ in $U(P_{\mathrm{Lie}})$, which we denote by $U(\bullet_{\mathrm{Lie}})$.
\label{prop:universal_enveloping_algebra_Construction_PROP}
\end{prop}
\begin{proof}
This follows immediately from Corollary~\ref{corollary:translation_cor}.
\end{proof}

In preparation for the universal property, we verify several basic identities.

\begin{lemma}
Let $\pi_{\bullet_{\mathrm{Lie}}}:\mathcal{T}(\bullet_{\mathrm{Lie}})\rightarrow U(\bullet_{\mathrm{Lie}})$ be the morphism in $U(P_{\mathrm{Lie}})$ obtained by repeatedly applying the multiplication $m_{\star}$ on $\bullet_{\mathrm{Lie}}\otimes \bullet_{\mathrm{Lie}}$. Let $\mathrm{sym}:\mathcal{T}(\bullet_{\mathrm{Lie}})\rightarrow S(\bullet_{\mathrm{Lie}})$ denote the direct sum of the symmetrizers $\mathrm{sym}_n:\bullet_{\mathrm{Lie}}^{\otimes n}\rightarrow S^n(\bullet_{\mathrm{Lie}})$. Then:
\begin{align}
    \pi_{\bullet_{\mathrm{Lie}}}\circ \mathrm{sym}&=\mathrm{sym}\\
    \mathrm{sym}\circ \pi_{\bullet_{\mathrm{Lie}}}&=\pi_{\bullet_{\mathrm{Lie}}}
\end{align}
\label{lemma:appendix_symmetrizers_pi_properties}
\end{lemma}

\begin{proof}
This is exactly the content of \cite[\S 1.3.7]{Quantum_Fields_Strings_Lie}, applied to the free Lie algebra on $\Z_{\geq 0}$ letters and transported via Corollary~\ref{corollary:translation_cor}.
\end{proof}

The following elementary lemma allows us to pass from an adjunction at the level of PROP categories to an adjunction between algebra objects internal to a category $\cC$.
\begin{lemma}
If $L:\cC\rightarrow \mathcal{E}$ and $R:\mathcal{E}\rightarrow \cC$ are functors such that $L\dashv R$, then precomposition defines functors
\[
R^*:\mathrm{Fun}(\cC,\mathcal{D})\rightarrow \mathrm{Fun}(\mathcal{E},\mathcal{D})
\qquad
L^*:\mathrm{Fun}(\mathcal{E},\mathcal{D})\rightarrow \mathrm{Fun}(\cC,\mathcal{D})
\]
such that:
\begin{equation}
    R^*\dashv L^*
\end{equation}
If the adjunction $L\dashv^{\otimes} R$ is monoidal, then $R^*\dashv L^*$ is also monoidal. 
\label{lemma:precomposition_adjunction}
\end{lemma}

By Lemma~\ref{lemma:precomposition_adjunction}, to prove the PBW theorem in the Ind-completion of a symmetric pseudo-tensor category it suffices to prove the corresponding adjunction in the completed PROP categories $U(\mathrm{P}_{\mathrm{Lie}})$ and $U(\mathrm{P}_{u\mathrm{Ass}})$.

\begin{lemma}
Denote the generating object of $U(\mathrm{P}_{\mathrm{Lie}})$ by $\bullet_{\mathrm{Lie}}$ and the generating object of $U(\mathrm{P}_{u\mathrm{Ass}})$ by $\bullet_{u\mathrm{Ass}}$. The commutator bracket induces a functor:
\begin{equation}
    \tilde{C}:U(P_{\mathrm{Lie}})\rightarrow U(P_{u\mathrm{Ass}})
\end{equation}
given by $\tilde{C}(\bullet_{\mathrm{Lie}}):=\bullet_{u\mathrm{Ass}}$, with the Lie bracket sent to the commutator bracket in $\mathrm{P}_{u\mathrm{Ass}}$. Similarly, the construction in Proposition~\ref{prop:universal_enveloping_algebra_Construction_PROP} induces a functor:
\begin{equation}
    \tilde{B}:U(P_{u\mathrm{Ass}})\rightarrow U(P_{\mathrm{Lie}})
\end{equation}
given by
\[
\tilde{B}(\bullet_{u\mathrm{Ass}}):=U(\bullet_{\mathrm{Lie}})=\bigl(S(\bullet_{\mathrm{Lie}}),m_{\star},1\bigr)
\]
These functors form an adjoint pair $\tilde{C}\dashv^{\otimes}\tilde{B}$.
\label{lemma:PROP_PBW_adjunction}
\end{lemma}

\begin{proof}
Since we are working with PROPs, it suffices to define the unit and counit on the monoidal generators and to verify the triangle identities there. Thus we seek morphisms:
\begin{equation}
    \epsilon_{\bullet_{u\mathrm{Ass}}}:(\tilde{C}\circ \tilde{B})(\bullet_{u\mathrm{Ass}})\rightarrow \bullet_{u\mathrm{Ass}}
\end{equation}
and
\begin{equation}
    \eta_{\bullet_{\mathrm{Lie}}}:\bullet_{\mathrm{Lie}}\rightarrow (\tilde{B}\circ \tilde{C})(\bullet_{\mathrm{Lie}})
\end{equation}
satisfying the zig-zag identities.

Observe that $(\tilde{B}\circ \tilde{C})(\bullet_{\mathrm{Lie}})$ is, as an object, simply $S(\bullet_{\mathrm{Lie}})$. We therefore take $\eta_{\bullet_{\mathrm{Lie}}}$ to be the canonical inclusion of $\bullet_{\mathrm{Lie}}$ into $S(\bullet_{\mathrm{Lie}})$. For $\epsilon_{\bullet_{u\mathrm{Ass}}}$, let
\[
\pi_{\bullet_{u\mathrm{Ass}}}:\mathcal{T}(\bullet_{u\mathrm{Ass}})\rightarrow \bullet_{u\mathrm{Ass}}
\]
be the canonical map obtained by repeatedly applying the multiplication on $\bullet_{u\mathrm{Ass}}$. We claim that:
\begin{equation}
    \epsilon_{\bullet_{u\mathrm{Ass}}}:=\pi_{\bullet_{u\mathrm{Ass}}}\circ \mathrm{sym}
\end{equation}
is the desired counit.

To check that $\epsilon_{\bullet_{u\mathrm{Ass}}}$ is a morphism of unital associative algebras, and hence a natural transformation, we proceed by induction on the standard filtration. Write:
\begin{equation}
    F^p(\bullet_{u\mathrm{Ass}}):=\bigoplus_{i=0}^p S^i(\bullet_{u\mathrm{Ass}})
\end{equation}
and for the filtered pieces set:
\begin{equation}
    m_{\star,(p,q)}:=m_{\star}\big|_{F^p\otimes F^q}:F^p(\bullet_{u\mathrm{Ass}})\otimes F^q(\bullet_{u\mathrm{Ass}})\rightarrow F^{p+q}(\bullet_{u\mathrm{Ass}})
\end{equation}
We must show that:
\begin{equation}
    \epsilon_{\bullet_{u\mathrm{Ass}}}\circ m_{\star,(p,q)}
    =
    m_{\bullet_{u\mathrm{Ass}}}\circ
    \bigl(\epsilon_{\bullet_{u\mathrm{Ass}}}\big|_{F^p(\bullet_{u\mathrm{Ass}})}\otimes
    \epsilon_{\bullet_{u\mathrm{Ass}}}\big|_{F^q(\bullet_{u\mathrm{Ass}})}\bigr)
\end{equation}

Since $m_{\star,(p,q)}$ lands in the symmetric algebra, it is enough to prove:
\begin{align}
    \pi_{\bullet_{u\mathrm{Ass}}}\circ m_{\star,(p,q)}
    &=
    m_{\bullet_{u\mathrm{Ass}}}\circ
    (\pi_{\bullet_{u\mathrm{Ass}}}\otimes \pi_{\bullet_{u\mathrm{Ass}}})\circ
    (\mathrm{sym}_p\otimes \mathrm{sym}_q)\\
    &=
    \pi_{\bullet_{u\mathrm{Ass}}}\circ (\mathrm{sym}_p\otimes \mathrm{sym}_q)
\end{align}

We prove this by applying Corollary~\ref{corollary:translation_cor} together with its associative counterpart, and then induct on $(p,q)$. The base cases $(0,0)$, $(1,0)$, and $(0,1)$ are immediate. For $(1,1)$, by definition:
\begin{equation}
    m_{\star,(1,1)}(x\otimes y):=\frac{1}{2}(x\otimes y+y\otimes x)+\frac{1}{2}[x,y]
\end{equation}
Equivalently:
\begin{align}
    m_{\star,(1,1)}(x\otimes y)
    &=\frac{1}{2}(x\otimes y+y\otimes x)+\frac{1}{2}(x*y-y*x)
\end{align}
Applying $\pi_{\bullet_{u\mathrm{Ass}}}$ therefore yields:
\begin{equation}
    (\pi_{\bullet_{u\mathrm{Ass}}}\circ m_{\star,(1,1)})(x\otimes y)=x*y
\end{equation}
which proves the case $(1,1)$.

Now assume the statement is known for $(1,q)$ whenever $q\leq n$. We claim it also holds for $(1,n+1)$.
By definition:
\begin{equation}
\label{eqn:inductive_step}
\begin{split}
    m_{\star,(1,k+1)} \bigl( x \otimes \operatorname{sym}(y_1 \otimes \dots \otimes y_{k+1}) \bigr)
    &=
    \operatorname{sym}(x \otimes y_1 \otimes \dots \otimes y_{k+1}) \\
    &\quad
    + \frac{1}{(k+2)!} \sum_{\sigma \in S_{k+1}} \sum_{i=1}^{k+1} (k-i+2)
    \bigl( y_{\sigma(1)} \star \dots \star [x, y_{\sigma(i)}] \star \dots \star y_{\sigma(k+1)} \bigr)
\end{split}
\end{equation}
Applying $\pi_{\bullet_{u\mathrm{Ass}}}$ to Equation~\eqref{eqn:inductive_step}, the first term becomes $\pi_{\bullet_{u\mathrm{Ass}}}\circ \mathrm{sym}_{k+2}$.

For the second summand, the induction hypothesis implies that for every $q<n+1$:
\begin{equation}
    \pi_{\bullet_{u\mathrm{Ass}}}\circ m_{\star,(1,q)}
    =
    \pi_{\bullet_{u\mathrm{Ass}}}\circ (\mathrm{sym}_1\otimes \mathrm{sym}_q)
\end{equation}
Moreover, the repeated multiplication map may be written as:
\begin{equation}
    \pi_{\star,(1,q-1)}:=
    m_{\star,(1,q-1)}\circ (\Id\otimes m_{\star,(1,q-2)})\circ \cdots \circ
    (\Id\otimes \cdots \otimes m_{\star,(1,1)})
\label{eqn:pi_star_formula}
\end{equation}
so that, after applying $\pi_{\bullet_{u\mathrm{Ass}}}$ and using the induction hypothesis, we obtain:
\begin{equation}
    \pi_{\bullet_{u\mathrm{Ass}}}\circ \pi_{\star,(1,q-1)}=\pi_{\bullet_{u\mathrm{Ass}}}.
\end{equation}

Thus the second term of Equation~\eqref{eqn:inductive_step} becomes:
\begin{align}
    \frac{1}{(k+2)!}\sum_{\sigma \in S_{k+1}}\sum_{i=1}^{k+1}(k+2-i)
    \bigl(y_{\sigma(1)}* \cdots * ([x,y_{\sigma(i)}])* \cdots * y_{\sigma(k+1)}\bigr)
    & \\
    =\frac{1}{(k+2)!}\sum_{\sigma \in S_{k+1}}\sum_{i=1}^{k+1}(k+2-i)
    \bigl(y_{\sigma(1)}* \cdots * (x* y_{\sigma(i)})* \cdots * y_{\sigma(k+1)}\bigr)
    & \\
    \quad -\frac{1}{(k+2)!}\sum_{\sigma \in S_{k+1}}\sum_{i=1}^{k+1}(k+2-i)
    \bigl(y_{\sigma(1)}* \cdots * (y_{\sigma(i)}* x)* \cdots * y_{\sigma(k+1)}\bigr). &
\end{align}
Rearranging the summations gives:
\begin{align}
\label{eqn:sum_0}
    \frac{1}{(k+2)!}\sum_{i=1}^{k+1}\sum_{\sigma \in S_{k+1}}(k+2-i)
    \bigl(y_{\sigma(1)}* \cdots * ([x,y_{\sigma(i)}])* \cdots * y_{\sigma(k+1)}\bigr)
    & \\
\label{eqn:sum_1}
    =\frac{1}{(k+2)!}\sum_{i=1}^{k+1}\sum_{\sigma \in S_{k+1}}(k+2-i)
    \bigl(y_{\sigma(1)}* \cdots * (x* y_{\sigma(i)})* \cdots * y_{\sigma(k+1)}\bigr)
    & \\
\label{eqn:sum_2}
    \quad -\frac{1}{(k+2)!}\sum_{i=1}^{k+1}\sum_{\sigma \in S_{k+1}}(k+2-i)
    \bigl(y_{\sigma(1)}* \cdots * (y_{\sigma(i)}* x)* \cdots * y_{\sigma(k+1)}\bigr) &
\end{align}

Consider the embedding $S_{k+1}\hookrightarrow S_{k+2}$ induced by the inclusion $\{2,\dots,k+2\}\hookrightarrow \{1,\dots,k+2\}$. Then:
\begin{equation}
    S_{k+2}/S_{k+1}=\bigcup_{i=1}^{k+2}(1\ i)S_{k+1}.
\end{equation}
Equation~\eqref{eqn:sum_1} is the sum over the first $k+1$ left cosets of the form $(1\ i)S_{k+1}$, while Equation~\eqref{eqn:sum_2} is the sum over the first $k+1$ cosets of the form $(1\ i+1)S_{k+1}$. Hence the two sums telescope, and for the overlapping terms the coefficient difference is:
\begin{equation}
    (k+2-i)-(k+2-(i-1))=-1.
\end{equation}
Therefore Equation~\eqref{eqn:sum_0} reduces to:
\begin{equation}
\label{eqn:result_inductive_coset_step}
\begin{split}
    \frac{k+1}{(k+2)!}\sum_{\sigma\in S_{k+1}} x* (y_{\sigma(1)}* \cdots * y_{\sigma(k+1)})
    &\\
    -\frac{1}{(k+2)!}\sum_{i=1}^{k}\sum_{\sigma \in S_{k+1}}
    \bigl(y_{\sigma(1)}* \cdots * (x* y_{\sigma(i)})* \cdots * y_{\sigma(k+1)}\bigr)
    &\\
    -\frac{1}{(k+2)!}\sum_{\sigma \in S_{k+1}}
    \bigl(y_{\sigma(1)}* \cdots * y_{\sigma(k+1)}* x\bigr) &
\end{split}
\end{equation}

Adding the first term on the right-hand side of Equation~\eqref{eqn:inductive_step} to Equation~\eqref{eqn:result_inductive_coset_step}, we find that the left-hand side of Equation~\eqref{eqn:inductive_step} equals:
\begin{equation}
    \frac{1}{(k+1)!}x* \Bigl(\sum_{\sigma \in S_{k+1}} y_{\sigma(1)}* \cdots * y_{\sigma(k+1)}\Bigr)
    =
    \pi_{\bullet_{u\mathrm{Ass}}}\circ (\mathrm{sym}_1\otimes \mathrm{sym}_{k+1})
\end{equation}
Hence the statement holds for all pairs $(1,q)$.

We now induct on $p$ to prove the general case. Recall from \cite[\S 1.3.7.3]{Quantum_Fields_Strings_Lie} that:
\begin{equation}
    m_{\star,(p,q)}\bigl(\mathrm{sym}(x_1\otimes\cdots \otimes x_p)\star Z\bigr)
    :=
    \frac{1}{p!}\sum_{\sigma\in S_p}x_{\sigma(1)}\star \cdots \star x_{\sigma(p)}\star Z
\end{equation}
The induction hypothesis allows us to split this into the $x$-terms and the $Z$-terms, and the conclusion follows from the identity $\pi\circ \mathrm{sym}=\mathrm{sym}$. Thus $\epsilon_{\bullet_{u\mathrm{Ass}}}$ is indeed a natural transformation.

It is clear that $\eta_{\bullet_{\mathrm{Lie}}}$ is a Lie algebra homomorphism, and hence a natural transformation.

It remains to verify the triangle identities:
\begin{align}
    1_{\tilde{C}(\bullet_{\mathrm{Lie}})}
    &=\epsilon_{\tilde{C}(\bullet_{\mathrm{Lie}})}\circ \tilde{C}(\eta_{\bullet_{\mathrm{Lie}}})
    \label{eqn:first_triangle_axiom}\\
    1_{\tilde{B}(\bullet_{u\mathrm{Ass}})}
    &= \tilde{B}(\epsilon_{\bullet_{u\mathrm{Ass}}})\circ \eta_{\tilde{B}(\bullet_{u\mathrm{Ass}})}
    \label{eqn:second_triangle_axiom}
\end{align}
Equation~\eqref{eqn:first_triangle_axiom} is immediate from the construction.

For Equation~\eqref{eqn:second_triangle_axiom}, note first that:
\begin{equation}
    \tilde{B}\bigl(\pi_{\bullet_{u\mathrm{Ass}}}\circ \mathrm{sym}_{\bullet_{u\mathrm{Ass}}}\bigr)
    =
    \pi_{\star}\circ \mathrm{sym}_{\tilde{B}(\bullet_{u\mathrm{Ass}})}
\end{equation}
where $\pi_{\star}$ denotes repeated multiplication in $\tilde{B}(\bullet_{u\mathrm{Ass}})$:
\[
\pi_{\star}: \mathcal{T}(\tilde{B}(\bullet_{u\mathrm{Ass}}))\rightarrow \tilde{B}(\bullet_{u\mathrm{Ass}})
\]
Thus it suffices to prove:
\begin{equation}
    \pi_{\star}\circ \mathrm{sym}_{\tilde{B}(\bullet_{u\mathrm{Ass}})}\circ \eta_{\tilde{B}(\bullet_{u\mathrm{Ass}})}
    =
    \mathrm{Id}_{\tilde{B}(\bullet_{u\mathrm{Ass}})}
\end{equation}
and this is immediate.
\end{proof}

\begin{theorem}
Let $\cC$ be a symmetric pseudo-tensor category over $R$ of characteristic $0$. Denote by $\mathrm{Lie}_{\IC}$ the category of Lie algebras internal to $\IC$, and similarly by $\mathrm{uAss}_{\IC}$ the category of unital associative algebras internal to $\IC$. Then the commutator functor:
\[
C:\mathrm{uAss}_{\IC}\rightarrow \mathrm{Lie}_{\IC}
\]
has a left adjoint:
\[
B:\mathrm{Lie}_{\IC}\rightarrow \mathrm{uAss}_{\IC}
\]
given by $B(\fr{g}):=U(\fr{g})$.
\label{theorem:categorical_PBW_appendix}
\end{theorem}

\begin{proof}
This follows from Lemmas~\ref{lemma:precomposition_adjunction} and \ref{lemma:PROP_PBW_adjunction}.
\end{proof}

\begin{corollary}
Let $\fr{g}$ be a Lie algebra in $\mathrm{Ind}(\cC)$ and let $A$ be a unital associative algebra in $\mathrm{Ind}(\cC)$. Suppose there is a morphism of Lie algebras:
\[
f:\fr{g}\rightarrow A^{\mathrm{com}}
\]
where $A^{\mathrm{com}}$ denotes the commutator Lie algebra. Then there exists a unique morphism of unital associative algebras:
\[
\tilde{f}:U(\fr{g})\rightarrow A
\]
such that, for the canonical map $\iota:\fr{g}\rightarrow U(\fr{g})$, one has $\tilde{f}\circ \iota=f$.
\end{corollary}

\begin{proof}
This is simply the adjunction restated as the universal property of $U(\fr{g})$.
\end{proof}

The following theorem allows us to pass between representations of $\fr{g}$ and representations of $U(\fr{g})$.

\begin{theorem}
Let $\mathrm{Mod}_{\mathrm{Ind}(\cC)}(\fr{g})$ denote the category of left $\fr{g}$-modules in the Ind-completion of a symmetric pseudo-tensor category $\cC$, and similarly let $\mathrm{Mod}_{\mathrm{Ind}(\cC)}(U(\fr{g}))$ denote the category of left $U(\fr{g})$-modules. These categories are equivalent as $R$-linear additive categories.
\label{theorem:representation_universal_enveloping_algebra}
\end{theorem}

\begin{proof}
The construction is essentially identical to the adjunction argument above, so we only sketch it. Given a Lie algebra module $(M,\rho_M)$, define
\[
\pi_{(M,\rho_M)}:\mathcal{T}(\fr{g})\otimes M\rightarrow M
\]
by repeatedly applying $\rho_M$. Then one checks, by the same inductive argument as above, that $\pi_{(M,\rho_M)}\circ (\mathrm{sym}\otimes \Id)$ equips $M$ with a $U(\fr{g})$-module structure. The converse direction is immediate.
\end{proof}

\subsection{Lie Algebra Triple}

Throughout this subsection, let $\cC$ denote a symmetric pseudo-tensor category over $R$ of characteristic $0$.

Recall that, by Definition~\ref{defn:S_graded_object}, a $\Z_{\geq 0}$-graded object $(X,F)$ satisfies
\begin{equation}F(T)\cong \bigoplus_{t\in T}F(\{t\})
\end{equation} canonically for every finite $T\subset \Z_{\geq 0}$. With this in mind, we define the associated graded as the following.
\begin{definition}
Let $X_F:=\mathrm{colim}(F)$ be a $\Z_{\geq 0}$-graded object arising from a functor $F:I_{\Z_{\geq 0}}\rightarrow \cC$. The associated graded of $X_F$ is defined by:
\begin{equation}
    \mathrm{gr}_F(X_F):=\bigoplus_{i=0}^{\infty}\frac{F(\{0,\dots,i\})}{F(\{0,\dots,i-1\})},
\end{equation}
with the convention that $F(\{0,\cdots,-1\})=0_{\cC}$.
\end{definition}

\begin{notation}
Since each $F(\{0,\dots,i\})$ is a direct summand of $F(\{0,\dots,i+1\})$, the quotient is well defined in a Karoubian category. For brevity, we write the induced filtration on $X$ as:
\begin{equation}
    F^iX:=F(\{0,\dots,i\}), \qquad i\in \Z_{\geq 0}
\end{equation}
\end{notation}

\begin{lemma}
Let $f:(X,F)\rightarrow (Y,G)$ be a morphism of $\Z_{\geq 0}$-graded objects that preserves the filtration $f|_{F^iX}:F^iX\rightarrow G^iY$. Then there is an induced morphism:
\begin{equation}
    \tilde{f}:\mathrm{gr}_F(X)\rightarrow \mathrm{gr}_G(Y)
\end{equation}
Furthermore, if $\tilde{f}$ is an isomorphism, then $f$ is an isomorphism.
\label{lemma:associated_graded}
\end{lemma}

\begin{proof}
Since $f$ preserves the filtration, there are maps $f_n:F^nX\rightarrow G^nY$ for all $n\geq 0$. For each $n$ we refine $F^nX$ and $G^nY$ as:
\begin{equation}
    F^nX=\bigoplus_{i=0}^n F(\{i\}),
    \qquad
    G^nY=\bigoplus_{i=0}^n G(\{i\})
\end{equation}
Hence for every $n\in \Z_{\geq 0}$ and every $k\leq n$ there are maps:
\begin{equation}
    f_{n,k}:F(\{n\})\rightarrow G(\{k\})
\end{equation}
such that:
\begin{equation}
    f_n=\bigoplus_{r,k=0}^n f_{r,k}
\end{equation}
Equivalently, one may view $f$ as a $\Z_{\geq 0}\times \Z_{\geq 0}$ matrix $(f_{r,k})_{r,k=0}^{\infty}$ with entries in $\mathrm{Hom}$-spaces. This matrix is upper triangular, and its diagonal is precisely $\tilde{f}=\sum_{n=0}^{\infty}f_{n,n}$. If $\tilde{f}$ is invertible, then each finite upper-triangular truncation $f_n$ is invertible. The inverses $f_n^{-1}$ are compatible and so $f$ is an isomorphism. 
\end{proof}

\begin{prop}
If $(X,F)$ and $(Y,G)$ are $\Z_{\geq 0}$-graded objects in $\IC$, then $X\otimes Y$ is again $\Z_{\geq 0}$-graded with filtration:
\begin{equation}
    (F\otimes G)^{i}(X\otimes Y)
    =
    \bigoplus_{k+s\leq i}
    \left(\frac{F^k(X)}{F^{k-1}(X)}\right)\otimes
    \left(\frac{G^s(Y)}{G^{s-1}(Y)}\right)
\end{equation}
and:
\begin{equation}
    \mathrm{gr}_F(X)\otimes \mathrm{gr}_G(Y)\cong \mathrm{gr}_{F\otimes G}(X\otimes Y)
\end{equation}
\end{prop}

\begin{proof}
This is essentially tautological from the definition of the tensor-product filtration.
\end{proof}

\begin{example}
Let $\fr{g}$ be a Lie algebra in $\IC$, and consider $U(\fr{g}):=(S(\fr{g}),m_{\star},1)$. This is a $\Z_{\geq 0}$-graded object whose algebra structure is compatible with the filtration:
\begin{equation}
    F^iU(\fr{g}):=\bigoplus_{k=0}^i S^k(\fr{g})
\end{equation}
By the construction of $m_{\star}$, one has:
\begin{equation}
    m_{\star}:F^iU(\fr{g})\otimes F^jU(\fr{g})\rightarrow F^{i+j}U(\fr{g})
\end{equation}
Furthermore, the defining formula for $m_{\star}$ implies that $\mathrm{gr}_F U(\fr{g})\cong S(\fr{g})$ as unital associative algebras.
\end{example}

\begin{definition}[Lie Algebra Triple]
A Lie algebra triple in a symmetric pseudo-tensor category $\cC$ over $R$ of characteristic $0$ is a tuple $(\fr{g},\fr{n}_-,\fr{p})$ such that each entry is a Lie algebra in $\IC$, and there is a decomposition of $\IC$-objects:
\begin{equation}
    \fr{g}=\fr{n}_-\oplus \fr{p}
\end{equation}
such that each summand is a Lie subalgebra.
\end{definition}
Notice a Lie algebra triple does not mean that $\fr{g}$ is isomorphic as a Lie algebra to the direct sum of two Lie algebras. Interaction between the summands is allowed. 

\begin{prop}
Let $(\fr{g},\fr{n}_-, \fr{p})$ be a Lie algebra triple in $\cC$. Then $U(\fr{n}_-)$ and $U(\fr{p})$ are direct summands of $U(\fr{g})$ as objects, and the multiplication map:
\begin{equation}
    m_{\fr{n}_-,\fr{p}}:U(\fr{n}_-)\otimes U(\fr{p})\rightarrow U(\fr{g})
\end{equation}
is an isomorphism of right $U(\fr{p})$-modules.
\label{prop:multiplication_iso_Lie_algebra_Triple}
\end{prop}

\begin{proof}
Let $\iota_{\fr{n}_-}:\fr{n}_-\rightarrow \fr{g}$ and $\iota_{\fr{p}}:\fr{p}\rightarrow \fr{g}$ denote the inclusion maps. These are morphisms of Lie algebras, so there are induced morphisms of associative algebras:
\[
U(\iota_{\fr{n}_-}):U(\fr{n}_-)\rightarrow U(\fr{g} )
\qquad
U(\iota_{\fr{p}}):U(\fr{p})\rightarrow U(\fr{g})
\]
The map under consideration is:
\begin{equation}
    m_{\fr{n}_-,\fr{p}}
    :=
    m_{\star}^{\fr{g}}\circ
    \bigl(U(\iota_{\fr{n}_-})\otimes U(\iota_{\fr{p}})\bigr)
\end{equation}
It is evidently a morphism of right $U(\fr{p})$-modules. By Lemma~\ref{lemma:associated_graded}, it is enough to show that its associated graded is an isomorphism. But:
\begin{equation}
    \mathrm{gr}(m_{\fr{n}_-,\fr{p}})
    =
    m_{S(\fr{g})}\circ
    \bigl(S(\iota_{\fr{n}_-})\otimes S(\iota_{\fr{p}})\bigr)
\end{equation}
Therefore it suffices to show that the multiplication map:
\begin{equation}
    S(\fr{n}_-)\otimes S(\fr{p})\rightarrow S(\fr{g})=S(\fr{n}_-\oplus \fr{p})
\end{equation}
is an isomorphism. This follows immediately from the decomposition:
\begin{equation}
    S^n(\fr{n}_-\oplus \fr{p})
    \cong
    \bigoplus_{k=0}^n S^k(\fr{n}_-)\otimes S^{n-k}(\fr{p})
\end{equation}
\end{proof}

\begin{prop}[Generalized Verma Modules]
\label{prop:generalized_verma_module_appendix}
Let $\fr{g}$ be a Lie algebra in $\IC$. For every Lie algebra triple $(\fr{g},\fr{n}_-,\fr{p})$, there is a natural restriction functor:
\begin{equation}
    \mathrm{Res}^{\fr{g}}_{\fr{p}}:\mathrm{Mod}_{\IC}(\fr{g})\rightarrow \mathrm{Mod}_{\IC}(\fr{p}).
\end{equation}
This restriction functor has a left adjoint:
\begin{equation}
    \mathrm{Ind}^{\fr{g}}_{\fr{p}}:\mathrm{Mod}_{\IC}(\fr{p})\rightarrow \mathrm{Mod}_{\IC}(\fr{g}),
\end{equation}
such that if $F:\cC\rightarrow \cD$ is a braided pseudo-tensor functor, then:
\begin{equation}
    F\circ \mathrm{Ind}_{\fr{p}}^{\fr{g}}
    =
    \mathrm{Ind}^{F(\fr{g})}_{F(\fr{p})}\circ F.
\end{equation}
Furthermore, these induction functors commute with base change along commutative algebras of characteristic $0$.
\end{prop}

\begin{proof}
Recall that $m_{\fr{g}}$ denotes the multiplication on $U(\fr{g})$, and let
\[
m^{-1}_{\fr{n}_-,\fr{p}}:U(\fr{g})\rightarrow U(\fr{n}_-)\otimes U(\fr{p})
\]
be the inverse of the right $U(\fr{p})$-module isomorphism from Proposition~\ref{prop:multiplication_iso_Lie_algebra_Triple}. For $X\in \mathrm{Mod}_{\IC}(\fr{p})$, define:
\begin{align}
    \mathrm{Ind}_{\fr{p}}^{\fr{g}}(X)
    &:=
    \bigl(U(\fr{n}_-)\otimes X,\mathrm{Ind}(\rho_X)\bigr)\\
    \mathrm{Ind}_{\fr{p}}^{\fr{g}}(f)
    &:=
    \Id_{U(\fr{n}_-)}\ox f
\end{align}
where the $U(\fr{g})$-module structure morphism is:
\begin{align}
    \mathrm{Ind}(\rho_X):U(\fr{g})\ox U(\fr{n}_-)\ox X&\rightarrow U(\fr{n}_-)\ox X\\
    \mathrm{Ind}(\rho_X)&:=
    (\Id_{U(\fr{n}_-)}\otimes \rho_X)\circ
    \bigl((m^{-1}_{\fr{n}_-,\fr{p}}\circ m_{\fr{g},\fr{n}_-})\ox \Id_X\bigr)
\end{align}
We first check that $\mathrm{Ind}_{\fr{p}}^{\fr{g}}(X)$ is indeed a left $U(\fr{g})$-module. Compute: 
\begin{align*}
    \mathrm{Ind}(\rho_X)\circ (\mathrm{Id}_{U(\fr{g})}\otimes \mathrm{Ind}(\rho_X))= \\ 
    \left((\Id_{U(\fr{n}_-)}\otimes \rho_X)\circ((m^{-1}_{\fr{n}_-,\fr{p}}\circ m_{\fr{g},\fr{n}_-})\ox \Id_{U(\fr{p})\ox X})\right)\circ \left( \mathrm{Id}_{U(\fr{g}) }\otimes ((\Id_{U(\fr{n}_-)}\otimes \rho_X)\circ((m^{-1}_{\fr{n}_-,\fr{p}}\circ m_{\fr{g},\fr{n}_-})\ox \Id_X))\right)=\\
    (\Id_{U(\fr{n}_-)}\otimes\rho_X\circ (\mathrm{Id}_{U(\fr{n}_-)}\otimes \Id_{U(\fr{p})}\otimes \rho_X))\circ((m^{-1}_{\fr{n}_-,\fr{p}}\circ m_{\fr{g},\fr{n}_-})\ox \Id_{U(\fr{p})\ox X}))\circ \left( \mathrm{Id}_{U(\fr{g}) }\otimes (((m^{-1}_{\fr{n}_-,\fr{p}}\circ m_{\fr{g},\fr{n}_-})\ox \Id_X))\right) =\\
   ( \Id_{U(\fr{n}_-)}\otimes \rho_X)\circ (\Id_{U(\fr{n}_-)}\otimes m_{\fr{p},\fr{p}}\otimes \Id_X)\circ((m^{-1}_{\fr{n}_-,\fr{p}}\circ m_{\fr{g},\fr{n}_-})\ox \Id_{U(\fr{p})\ox X}))\circ \left( \mathrm{Id}_{U(\fr{g}) }\otimes (((m^{-1}_{\fr{n}_-,\fr{p}}\circ m_{\fr{g},\fr{n}_-})\ox \Id_X))\right)\\
\end{align*}
Here we have used naturality of the tensor product and the fact that $\rho_X$ is a $U(\fr{p})$-module structure. Since $m_{\fr{n}_-,\fr{p}}$ is a right $U(\fr{p})$-module homomorphism, we have:
\begin{equation}
    (\Id_{U(\fr{n}_-)}\ox m_{\fr{p},\fr{p}})\circ
    (m^{-1}_{\fr{n}_-,\fr{p}}\otimes \Id_{U(\fr{p})})
    =
    m_{\fr{n}_-,\fr{p}}^{-1}\circ m_{\fr{g},\fr{p}}
\end{equation}
Therefore the previous expression simplifies to:
\begin{align}
    (\Id_{U(\fr{n}_-)}\otimes \rho_X)\circ
    (m_{\fr{n}_-,\fr{p}}^{-1}\circ m_{\fr{g},\fr{p}}\otimes \Id_X)\circ
    (m_{\fr{g},\fr{n}_-}\ox \Id_{U(\fr{p})\ox X})\circ
    \left(\mathrm{Id}_{U(\fr{g})}\otimes
    \bigl((m^{-1}_{\fr{n}_-,\fr{p}}\circ m_{\fr{g},\fr{n}_-})\ox \Id_X\bigr)\right)
\end{align}

On the other hand, associativity of multiplication implies:
\begin{equation}
    m_{\fr{g},\fr{p}}\circ (m_{\fr{g},\fr{n}_-}\otimes \Id_{U(\fr{p})})
    =
    m_{\fr{g},\fr{g}}\circ (\Id_{U(\fr{g})}\otimes m_{\fr{n}_-,\fr{p}})
\end{equation}
and
\begin{equation}
    m_{\fr{g},\fr{g}}\circ (\Id\otimes m_{\fr{g},\fr{n}_-})
    =
    m_{\fr{g},\fr{n}_-}\circ (m_{\fr{g},\fr{g}}\otimes \Id_{U(\fr{n}_-)})
\end{equation}
Using these identities, we obtain:
\begin{align}
    (\Id_{U(\fr{n}_-)}\otimes \rho_X)\circ
    (m_{\fr{n}_-,\fr{p}}^{-1}\circ m_{\fr{g},\fr{n}_-}\otimes \Id_X)\circ
    (m_{\fr{g},\fr{g}}\otimes \Id_{U(\fr{n}_-)}\ox \Id_X)
    =
    \mathrm{Ind}(\rho_X)\circ
    (m_{\fr{g},\fr{g}}\otimes \Id_{\mathrm{Ind}^{\fr{g}}_{\fr{p}}(X)})
\end{align}
The unit axiom is immediate, so $\mathrm{Ind}_{\fr{p}}^{\fr{g}}(X)$ is indeed a left $U(\fr{g})$-module. The verification that $\mathrm{Ind}^{\fr{g}}_{\fr{p}}(f)$ is a module morphism is similar and is omitted.

To prove the adjunction, define:
\[
\Phi_{X,Y}:\mathrm{Hom}_{\mathrm{Mod}_{\fr{g}}(\IC)}(\mathrm{Ind}^{\fr{g}}_{\fr{p}}(X),Y)\rightarrow \mathrm{Hom}_{\mathrm{Mod}_{\fr{p}}(\IC)}(X,\mathrm{Res}_{\fr{p}}^{\fr{g}}(Y))
\]
by:
\begin{equation}
    \Phi_{X,Y}(f):=f\circ (1_{\fr{n}_-}\otimes \Id_X):X\rightarrow Y
\end{equation}
where $1_{\fr{n}_-}$ denotes the unit of $U(\fr{n}_-)$. Its inverse is:
\[
\Phi^{-1}_{X,Y}:\mathrm{Hom}_{\mathrm{Mod}_{\fr{p}}(\IC)}(X,\mathrm{Res}_{\fr{p}}^{\fr{g}}(Y))
\rightarrow
\mathrm{Hom}_{\mathrm{Mod}_{\fr{g}}(\IC)}(\mathrm{Ind}^{\fr{g}}_{\fr{p}}(X),Y)
\]
given by:
\begin{equation}
    \Phi^{-1}_{X,Y}(g):=\rho_Y\circ (\Id_{\fr{n}_-}\otimes g)
\end{equation}
Indeed:
\begin{align}
    (\Phi_{X,Y}\circ \Phi^{-1}_{X,Y})(g)
    &=
    \rho_Y\circ (\Id_{\fr{n}_-}\otimes g)\circ (1_{\fr{n}_-}\otimes \Id_X)
    =
    \rho_Y\circ (1_{\fr{n}_-}\otimes \Id_Y)\circ g
    =
    g\\
    (\Phi^{-1}_{X,Y}\circ \Phi_{X,Y})(f)
    &=
    \rho_Y\circ (\Id_{\fr{n}_-}\otimes f)\circ
    (\Id_{\fr{n}_-}\otimes 1_{\fr{n}_-}\otimes \Id_X)
\end{align}

To see that the second composite equals $f$, note that:
\begin{align}
    &(\mathrm{Ind}(\rho_X)\circ (U(\iota_{\fr{n}_-})\otimes \Id_{U(\fr{n}_-)\otimes X}))
    \circ (\Id_{\fr{n}_-}\otimes 1_{\fr{n}_-}\otimes \Id_X) \\
    &\qquad =
    (\Id_{U(\fr{n}_-)}\otimes \rho_X)\circ
    \bigl((m^{-1}_{\fr{n}_-,\fr{p}}\circ m_{\fr{g},\fr{n}_-})\ox \Id_X\bigr)\circ
    (U(\iota_{\fr{n}_-})\otimes \Id_{U(\fr{n}_-)\otimes X})\circ
    (\Id_{\fr{n}_-}\otimes 1_{\fr{n}_-}\otimes \Id_X)
\end{align}
Moreover:
\begin{equation}
    m_{\fr{g},\fr{n}_-}\circ \bigl(U(\iota_{\fr{n}_-})\otimes 1_{\fr{n}_-}\bigr)=U(\iota_{\fr{n}_-})
\end{equation}
and since:
\begin{equation}
    m_{\fr{n}_-,\fr{p}}\circ (\Id_{U(\fr{n}_-)}\otimes 1_{\fr{p}})
    =
    U(\iota_{\fr{n}_-})
\end{equation}
we get:
\begin{equation}
    m^{-1}_{\fr{n}_-,\fr{p}}\circ U(\iota_{\fr{n}_-})
    =
    \Id_{U(\fr{n}_-)}\otimes 1_{\fr{p}}
\end{equation}
Using the unit axiom $\rho_X\circ (1_{\fr{p}}\otimes \Id_X)=\Id_X$, it follows that:
\begin{equation}
    (\Phi^{-1}_{X,Y}\circ \Phi_{X,Y})(f)=f
\end{equation}
Hence $\Phi_{X,Y}$ is an isomorphism, and one checks directly that it is natural in both variables.

Finally, induction is preserved by symmetric pseudo-tensor functors as the construction depends only on universal enveloping algebras, tensor products, and the braiding. The same argument shows that induction commutes with base change along commutative algebras of characteristic $0$.
\end{proof}